\documentclass[opre,nonblindrev]{informs3-VVM-noheader}

\SingleSpacedXI
\usepackage{natbib}
 \bibpunct[, ]{(}{)}{,}{a}{}{,}%
 \def\bibfont{\small}%
 \usepackage{url}
\usepackage{amsmath}
\usepackage{amssymb}
\usepackage{color}

\usepackage{algorithm}
\usepackage{algorithmic}

\usepackage{setspace}
\usepackage{ctable}
\usepackage[export]{adjustbox}

\usepackage{natbib}
 \bibpunct[, ]{(}{)}{,}{a}{}{,}%

\def\zb{\mathbf{z}}
\def\xb{\mathbf{x}}
\def\Xb{\mathbf{X}}

\def\cb{\mathbf{c}}

\def\Xcal{\mathcal{X}}

\def\Ncal{\mathcal{N}}
\def\Ccal{\mathcal{C}}
\def\Acal{\mathcal{A}}

\def\Pcal{\mathcal{P}}
\def\Tcal{\mathcal{T}}
\def\Dcal{\mathcal{D}}
\def\Ncal{\mathcal{N}}
\def\Scal{\mathcal{S}}
\def\Mcal{\mathcal{M}}

\def\Halmos{$\square$}
\def\Ibb{\mathbb{I}}

\def\var{\mathbf{V}}
\def\cond{\mathbf{C}}
\def\leftleaves{\mathbf{left}}
\def\rightleaves{\mathbf{right}}
\def\leftsplits{\mathbf{LS}}
\def\rightsplits{\mathbf{RS}}

\def\leaves{\mathbf{leaves}}
\def\splits{\mathbf{splits}}
\def\root{\mathbf{root}}
\def\leftchild{\mathbf{leftchild}}
\def\rightchild{\mathbf{rightchild}}
\def\untestedVars{\mathbf{untestedVariables}}

\def\avg{\mathrm{avg}}

\def\xb{\mathbf{x}}
\def\yb{\mathbf{y}}
\def\zerob{\mathbf{0}}

\def\getLeaf{\textsc{GetLeaf}}

\def\alphab{\boldsymbol \alpha}
\def\betab{\boldsymbol \beta}
\def\thetab{\boldsymbol \theta}

\def\Normal{N}
\def\IW{IW}
\def\pb{\mathbf{p}}
\def\fb{\mathbf{f}}
\def\db{\mathbf{d}}
\def\vecop{\mathrm{vec}}

\def\ie{\emph{i.e.}}

 \newcommand{\problemeqref}[1]{\eqref{#1}}

\TheoremsNumberedThrough     %
\ECRepeatTheorems

\EquationsNumberedThrough    %
\begin{document}
\RUNTITLE{Optimization of Tree Ensembles}
\RUNAUTHOR{Mi\v{s}i\'{c}}

\TITLE{Optimization of Tree Ensembles}

\ARTICLEAUTHORS{%
\AUTHOR{Velibor V. Mi\v{s}i\'{c}}
\AFF{Anderson School of Management, University of California, Los Angeles, 110 Westwood Plaza, Los Angeles, CA, 90095, \EMAIL{velibor.misic@anderson.ucla.edu}} %
} %

\ABSTRACT{%
Tree ensemble models such as random forests and boosted trees are among the most widely used and practically successful predictive models in applied machine learning and business analytics. Although such models have been used to make predictions based on exogenous, uncontrollable independent variables, they are increasingly being used to make predictions where the independent variables are controllable and are also decision variables. In this paper, we study the problem of tree ensemble optimization: given a tree ensemble that predicts some dependent variable using controllable independent variables, how should we set these variables so as to maximize the predicted value? We formulate the problem as a mixed-integer optimization problem. We theoretically examine the strength of our formulation, provide a hierarchy of approximate formulations with bounds on approximation quality and exploit the structure of the problem to develop two large-scale solution methods, one based on Benders decomposition and one based on iteratively generating tree split constraints. We test our methodology on real data sets, including two case studies in drug design and customized pricing, and show that our methodology can efficiently solve large-scale instances to near or full optimality, and outperforms solutions obtained by heuristic approaches. %
}%

\KEYWORDS{tree ensembles; random forests; mixed-integer optimization; drug design; customized pricing.} %

\maketitle

\section{Introduction}

\begin{figure}
\centering
\includegraphics[width=0.7\textwidth]{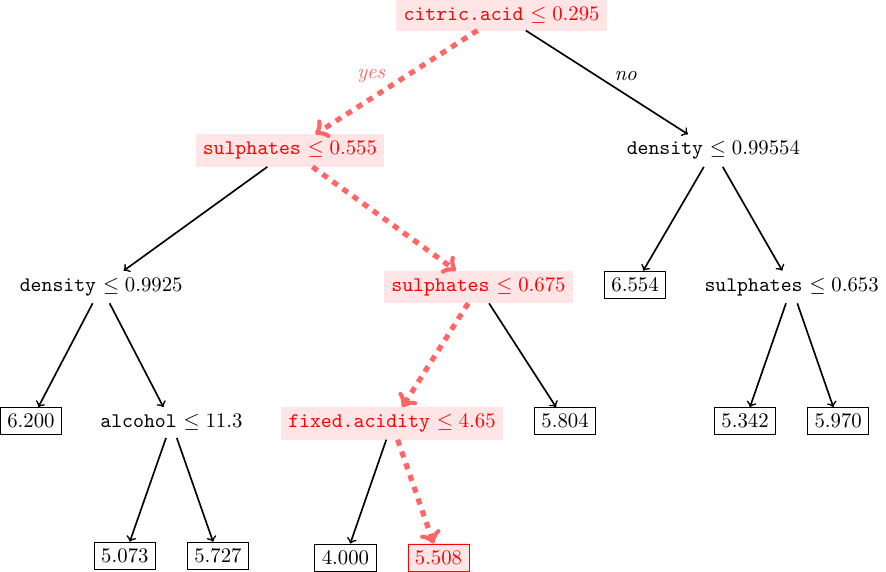} 
\caption{Example of a decision tree based on the \texttt{winequalityred} data set (see Section~\ref{subsec:computationalexperiments_background}). The goal is to predict the quality rating of a (red) wine using chemical properties of the wine. The shaded nodes and dashed edges indicate how an observation with \texttt{citric.acid} $= 0.22$, \texttt{density} $= 0.993$, \texttt{sulphates} $= 0.63$, \texttt{alcohol} $= 10.6$, \texttt{fixed.acidity} $= 4.9$, is mapped to a prediction (value of 5.508). \label{figure:decision_tree_example} }
\end{figure}

A decision tree is a model used for predicting a dependent variable $Y$ using a collection of independent variables $\Xb = (X_1, \dots, X_n)$. To make a prediction, we start at the root of the tree, and check a query (e.g., ``Is $X_3 \leq 5.6$?''); we then proceed to the left child node if the query is true, and to the right child node if the query is false. We then check the new node's query; the process continues until we reach a leaf node, where the tree outputs a prediction. A generalization of this type of model, called a \emph{tree ensemble model}, involves making this type of prediction from each of a collection of trees and aggregating the individual predictions into a single prediction (for example, by taking a weighted average of the predictions for a regression setting, or by taking a majority vote of the predictions for a classification setting). An example of a decision tree and a prediction being made is given in Figure~\ref{figure:decision_tree_example}.

Many types of tree ensemble models have been proposed in the machine learning literature; the most notable examples are random forests and boosted trees. Tree ensemble models are extremely attractive due to their ability to model complex, nonlinear relationships between the independent variables $\Xb$ and the dependent variable $Y$. %
As a result, tree ensemble models have gained or are gaining widespread popularity in a number of application areas, including chemistry \citep{svetnik2003random}, genomics \citep{diaz2006gene},  economics \citep{varian2014big,bajari2015machine}, marketing \citep{lemmens2006bagging} and operations management \citep{ferreira2015analytics}.  %

In many applications of tree ensemble models and predictive models in general, the independent variables that are used for prediction are exogenous and beyond our control as the decision maker. For example, one might build a random forest model to predict whether a patient is at risk of developing a disease based on the patient's age, blood pressure, family history of the disease and so on. Clearly, features such as age and family history are not amenable to intervention. Such models are typically used for some form of post hoc action or prioritization. In the disease risk example, one can use the random forest model to rank patients by decreasing predicted risk of an imminent acute event, and this ranking can be used to determine which patients should receive some intervention, e.g., closer monitoring, treatment with a particular drug and so on. 

In an increasing number of predictive modeling applications, however, some of the independent variables are controllable; that is to say, those independent variables are also \emph{decision} variables. We provide a couple of examples:
\begin{enumerate}
\item \textbf{Design of drug therapies}. In a healthcare context, one might be interested in building a model to predict patient response given a particular drug therapy, and then using such a model to find the optimal drug therapy. The dependent variable is some metric of the patient's health, and the independent variables may specify the drug therapy (which drugs and in what doses) and characteristics of the patient. A recent example of such an approach can be found in \cite{bertsimas2016analytics}, which considers the design of clinical trials for combination drug chemotherapy for gastric cancer; the first step involves estimating a model (specifically, a ridge regression model) of patient survival and toxicity using information about the drug therapy, and the second step involves solving an optimization problem to find the drug therapy that maximizes the predicted survival of the given patient group subject to a constraint on the predicted toxicity. 
\item \textbf{Pricing/promotion planning}. In marketing and operations management, a fundamental problem is that of deciding which products should be promoted when and at what price. In such a context, the data might consist of weekly sales of each product in a category, and the prices of the products in that category during that week and in previous weeks; one might use this data to build a predictive model of demand as a function of the prices, and then optimize such a model to decide on a promotion schedule \citep[for a recent example see, e.g.,][]{cohen2017impact}. 
\end{enumerate}

In this paper, we seek to unlock the prescriptive potential of tree ensemble models by considering the problem of \emph{tree ensemble optimization}. This problem is concerned with the following question: given a tree ensemble model that predicts some quantity $Y$ using a set of controllable independent variables $\Xb$, how should we set the independent variables $\Xb$ so as to maximize the predicted value of $Y$? This problem is of significant practical interest because it allows us to leverage the high accuracy afforded by tree ensemble models to obtain high quality decisions. At the same time, the problem is challenging, due to the highly nonlinear and large-scale nature of tree ensemble models.

We make the following contributions:
\begin{enumerate}
\item We propose the tree ensemble optimization problem and we show how to formulate this problem as a mixed-integer optimization (MIO) problem. The formulation can accommodate independent variables that are discrete, categorical variables as well as continuous, numeric variables. To the best of our knowledge, the problem of optimizing an objective function described as the prediction of a tree ensemble has not been previously proposed in either the machine learning or the operations research community. 
\item From a theoretical standpoint, we develop a number of results that generally concern the tractability of the formulation. First, we prove that the tree ensemble optimization problem is in general NP-Hard. We then show that our proposed MIO formulation offers a tighter relaxation of the problem than an alternate MIO formulation, obtained by applying a standard linearization to a binary polynomial formulation of the problem. We develop a hierarchy of approximate formulations for the problem, obtained by truncating each tree in the ensemble to a depth $d$ from the root node. We prove that the objective value of such an approximate formulation is an upper bound that improves as $d$ increases, and show how to construct a complementary a priori lower bound that depends on the variability of each tree's prediction below the truncation depth $d$. 
\item From a solution methodology standpoint, we present two different strategies for tackling large-scale instances of our MIO formulation. The first is based on solving a Benders reformulation of the problem using constraint generation. Here, we analyze the structure of the Benders subproblem and show that it can be solved efficiently. The second is based on applying lazy constraint generation directly to our MIO formulation. For this approach, we propose an efficient procedure for identifying violated constraints, which involves simply traversing each tree in the ensemble. 
\item We evaluate the effectiveness of our formulation and solution methods computationally using an assortment of real data sets. We show that the full MIO formulation can be solved to full optimality for small to medium sized instances within minutes, and that our formulation is significantly stronger in terms of relaxation bound and solution time than the aforementioned alternate formulation. We also show that our approach often significantly outperforms a simple local search heuristic that does not guarantee optimality. Lastly, we show that our customized solution methods can drastically reduce the solution time of our formulation.  
\item We provide a deeper showcase of the utility of our approach in two applications. The first is a case study in drug design, using a publicly available data set from Merck Research Labs \citep{ma2015deep}. Here, we show that our approach can optimize large-scale tree ensemble models with thousands of independent variables to full or near optimality within a two hour time limit, and can be used to construct a Pareto frontier of compounds that efficiently trade off predicted performance and similarity to existing, already-tested compounds. The second is a case study in customized pricing using a supermarket scanner data set \citep{montgomery1997creating}. Here, we show that a random forest model leads to considerable improvements in out-of-sample prediction accuracy over two state-of-the-art models based on hierarchical Bayesian regression, and that our optimization approach can find provably optimal prices at the individual store level within seconds. %
\end{enumerate}
The rest of the paper is organized as follows. In Section~\ref{sec:literaturereview}, we survey some of the related literature to this work. In Section~\ref{sec:model}, we present our formulation of the tree ensemble optimization problem as an MIO problem, and provide theoretical results on the structure of this problem. In Section~\ref{sec:solutionmethods}, we present two solution approaches for large-scale instances of the tree ensemble optimization problem. In Section~\ref{sec:computationalexperiments}, we present the results of our computational experiments with real data sets. In Section~\ref{sec:drugdesign} we present our case study in drug design and in Section~\ref{sec:pricing_summary}, we summarize our case study in customized pricing (described in greater detail in Section~\ref{sec:pricing} in the electronic companion). We give concluding remarks in Section~\ref{sec:conclusion}. %

\section{Literature review}
\label{sec:literaturereview}

Decision tree models became popular in machine learning with the introduction of two algorithms, ID3 \citep[iterative dichotomiser; see][]{quinlan1986induction} and CART \citep[classification and regression tree; see][]{breiman1984classification}. Decision tree models gained popularity due to their interpretability, but were found to be generally less accurate than other models such as linear and logistic regression. A number of ideas were thus consequently proposed for improving the accuracy of tree models, which are all based on constructing an \emph{ensemble} of tree models. \cite{breiman1996bagging} proposed the idea of bootstrap aggregation, or \emph{bagging}, where one builds a collection of predictive models, each trained with a bootstrapped sample of the original training set; the predictions of each model are then aggregated into a single prediction (by majority vote for classification and by averaging for regression). The motivation for bagging is that it reduces the prediction error for predictive models that are unstable/highly sensitive to the training data (such as CART); indeed, \cite{breiman1996bagging} showed that bagged regression trees can be significantly better than ordinary regression trees in out-of-sample prediction error. Later, \cite{breiman2001random} proposed the random forest model, where one builds a collection of bagged CART trees for which the subset of features selected for splitting at each node of each tree is randomly sampled from the set of all features \citep[the so-called \emph{random subspace} method; see][]{ho1998random}. Concurrently, other research has considered the idea of \emph{boosting} \citep{schapire2012boosting}, wherein one iteratively builds a weighted collection of basic predictive models (such as CART trees), with the goal of reducing the prediction error with each iteration. %

Tree ensembles occupy a central place in machine learning because they generally work very well in practice. In a systematic comparison of 179 different prediction methods on a broad set of benchmark data sets, \cite{fernandez2014we} found that random forests achieved best or near-best performance over all of these data sets. Boosted trees have been similarly successful: on the data science competition website Kaggle, one popular implementation of boosted trees, \texttt{XGBoost}, was used in more than half of the winning solutions in the year 2015 \citep{chen2016xgboost}. There exist robust and open source software implementations of many tree ensemble models. For boosted trees, the R package \texttt{gbm} \citep{ridgeway2006gbm} and \texttt{XGBoost} are widely used; for random forests, the R package \texttt{randomForest} \citep{liaw2002classification} is extremely popular. 

There has also been a significant concurrent effort to develop a theoretical foundation for tree ensemble methods; we briefly survey some of this work for random forests. For random forests, the original paper of \cite{breiman2001random} developed an upper bound on the generalization error of a random forest. Later research studied the consistency of simplified versions of the random forest model \citep[e.g.,][]{biau2008consistency} as well as the original model \citep[e.g.,][]{scornet2015consistency}. %
For an overview of recent theoretical advances in random forests, the reader is referred to \cite{Biau2016}.

At the same time, there is an increasing number of papers originating in operations research where a predictive model is used to represent the effect of the decision, and one solves an optimization problem to find the best decision with respect to this predictive model. Aside from the papers already mentioned in clinical trials and promotion planning, we mention two other examples. In pricing, data on historical prices and demand observed at those prices is often used to build a predictive model of demand as a function of price and to then determine the optimal price \citep[e.g.,][]{besbes2010testing}. In assortment optimization, \cite{misic2016data} considers a two-step approach, where the first step involves estimating a ranking-based choice model from historical data, and the second involves solving an MIO model to optimize the revenue predicted under that model.

In the research literature where predictive models are used to understand and subsequently optimize decisions, the closest paper conceptually to this one is the paper of \cite{ferreira2015analytics}. This paper considers the problem of pricing weekly sales for an online fashion retailer. To solve the problem, the paper builds a random forest model of the demand of each style included in a sale as a function of the style's price and the average price of the other styles. The paper then formulates an MIO problem to determine the optimal prices for the styles to be included in the sale, where the revenue is based on the price and the demand (as predicted by the random forest) of each style. The MIO formulation does not explicitly model the random forest prediction mechanism -- instead, one computes the random forest prediction for each style at each of its possible prices and at each possible average price of the other styles, and these predictions enter the MIO model as coefficients in the objective function. (The predictions could just as easily have come from a different form of predictive model, without changing the structure of the optimization problem.) In contrast, our MIO formulation explicitly represents the structure of each tree in the ensemble, allowing the prediction of each tree to be determined through the variables and constraints of the MIO. Although the modeling approach of \cite{ferreira2015analytics} is feasible for their pricing problem, it is difficult to extend this approach when there are many independent variables, as one would need to enumerate all possible combinations of values for them and compute the tree ensemble's prediction for each combination of values. To the best of our knowledge, our paper is the first to conceptualize the problem of how to optimize an objective function that is given by a tree ensemble.

Methodologically, the present paper is most related to the paper of \cite{bertsimas2018exact}, which considers the problem of designing a product line or assortment to optimize revenue under a ranking-based choice model. The ranking-based model considered in \cite{bertsimas2018exact} can be understood as a type of tree ensemble model; as such, the MIO formulation of \cite{bertsimas2018exact} can be regarded as a special case of the more general formulation that we analyze here. Some of the theoretical results found in \cite{bertsimas2018exact} -- specifically, those results on the structure of the Benders cuts -- are generalized in the present paper to tree ensemble models. Despite this similarity, the goals of the two papers are different. \cite{bertsimas2018exact} consider an optimization approach specifically for product line decisions, whereas in the present paper, our methodology can be applied to \emph{any} tree ensemble model, thus spanning a significantly broader range of application domains. Indeed, later in the paper we will present two different case studies -- one on drug design (Section~\ref{sec:drugdesign}) and one on customized pricing (Section~\ref{sec:pricing}) -- to illustrate the broad applicability of tree ensemble optimization.  %

Finally, we note that there is a growing literature on the use of mixed-integer optimization for the purpose of estimating decision tree models and other forms of statistical models. For example, \cite{bertsimas2017optimal} consider an exact MIO approach to constructing CART trees, while \cite{bertsimas2015algorithmic} consider an MIO approach to model selection in linear regression. While the present paper is related to this previous work in that it also leverages the technology of MIO, the goal of the present paper is different. The above papers focus on the estimation of trees and other statistical models, whereas our paper is focused on optimization, namely, how to determine the optimal \emph{decision} with respect to a given, fixed tree ensemble model.

\section{Model} 
\label{sec:model}

We begin by providing some background on tree ensemble models in Section~\ref{subsec:model_background} and defining the tree ensemble optimization problem. We then present our mixed-integer optimization model in Section~\ref{subsec:model_optimizationmodel}. We provide results on the strength of our formulation in Section~\ref{subsec:model_theoreticalproperties}. Finally, in Section~\ref{subsec:model_depthapproximation}, we describe a hierarchy of approximate formulations based on depth truncation. 

\subsection{Background}
\label{subsec:model_background}

In this section, we provide some background on tree ensemble models. We are given the task of predicting a dependent variable $Y$ using the independent variables $X_1, \dots, X_n$; for convenience, we use $\Xb$ to denote the vector of independent variables. We let $\Xcal_i$ denote the domain of independent variable $i$ and let $\Xcal = \prod_{i=1}^n \Xcal_i$ denote the domain of $\Xb$. An independent variable $i$ may be a numeric variable or a categorical variable. 

A decision tree is a model for predicting the dependent variable $Y$ using the independent variable $\Xb$ by checking a collection of \emph{splits}. A split is a condition or query on a single independent variable that is either true or false. More precisely, for a numeric variable $i$, a split is a query of the form
\begin{equation*}
\text{Is $X_i \leq a$?}
\end{equation*}
for some $a \in \mathbb{R}$. For a categorical variable $i$, a split is a query of the form
\begin{equation*}
\text{Is $X_i \in A$?}
\end{equation*}
where $A \subseteq \Xcal_i$ is a set of levels of the categorical variable. The splits are arranged in the form of a tree, with each split node having two child nodes. The left child corresponds to the split condition being true, while the right child corresponds to the condition being false. To make a prediction for an observation with the independent variable $\Xb$, we start at the root of the tree and check whether $\Xb$ satisfies the split condition; if it is true, we move to the left child, and if it is false, we move to the right child. At the new node, we check the split again, and move again to the corresponding node. This process continues until we reach a leaf of the tree. The prediction that we make is the value corresponding to the leaf we have reached. 

In this paper, we will focus on predictive models that are ensembles or collections of decision trees. We assume that there are $T$ trees, where each tree is indexed from 1 to $T$. Each tree $t$ has a weight $\lambda_t$, and its prediction is denoted by the function $f_t$; for the independent variable $\Xb$, the prediction of tree $t$ is $f_t(\Xb)$. For an observation with independent variable $\Xb$, the prediction of the ensemble of trees is given by 
\begin{equation*}
\sum_{t=1}^T \lambda_t f_t( \Xb). 
\end{equation*}

The optimization problem that we would like to solve is to find the value of the independent variable $\Xb$ that maximizes the ensemble prediction:
\begin{equation}
\underset{\Xb \in \Xcal}{\text{maximize}} \, \sum_{t=1}^T \lambda_t f_t( \Xb). \label{prob:TreeEnsembleOpt_abstract}
\end{equation}
We shall make two key assumptions about the tree ensemble model $\sum_{t=1}^T \lambda_t f_t(\cdot)$ and our tree ensemble optimization problem~\eqref{prob:TreeEnsembleOpt_abstract}: 
\begin{enumerate}
\item First, \emph{we assume that we are only making a single, one-time decision and that the tree ensemble model is fixed}. Extending our approach to the multistage setting is an interesting direction for future research. 
\item Second, \emph{we assume that the tree ensemble model is an accurate representation of the outcome when we make the decision $\Xb$}. In practice, some care must be taken here because depending on how the tree ensemble model is estimated and the nature of the data, the prediction $\sum_{t=1}^T \lambda_t f_t(\Xb)$ may not necessarily be an accurate estimate of the causal effect of setting the independent variables to $\Xb$. This issue has been the focus of some recent work in prescriptive analytics \citep[see][]{bertsimas2016pricing,kallus2016recursive}. Our goal in this paper is to address only the question of \emph{optimization} -- how to efficiently and scalably optimize a tree ensemble function $\sum_{t=1}^T \lambda_t f_t(\cdot)$ -- which is \emph{independent} of such statistical questions. As such, we will assume that the tree ensemble model we are given at the outset is beyond suspicion. 
\end{enumerate}

Problem~\eqref{prob:TreeEnsembleOpt_abstract} is very general, and one question we may have is whether it is theoretically tractable or not. Our first theoretical result answers this question in the negative. 
\begin{proposition}
The tree ensemble optimization problem~\eqref{prob:TreeEnsembleOpt_abstract} is NP-Hard.  \label{prop:TEO_NPComplete}
\end{proposition}
The proof of Proposition~\ref{prop:TEO_NPComplete}, given in Section~\ref{appendix:proof_TEO_NPComplete} of the e-companion, uses a reduction from the minimum vertex cover problem.

\subsection{Optimization model}
\label{subsec:model_optimizationmodel}

We now present an MIO formulation of \problemeqref{prob:TreeEnsembleOpt_abstract}. Before we present the model, we will require some additional notation. We let $\Ncal$ denote the set of numeric variables and $\Ccal$ denote the set of categorical variables; we have that $\Ncal \cup \Ccal = \{1,\dots, n\}$. 

For each numeric variable $i \in \Ncal$, let $\Acal_i$ denote the set of unique split points, that is, the set of values $a$ such that $X_i \leq a$ is a split condition in some tree in the ensemble $\{ f_t \}_{t=1}^T$. Let $K_i = | \Acal_i |$ be the number of unique split points. Let $a_{i,j}$ denote the $j$th smallest split point of variable $i$, so that $a_{i,1} < a_{i,2} < \dots < a_{i,K_i}$.  %

For each categorical variable $i \in \Ccal$, recall that $\Xcal_i$ is the set of possible values of $i$. For convenience, we use $K_i$ in this case to denote the size of $\Xcal_i$ (\ie, $K_i = | \Xcal_i|$) and use the values $1, 2, \dots, K_i$ to denote the possible levels of variable $i$. 

Let $\leaves(t)$ be the set of leaves or terminal nodes of tree $t$. Let $\splits(t)$ denote the set of splits of tree $t$ (non-terminal nodes). Recall that the left branch of the split corresponds to ``yes'' or ``true'' to the split query, and the right branch corresponds to ``no'' or ``false''. Therefore, for each split $s$ in $S_t$, we let $\leftleaves(s)$ be the set of leaves that are accessible from the left branch (all of the leaves for which the condition of split $s$ must be true), and $\rightleaves(s)$ be the set of leaves that are accessible from the right branch (all of the leaves for which the condition of split $s$ must be false). For each split $s$, we let $\var(s) \in \{1,\dots, n\}$ denote the variable that participates in split $s$, and let $\cond(s)$ denote the set of values of variable $i$ that participate in the split query of $s$. Specifically, if $\var(s)$ is numeric, then $\cond(s) = \{ j\}$ for some $j \in \{1,\dots, K_{\var(s)}\}$, which corresponds to the split query $X_i \leq a_{i,j}$. If $\var(s)$ is categorical, then $\cond(s) \subseteq \{1,\dots, K_{\var(s)}\}$, which corresponds to the query $X_i \in \cond(s)$. For each leaf $\ell \in \leaves(t)$, we use $p_{t,\ell}$ to denote the prediction that tree $t$ makes when an observation reaches leaf $\ell$.

We now define the decision variables of the problem. There are two sets of decision variables. The first set is used to specify the independent variable value $\Xb$. For each categorical independent variable $i \in \Ccal$ and each category/level $j \in \Xcal_i$, we let $x_{i,j}$ be 1 if independent variable $i$ is set to level $j$, and 0 otherwise. For each numeric independent variable $i \in \Ncal$ and each $j \in \{1,\dots,K_i\}$, we let $x_{i,j}$ be 1 if independent variable $i$ is set to a value less than or equal to the $j$th split point, and 0 otherwise. Mathematically, 
\begin{align*}
x_{i,j} & = \Ibb\{ X_i = j \}, \quad \forall i \in \Ccal, \ j \in \{1,\dots,K_i\}, \\
x_{i,j} & = \Ibb\{ X_i \leq a_{i,j} \}, \quad \forall i \in \Ncal, \ j \in \{1,\dots, K_i\}.
\end{align*}
We use $\xb$ to denote the vector of $x_{i,j}$ values. 

The second set of decision variables is used to specify the prediction of each tree $t$. For each tree $t$ and each leaf $\ell \in \leaves(t)$, we let $y_{t,\ell}$ be a binary decision variable that is 1 if the observation encoded by $\xb$ belongs to/falls into leaf $\ell$ of tree $t$, and 0 otherwise. 

With these definitions, the MIO can be written as follows:
{\allowdisplaybreaks
\begin{subequations}
\begin{alignat}{2}
& \underset{\xb,\yb}{\text{maximize}} \quad & & \sum_{t=1}^T \sum_{\ell \in \leaves(t)}  \lambda_t \cdot p_{t,\ell} \cdot y_{t,\ell} \label{prob:TEOMIO_objective}\\
& \text{subject to}  & & \sum_{\ell \in \leaves(t)} y_{t,\ell} = 1, \quad \forall \ t \in \{1,\dots,T\},  \label{prob:TEOMIO_ysumtoone}\\
&  & & \sum_{\ell \in \leftleaves(s)} y_{t,\ell} \leq \sum_{j \in \cond(s)} x_{\var(s),j}, \nonumber \\
& & & \qquad  \forall \ t \in \{1,\dots, T\},\ s \in \splits(t), \label{prob:TEOMIO_left} \\[0.5em]
& & &  \sum_{\ell \in \rightleaves(s)} y_{t,\ell} \leq 1 - \sum_{j \in \cond(s)} x_{\var(s),j},  \nonumber \\
& & & \qquad \forall \ t \in \{1,\dots, T\},\ s \in \splits(t), \label{prob:TEOMIO_right} \\
& & & \sum_{j=1}^{K_i} x_{i,j} = 1, \quad \forall \ i \in \Ccal, \label{prob:TEOMIO_categorical}\\
& & & x_{i,j} \leq x_{i,j+1}, \quad \forall \ i \in \Ncal, \ j \in \{1,\dots, K_i-1\},  \label{prob:TEOMIO_numeric}\\
& & & x_{i,j} \in \{0,1\}, \quad \forall \ i \in \{1,\dots, n\}, \ j \in \{1,\dots, K_i\} \label{prob:TEOMIO_xbinary}\\
& & & y_{t,\ell} \geq 0, \quad \forall \ t \in \{1,\dots, T\}, \ \ell \in \leaves(t).  \label{prob:TEOMIO_ycontinuous} 
\end{alignat}
\label{prob:TEOMIO}%
\end{subequations}}%
The constraints have the following meaning. Constraint~\eqref{prob:TEOMIO_ysumtoone} ensures that the observation falls in exactly one of the leaves of each tree $t$. Constraint~\eqref{prob:TEOMIO_left} ensures that, if $\sum_{j \in \cond(s)} x_{\var(s), j} = 0$, then $y_{t,\ell}$ is forced to zero for all $\ell \in \leftleaves(s)$; in words, the condition of the split is false, so the observation cannot fall into any leaf to the left of split $s$, as this would require the condition to be true. Similarly, constraint~\eqref{prob:TEOMIO_right} ensures that if the condition of split $s$ \emph{is} satisfied, then $y_{t,\ell}$ is forced to zero for all $\ell \in \rightleaves(s)$; in words, the condition of the split is true, so the observation cannot fall into any leaf to the right of split $s$, as this would require the condition to be false. Constraint~\eqref{prob:TEOMIO_categorical} ensures that for each categorical variable $i \in \Ccal$, exactly one of the $K_i$ levels is selected. Constraint~\eqref{prob:TEOMIO_numeric} requires that if the numeric independent variable $i$ is less than or equal to the $j$th lowest split point, then it must also be less than or equal to the $(j+1)$th lowest split point. Constraint~\eqref{prob:TEOMIO_xbinary} defines each $x_{i,j}$ to be binary, while constraint~\eqref{prob:TEOMIO_ycontinuous} defines each $y_{t,\ell}$ to be nonnegative. The objective represents the prediction of the ensemble of trees on the observation that is encoded by $\xb$. 

We now comment on several features of the model. The first is that the $y_{t,\ell}$ variables, despite having a binary meaning, are defined as continuous variables. The reason for this is that when $\xb$ is binary, the constraints automatically force $\yb$ to be binary.  We will formally state this result later (Proposition~\ref{prop:SubPrimalBinaryOptimal} of Section~\ref{subsec:solutionmethods_bendersdecomposition}). As a result, the only binary variables are those in $\xb$, of which there are $\sum_{i=1}^n K_i$. Recall that for categorical independent variables, $K_i$ is the number of levels, whereas for numeric independent variables, $K_i$ is the number of unique split points found in the tree ensemble $\{ f_t \}_{t=1}^T$.

The second is that our formulation does not model the exact value of each numeric independent variable $i \in \Ncal$. In contrast, the formulation only models where the variable is in relation to the unique split points in $\Acal_i$ -- for example, $x_{i,1} = 1$ indicates that independent variable $i$ is set to be less than or equal to the first lowest split point. The reason for this is that each decision tree function $f_t(\cdot)$ is a piecewise constant function and therefore the tree ensemble function $\sum_{t=1}^T \lambda_t f_t(\Xb)$ is also piecewise constant. Thus, for the purpose of optimizing the function $\sum_{t=1}^T \lambda_t f_t(\Xb)$, it is not necessary to explicitly maintain the value $X_i$ of each numeric independent variable $i$. 

The third is that numeric independent variables are modeled in terms of an inequality, that is, $x_{i,j} = \Ibb\{X_i \leq a_{i,j}\}$. Alternatively, one could model numeric independent variables by using $x_{i,j}$ to represent whether $X_i$ is between two consecutive split points, e.g., %
\begin{align*}
x_{i,1} & = \Ibb\{ X_i \leq a_{i,1} \}, \\
x_{i,2} & = \Ibb\{ a_{i,1} < X_i \leq a_{i,2} \}, \\
x_{i,3} & = \Ibb\{ a_{i,2} < X_i \leq a_{i,3} \}, \\
& \vdots \\
x_{i,K_i} & = \Ibb\{ a_{i,K_i - 1} < X_i \leq a_{i, K_i} \}, \\
x_{i,K_i+1} & = \Ibb\{ a_{i,K_i} > X_i  \}.
\end{align*}
One would then re-define the set $\cond(s)$ for each split involving the variable $i$ so as to include all of the relevant $j$ values under this new encoding. The advantage of our choice of encoding -- using $x_{i,j} = \Ibb \{ X_i \leq a_{i,j} \}$ -- is that the resulting formulation enhances the power of branching on fractional values of $x_{i,j}$ and leads to more balanced branch-and-bound trees \citep{vielma2015mixed}. This type of encoding has been used successfully in scheduling and transportation applications \citep[so-called ``by'' variables, representing an event happening by some period $t$; see for example][]{bertsimas2011integer}; for further details, the reader is referred to \cite{vielma2015mixed}.

\subsection{Theoretical properties}
\label{subsec:model_theoreticalproperties}

We now compare our formulation against an alternate MIO formulation of the tree ensemble optimization problem, which involves relating the $\yb$ and $\xb$ variables in a different way.

In particular, for any leaf $\ell$ of any tree $t$, let $\leftsplits(\ell)$ be the set of splits for which leaf $\ell$ is on the left side (\ie, $s$ such that $\ell \in \leftleaves(s)$), and $\rightsplits(\ell)$ be the set of splits for which leaf $\ell$ is on the right side (\ie, $s$ such that $\ell \in \rightleaves(s)$). The tree ensemble optimization problem can then be formulated as the following problem:
\begin{subequations}
\begin{alignat}{2}
& \underset{\xb}{\text{maximize}} \quad & & \sum_{t=1}^T \sum_{\ell \in \leaves(t)} \lambda_t \cdot p_{t,\ell} \cdot \prod_{s \in \leftsplits(\ell)} \left( \sum_{j \in \cond(s)} x_{\var(s), j} \right) \cdot \prod_{s \in \rightsplits(\ell)} \left(1 - \sum_{j \in \cond(s)} x_{\var(s),j}  \right) \\
& \text{subject to} & & \text{constraints~\eqref{prob:TEOMIO_categorical}-\eqref{prob:TEOMIO_xbinary}}.
\end{alignat}
\label{prob:TEOPolynomial}%
\end{subequations}
The above problem is a binary polynomial problem. Note that the product term, $\prod_{s \in \leftsplits(\ell)} \left( \sum_{j \in \cond(s)} x_{\var(s), j} \right) \cdot \prod_{s \in \rightsplits(\ell)} \left(1 - \sum_{j \in \cond(s)} x_{\var(s),j}  \right)$, is exactly 1 if the observation is mapped to leaf $\ell$ of tree $t$, and 0 otherwise. The standard linearization of \problemeqref{prob:TEOPolynomial} \cite[see][]{crama1993concave} is the following MIO:
{\allowdisplaybreaks
\begin{subequations}
\begin{alignat}{2}
& \underset{\xb,\yb}{\text{maximize}} \quad & & \sum_{t=1}^T \sum_{\ell \in \leaves(t)} \lambda_t \cdot p_{t,\ell} \cdot y_{t,\ell} \\
& \text{subject to} & & y_{t,\ell} \leq \sum_{j \in \cond(s)} x_{\var(s),j}, \quad \forall t \in \{1,\dots,T\}, \ \ell \in \leaves(t), \ s \in \leftsplits(\ell), \label{prob:TEOLinearized_left}\\
& & & y_{t,\ell} \leq 1 - \sum_{j \in \cond(s)} x_{\var(s),j},  \quad \forall t \in \{1,\dots,T\}, \ \ell \in \leaves(t), \ s \in \rightsplits(\ell), \label{prob:TEOLinearized_right} \\
& & & y_{t,\ell} \geq \sum_{s \in \leftsplits(\ell)} \left( \sum_{j \in \cond(s)} x_{\var(s),j} \right)  + \sum_{s \in \rightsplits(\ell)} \left(  1 - \sum_{j \in \cond(s)} x_{\var(s),j} \right)  - ( | \leftsplits(\ell)| + |\rightsplits(\ell)| - 1), \nonumber \\
& & & \qquad \forall t \in \{1,\dots,T\}, \ \ell \in \leaves(t), \label{prob:TEOLinearized_lower}  \\
& & & \text{constraints~\eqref{prob:TEOMIO_categorical}-\eqref{prob:TEOMIO_ycontinuous}}. 
\end{alignat}
\label{prob:TEOLinearized}
\end{subequations}
}
Let $Z^*_{LO}$ be the optimal value of the linear optimization (LO) relaxation of \problemeqref{prob:TEOMIO} and let $Z^*_{LO,StdLin}$ be the optimal value of the LO relaxation of \problemeqref{prob:TEOLinearized}. The following result relates the two optimal values.
\begin{proposition}
$Z^*_{LO} \leq Z^*_{LO,StdLin}$. \label{prop:relaxationstrength}
\end{proposition}
The proof of Proposition~\ref{prop:relaxationstrength} (see Section~\ref{appendix:proof_relaxationstrength}) consists of showing that an optimal solution of the relaxation of \problemeqref{prob:TEOMIO} is a feasible solution for the relaxation of \problemeqref{prob:TEOLinearized} and achieves an objective value of exactly $Z^*_{LO}$. The significance of Proposition~\ref{prop:relaxationstrength} is that it establishes that formulation~\eqref{prob:TEOMIO} is a stronger formulation of the tree ensemble optimization problem than formulation~\eqref{prob:TEOLinearized}. This is desirable from a practical perspective, as stronger MIO formulations are generally faster to solve than weaker MIO formulations. We shall see in Section~\ref{subsec:computationalexperiments_fullmio} that the difference in relaxation bounds can be substantial and that formulation~\eqref{prob:TEOLinearized} is significantly less tractable than our formulation~\eqref{prob:TEOMIO}. 

\subsection{Depth $d$ approximation}

\label{subsec:model_depthapproximation}

In this section, we describe a hierarchy of relaxations of \problemeqref{prob:TEOMIO} that are based on approximating each tree in the ensemble up to a particular depth. 

The motivation for this hierarchy of relaxations comes from the following observation regarding the size of \problemeqref{prob:TEOMIO}. 
In particular, a key driver of the size of \problemeqref{prob:TEOMIO} is the number of left and right split constraints \eqref{prob:TEOMIO_left} and \eqref{prob:TEOMIO_right}, respectively; these constraints are enforced for every single split in each tree in the ensemble. For a large number of trees that are deep (and thus have many splits), the resulting number of left and right split constraints will be large. At the same time, it may be reasonable to expect that if we do not represent each tree to its full depth, but instead only represent each tree up to some depth $d$ and only include splits that occur before (and including) depth $d$, then we might still obtain a reasonable solution to \problemeqref{prob:TEOMIO}. In this section, we rigorously define this hierarchy of approximate formulations, and provide theoretical guarantees on how close such approximations are to the original formulation. 

Let $\Omega = \{ (t,s) \, | \, t \in \{1,\dots, T\}, \ s \in \splits(t) \}$ be the set of tree-split pairs. Let $\bar{\Omega} \subseteq \Omega$ be a subset of all possible tree-split pairs. The $\bar{\Omega}$ tree ensemble problem is defined as problem~\eqref{prob:TEOMIO} where constraints~\eqref{prob:TEOMIO_left} and \eqref{prob:TEOMIO_right} are restricted to $\bar{\Omega}$: 
\begin{subequations}
\begin{alignat}{2}
& \underset{\xb,\yb}{\text{maximize}} \quad & & \sum_{t=1}^T \sum_{\ell \in \leaves(t)}  \lambda_t \cdot p_{t,\ell} \cdot y_{t,\ell} \label{prob:TEOMIO_barOmega_objective}\\
& \text{subject to}  & & \sum_{\ell \in \leftleaves(s)} y_{t,\ell} \leq \sum_{j \in \cond(s)} x_{\var(s),j}, \quad \forall \ (t,s) \in \bar{\Omega}, \label{prob:TEOMIO_barOmega_left} \\[0.5em]
& & &  \sum_{\ell \in \rightleaves(s)} y_{t,\ell} \leq 1 - \sum_{j \in \cond(s)} x_{\var(s),j},  \quad \forall \ (t,s) \in \bar{\Omega} , \label{prob:TEOMIO_barOmega_right} \\
& & & \text{constraints~\eqref{prob:TEOMIO_ysumtoone}, \eqref{prob:TEOMIO_categorical} - \eqref{prob:TEOMIO_ycontinuous}}.
\end{alignat}
\label{prob:TEOMIO_barOmega}%
\end{subequations} %
Problem~\eqref{prob:TEOMIO_barOmega} is obtained by removing constraints from \problemeqref{prob:TEOMIO}, while retaining the same decision variables. Solving problem~\eqref{prob:TEOMIO_barOmega} for a fixed $\bar{\Omega}$ will result in a solution $(\xb, \yb)$ that is only guaranteed to satisfy constraints~\eqref{prob:TEOMIO_left} and \eqref{prob:TEOMIO_right} for those $(t,s)$ pairs in $\bar{\Omega}$. As a result, for any $\bar{\Omega} \subseteq \Omega$, the objective value of \problemeqref{prob:TEOMIO_barOmega} will be a valid upper bound on \problemeqref{prob:TEOMIO}. We will now define a collection of subsets $\bar{\Omega}$, obtained by truncating each tree at a particular depth, for which we also have an easily-computable accompanying lower bound.

For any $d \in \mathbb{Z}_+$, let $\bar{\Omega}_d$ be the set of all tree-split pairs where the split is at a depth $d' \leq d$ (a depth of 1 corresponds to the split at the root node). Let $Z^*_{MIO,d}$ denote the objective value of \problemeqref{prob:TEOMIO_barOmega} with $\bar{\Omega}_d$, \ie, all splits up to and including depth $d$, and let $d_{\max}$ be the maximum depth of any split in any tree of the ensemble. Let $Z^*_{MIO}$ be the objective value of \problemeqref{prob:TEOMIO}, where the split constraints are up to the full depth of each tree; note that $Z^*_{MIO} = Z^*_{MIO,d_{\max}}$. The following result establishes that $Z^*_{MIO,d}$ is an upper bound on $Z^*_{MIO}$ that becomes tighter as $d$ increases:
\begin{proposition}
$Z^*_{MIO,1} \geq Z^*_{MIO,2} \geq \dots \geq Z^*_{MIO,d_{\max}} = Z^*_{MIO}.$ \label{prop:depth_ordering}
\end{proposition}
We now show how to construct a complementary lower bound. %

Fix some depth $d$. Let $\splits(t,d)$ denote the set of splits at depth $d$; if the depth of the tree is strictly less than $d$, $\splits(t,d)$ is empty. Let us define the constant $\delta_{t,s}$ for each split $s \in \splits(t,d)$ of each tree $t$ as 
\begin{equation}
\delta_{t,s} = \max \left\{ \max_{\ell \in \leftleaves(s)} p_{t,\ell} - \min_{\ell \in \leftleaves(s)} p_{t,\ell}, \max_{\ell \in \rightleaves(s)} p_{t,\ell} - \min_{\ell \in \rightleaves(s)} p_{t,\ell} \right\}.
\end{equation}
The constant $\delta_{t,s}$ is an upper bound on the maximum error possible (due to the depth $d$ truncation of the split constraints) in the prediction of tree $t$ for the observation encoded by $\xb$, given that the observation reaches split $s$. Stated differently, if the observation $\xb$ reaches split $s$, the error from omitting all splits below $s$ cannot be more than the difference between the highest and lowest leaf prediction on the left subtree or the right subtree at $s$.

We define $\Delta_t$ as the maximum of the $\delta_{t,s}$ values over all the depth $d$ splits of tree $t$:
\begin{equation*}
\Delta_t = \max_{s \in \splits(t,d)} \delta_{t,s}. 
\end{equation*}
(In the case that $\splits(t,d)$ is empty, we set $\Delta_t = 0$.) %

Before stating our approximation guarantee, we note that given a solution $(\xb,\yb)$ that solves problem~\eqref{prob:TEOMIO_barOmega} with $\bar{\Omega}_d$, it is possible to find a solution $(\xb,\tilde{\yb})$ that is a feasible solution for the full depth \problemeqref{prob:TEOMIO}. %
Our approximation guarantee is given below. 
\begin{theorem}
Suppose that $\lambda_t \geq 0$ for all $t \in \{1,\dots,T\}$ and $d \in \mathbb{Z}_+$. Let $(\xb, \yb)$ be an optimal solution of problem~\eqref{prob:TEOMIO_barOmega} with $\bar{\Omega}_d$. Let $Z_d$ be the true objective of $\xb$ when embedded within the full-depth problem~\eqref{prob:TEOMIO_barOmega}. We then have
\begin{equation*}
Z^*_{MIO,d} - \sum_{t=1}^T \lambda_t \Delta_t \leq Z_d \leq Z^*_{MIO} \leq Z^*_{MIO,d}. 
\end{equation*}
\label{theorem:depth_approximation}
\end{theorem}
The above theorem, which we prove in Section~\ref{appendix:proof_depth_approximation}, provides a guarantee on how suboptimal the $\bar{\xb}$ solution, derived from the depth $d$ problem~\eqref{prob:TEOMIO_barOmega} with $\bar{\Omega}_d$, is for the true (full depth) problem~\eqref{prob:TEOMIO}. Note that the requirement of $\lambda_t$ being nonnegative is not particularly restrictive, as we can always make $\lambda_t$ of a given tree $t$ positive by negating the leaf predictions $p_{t,\ell}$ of that tree. This result is of practical relevance because it allows the decision maker to judiciously trade-off the complexity of the problem (represented by the depth $d$) against an a priori guarantee on the quality of the approximation. Moreover, the quantity $\sum_{t=1}^T \lambda_t \Delta_t$, which bounds the difference between $Z^*_{MIO,d}$ and $Z_d$, can be easily computed from each tree, allowing the bound to be readily implemented in practice. We shall see in Section~\ref{subsec:computationalexperiments_depth} that although the lower bound can be rather conservative for small values of $d$, the true objective value of $\bar{\xb}$ is often significantly better.

Theorem~\ref{theorem:depth_approximation} suggests an interesting question regarding how the trees that are input to the tree ensemble optimization problem are estimated. In particular, for any given tree $t$, it may be possible to re-arrange the splits in the tree in a way that one obtains a new tree, $t'$, such that $t'$ gives identical predictions to $t$ for all $\Xb \in \Xcal$. We may therefore ask: is there a way to re-arrange the splits of the trees so as to obtain a tighter upper bound $Z^*_{MIO,d}$, and a smaller approximation error bound $\sum_{t=1}^T \lambda_t \Delta_t$? This question constitutes an interesting direction for future work.

\section{Solution methods}
\label{sec:solutionmethods}

The optimization model in Section~\ref{subsec:model_optimizationmodel}, although tractable for small to medium instances, can be difficult to solve directly for large instances. In this section, we present two solution approaches for tackling large-scale instances of \problemeqref{prob:TEOMIO}. In Section~\ref{subsec:solutionmethods_bendersdecomposition}, we present an approach based on Benders decomposition. In Section~\ref{subsec:solutionmethods_primalconstraintgeneration}, we present an alternate approach based on iteratively generating the split constraints.

\subsection{Benders decomposition}
\label{subsec:solutionmethods_bendersdecomposition}

The first solution approach that we present is Benders decomposition. Recall that in \problemeqref{prob:TEOMIO}, we have two sets of variables, $\xb$ and $\yb$; furthermore, $\yb$ can be further partitioned as $\yb = (\yb_1, \yb_2, \dots, \yb_T)$, where $\yb_t$ is the collection of $y_{t,\ell}$ variables corresponding to tree $t$. For any two trees $t, t'$ with $t \neq t'$, notice that the variables $\yb_t$ and $\yb_{t'}$ do not appear together in any constraints; they are only linked through the $\xb$ variables. 

The above observation suggests a Benders reformulation of \problemeqref{prob:TEOMIO}. Re-write \problemeqref{prob:TEOMIO} as follows:
\begin{subequations}
\begin{alignat}{2}
& \underset{\xb}{\text{maximize}} \quad & & \sum_{t=1}^T \lambda_t G_t(\xb), \\
& \text{subject to} & & \text{constraints~\eqref{prob:TEOMIO_categorical} - \eqref{prob:TEOMIO_xbinary}},
\end{alignat}
\label{prob:TEOMIO_xonly}%
\end{subequations}
where $G_t(\xb)$ is defined as the optimal value of the following subproblem:
\begin{subequations}
\begin{alignat}{2}
G_t(\xb) = \quad & \underset{\yb_t}{\text{maximize}} \quad & & \sum_{\ell \in \leaves(t)} p_{t,\ell} \cdot y_{t,\ell} \\
& \text{subject to} & & \sum_{\ell \in \leaves(t)} y_{t,\ell} = 1, \label{prob:SubPrimal_ysumtoone}\\
&  & & \sum_{\ell \in \leftleaves(s)} y_{t,\ell} \leq \sum_{j \in \cond(s)} x_{\var(s),j}, \quad  \forall \ s \in \splits(t), \label{prob:SubPrimal_left} \\[0.5em]
& & &  \sum_{\ell \in \rightleaves(s)} y_{t,\ell} \leq 1 - \sum_{j \in \cond(s)} x_{\var(s),j},  \quad \forall \ s \in \splits(t), \label{prob:SubPrimal_right} \\
& & & y_{t,\ell} \geq 0, \quad \forall\  \ell \in \leaves(t) \label{prob:SubPrimal_ynonnegative}.
\end{alignat}
\label{prob:SubPrimal}%
\end{subequations}
The first result we will prove is on the form of the optimal solution to \problemeqref{prob:SubPrimal}. To do this, we first provide a procedure in Algorithm~\ref{algorithm:SubPrimalGetLeaf} for determining the leaf of the solution encoded by $\xb$ for tree $t$. For ease of exposition, we will denote this procedure applied to a particular observation encoded by $\xb$ and a given tree $t$ as $\getLeaf(\xb, t)$. We use $\leftchild(\nu)$ to denote the left child of a split node $\nu$, $\rightchild(\nu)$ to denote the right child of a split node $\nu$, and $\root(t)$ to denote the root split node of tree $t$. Note that if $d_{\max}$ is the maximum depth of the trees in the ensemble (\ie, the maximum depth of any split in any tree), then the time complexity of $\getLeaf(\xb,t)$ is $O(d_{\max})$, as we will traverse at most $d_{\max}$ split nodes to reach a leaf.
\begin{algorithm}
\small
\caption{$\getLeaf$ procedure for determining the leaf to which tree $t$ maps $\xb$.}
\label{algorithm:SubPrimalGetLeaf}
\begin{algorithmic}
\STATE Initialize $\nu \leftarrow \root(t)$ 
\WHILE{$\nu \notin \leaves(t)$}
	\IF{$\sum_{j \in \cond(\nu)} x_{\var(\nu),j} =1$}
	\STATE $\nu \leftarrow \leftchild(\nu)$
	\ELSE
	\STATE $\nu \leftarrow \rightchild(\nu)$
	\ENDIF
\ENDWHILE
\RETURN $\nu$
\end{algorithmic}
\end{algorithm}

Having defined $\getLeaf$, we now present our first theoretical result (see Section~\ref{appendix:proof_SubPrimalBinaryOptimal} for the proof).
\begin{proposition}
Let $\xb \in \{0,1\}^{\sum_{i=1}^n K_i}$ be a feasible solution of \problemeqref{prob:TEOMIO_xonly}. Let $\ell^* = \getLeaf(\xb,t)$ be the leaf into which $\xb$ falls, and let $\yb_t$ be the solution to \problemeqref{prob:SubPrimal} defined as 
\begin{equation*}
y_{t,\ell} = \left\{ \begin{array}{ll} 1 & \text{if} \ \ell = \ell^*, \\ 0 & \text{otherwise}. \end{array} \right. 
\end{equation*}
The solution $\yb_t$ is the only feasible solution and therefore, the optimal solution of \problemeqref{prob:SubPrimal}.
\label{prop:SubPrimalBinaryOptimal}
\end{proposition}
Since \problemeqref{prob:SubPrimal} is feasible and has a finite optimal value, then by LO strong duality the optimal value of \problemeqref{prob:SubPrimal} is equal to the optimal value of its dual. The dual of subproblem~\eqref{prob:SubPrimal} is
\begin{subequations}
\begin{alignat}{2}
& \underset{\alphab_t, \betab_t, \gamma_t}{\text{minimize}} \quad & & \sum_{s \in \splits(t)} \alpha_{t,s} \left[ \sum_{j \in \cond(s)} x_{\var(s),j} \right] \nonumber \\
& & & + \sum_{s \in \splits(t)} \beta_{t,s} \left[ 1 - \sum_{j \in \cond(s)} x_{\var(s),j} \right] + \gamma_t \\[0.5em]
& \text{subject to} & & \sum_{s: \ell \in \leftleaves(s)} \alpha_{t,s} + \sum_{s: \ell \in \rightleaves(s)} \beta_{t,s} + \gamma_t \geq p_{t,\ell}, \quad \forall \ \ell \in \leaves(t), \label{prob:SubDual_constraint} \\
& & & \alpha_{t,s}, \beta_{t,s} \geq 0, \quad \forall \ s \in \splits(t). 
\end{alignat}
\label{prob:SubDual}%
\end{subequations}
Letting $\Dcal_{t}$ denote the set of dual feasible solutions $(\alphab_t, \betab_t, \gamma_t)$ for subproblem $t$, we can re-write \problemeqref{prob:TEOMIO_xonly} as 
\begin{subequations}
\begin{alignat}{2}
& \underset{\xb, \thetab}{\text{maximize}} & & \sum_{t=1}^T  \lambda_t \theta_t \\
& \text{subject to} & & \sum_{s \in \splits(t)} \alpha_{t,s} \left[ \sum_{j \in \cond(s)} x_{\var(s),j} \right] + \sum_{s \in \splits(t)} \beta_{t,s} \left[ 1 - \sum_{j \in \cond(s)} x_{\var(s),j} \right] + \gamma_t \geq \theta_t, \nonumber \\
& & & \quad \forall \ (\alphab_t, \betab_t, \gamma_t) \in \Dcal_t, \ t \in \{1,\dots,T\}, \label{prob:TEOMIO_benders_main} \\
& & & \text{constraints~\eqref{prob:TEOMIO_categorical} - \eqref{prob:TEOMIO_xbinary}}. 
\end{alignat}
\label{prob:TEOMIO_benders}%
\end{subequations}
We can now solve \problemeqref{prob:TEOMIO_benders} using constraint generation. In particular, we start with constraint~\eqref{prob:TEOMIO_benders_main} enforced for a subset of dual solutions $\bar{\Dcal}_t \subseteq \Dcal_t$, and solve \problemeqref{prob:TEOMIO_benders}. This will yield a candidate integer solution $\xb$. We then solve problem~\eqref{prob:SubDual} for each tree $t$ to determine if there exists a solution $(\alphab_t, \betab_t, \gamma_t) \in \Dcal_t$ for which constraint~\eqref{prob:TEOMIO_benders_main} is violated. If so, we add the constraint and solve the problem again. Otherwise, if no such $(\alphab_t, \betab_t, \gamma_t)$ is found for any tree, then the current $\xb$ solution is optimal. 

In the above constraint generation scheme, the key step is to find a dual subproblem solution $(\alphab_t, \betab_t, \gamma_t)$ whose constraint~\eqref{prob:TEOMIO_benders_main} is violated. With this motivation, we now prove a result on the structure of an optimal solution to the dual problem~\eqref{prob:SubDual} (see Section~\ref{appendix:proof_SubDualBinaryOptimal} for the proof).
\begin{proposition}
Let $\xb \in \{0,1\}^{\sum_{i=1}^n K_i}$ be a feasible solution of problem~\eqref{prob:TEOMIO_xonly}. Let $\ell^* = \getLeaf(\xb,t)$. An optimal solution of dual subproblem~\eqref{prob:SubDual} is then given as follows:
\begin{align*}
& \alpha_{t,s} = \left\{ \begin{array}{ll} \max \left\{ \displaystyle \max_{\ell \in \leftleaves(s)} (p_{t,\ell} - p_{t,\ell^*}), 0 \right\} & \text{if}\ s \in \rightsplits(\ell^*), \\ 0 & \text{otherwise},  \end{array}   \right. \\
& \beta_{t,s} = \left\{ \begin{array}{ll} \max \left\{ \displaystyle \max_{\ell \in \rightleaves(s)} (p_{t,\ell} - p_{t,\ell^*}), 0 \right\} & \text{if}\ s \in \leftsplits(\ell^*), \\
 0 & \text{otherwise}, \end{array} \right. \\
& \gamma_t = p_{t,\ell^*}.
\end{align*}\\[-2em]
\label{prop:SubDualBinaryOptimal}
\end{proposition}
The value of Proposition~\ref{prop:SubDualBinaryOptimal} is that we can check for violated constraints in problem~\eqref{prob:TEOMIO_benders} through a simple calculation, without invoking an LO solver. 

We remark on a couple of important aspects of how Proposition~\ref{prop:SubDualBinaryOptimal} is used within our constraint generation approach. First, Proposition~\ref{prop:SubDualBinaryOptimal} requires that the candidate solution $\xb$ to the master problem~\eqref{prob:TEOMIO_benders} be an integer solution. Thus, our constraint generation procedure can only be used to generate Benders cuts for candidate integer solutions of \problemeqref{prob:TEOMIO_benders}. It cannot be used to generate cuts for fractional solutions, such as those arising from the LO relaxation of \problemeqref{prob:TEOMIO_benders}.

Second, we solve \problemeqref{prob:TEOMIO_benders} using branch-and-bound, with none of the Benders constraints~\eqref{prob:TEOMIO_benders_main} enforced. At each integer solution $\xb$ encountered in the branch-and-bound tree, we check whether it satisfies constraint~\eqref{prob:TEOMIO_benders_main} for each tree $t$. This is accomplished by solving the subproblem~\eqref{prob:SubDual} using Proposition~\ref{prop:SubDualBinaryOptimal} to produce a dual solution $(\alphab_t, \betab_t, \gamma_t)$ and checking whether constraint~\eqref{prob:TEOMIO_benders_main} is violated at that dual solution. This procedure is executed for each tree $t$, and any violated constraints are added to every node in the current branch-and-bound tree. (Note that each integer solution encountered in the branch-and-bound process may result in more than one constraint being added if constraint~\eqref{prob:TEOMIO_benders_main} is violated for more than one tree in the ensemble.) In our numerical experiments, we implement this solution approach using \emph{lazy constraint generation}, a computational paradigm implemented in modern solvers such as Gurobi that allows for constraints to be added as needed in a single branch-and-bound tree \citep{gurobi}.

\subsection{Split constraint generation}
\label{subsec:solutionmethods_primalconstraintgeneration}

Recall from Section~\ref{subsec:model_depthapproximation} that when there is a large number of trees and each tree is deep, the total number of splits will be large, and the number of left and right split constraints will be large. However, for a given encoding $\xb$, observe that we do not need all of the left and right split constraints in order for $\yb$ to be completely determined by $\xb$. As an example, suppose for a tree $t$ that $s$ is the root split, and $\sum_{j \in \cond(s)} x_{\var(s),j} = 1$ (\ie, we take the left branch of the root split). In this case, the right split constraint~\eqref{prob:TEOMIO_right} will force all $y_{t,\ell}$ to zero for $\ell \in \rightleaves(s)$. It is clear that in this case, it is not necessary to include any left or right split constraint for any split node $s'$ that is to the right of split $s$, because all of the $y_{t,\ell}$ values that could be affected by those constraints are already fixed to zero.

This suggests an alternate avenue to solving \problemeqref{prob:TEOMIO}, based on iteratively generating the left and right split constraints. Rather than attempting to solve the full problem~\eqref{prob:TEOMIO} with all of the left and right split constraints included in the model, start with a subset of left split constraints and a subset of right split constraints, and all of the decision variables ($x_{i,j}$ for all $i, j$, and $y_{t,\ell}$ for all $t, \ell$) of the full problem. This leads to a relaxation of \problemeqref{prob:TEOMIO}, which we then solve. For the resulting solution $(\xb,\yb)$, determine whether there exist any tree-split pairs $(t,s)$ for which the left split constraint~\eqref{prob:TEOMIO_left} or right split constraint~\eqref{prob:TEOMIO_right} are violated. If a violated left or right split constraint is found, add the corresponding left or right constraint to the formulation, and solve it again. Repeat the procedure until no violated constraints are found, at which point we terminate with the current solution as the optimal solution. 

The key implementation question in such a proposal is: how do we efficiently determine violated constraints? The answer to this question is given in the following proposition. 
\begin{proposition}
Let $(\xb,\yb) \in \{0,1\}^{\sum_{i=1}^n K_i} \times \mathbb{R}^{ \sum_{t=1}^T | \leaves(t) | }$ be a candidate solution to \problemeqref{prob:TEOMIO} that satisfies constraints~\eqref{prob:TEOMIO_ysumtoone} and constraints \eqref{prob:TEOMIO_categorical} to \eqref{prob:TEOMIO_ycontinuous}. Let $t \in \{1,\dots,T\}$. The solution $(\xb,\yb)$ satisfies constraints~\eqref{prob:TEOMIO_left} and \eqref{prob:TEOMIO_right} for all $s \in \splits(t)$ if and only if it satisfies constraint~\eqref{prob:TEOMIO_left} for $s \in \rightsplits(\ell^*)$ and constraint~\eqref{prob:TEOMIO_right} for $s \in \leftsplits(\ell^*)$, where $\ell^* = \getLeaf(\xb,t)$. 
\label{prop:TreeTraversalFeasibility}
\end{proposition} 
Proposition~\ref{prop:TreeTraversalFeasibility} (see Section~\ref{appendix:proof_TreeTraversalFeasibility} for the proof) states that, to check whether solution $(\xb,\yb)$ satisfies the split constraints for tree $t$, it is only necessary to check the split constraints for those splits that are traversed when the observation encoded by $\xb$ is mapped to a leaf by the action of $\getLeaf$. This is a simple but extremely useful result, because it implies that we can check for violated constraints simply by traversing the tree, in the same way that we do to find the leaf of $\xb$. %

Algorithm~\ref{algorithm:ConstraintGeneration} provides the pseudocode of this procedure. This algorithm involves taking the observation encoded by $\xb$ and walking it down tree $t$, following the splits along the way. For each split we encounter, we determine whether we should proceed to the left child ($\sum_{j \in \cond(s)} x_{\var(s),j} = 1$) or to the right child ($\sum_{j \in \cond(s)} x_{\var(s), j} = 0$). If we are going to the left ($s \in \leftsplits(\ell^*)$ or equivalently, $\ell^* \in \leftleaves(s)$), then we check that $y_{t,\ell}$ is zero for all the leaves to the right of split $s$ (constraint~\eqref{prob:TEOMIO_right}). If we are going to the right ($s \in \rightsplits(\ell^*)$ or equivalently, $\ell^* \in \rightleaves(s)$), then we check that $y_{t,\ell}$ is zero for all the leaves to the left of split $s$ (constraint~\eqref{prob:TEOMIO_left}). In words, we are traversing the tree as we would to make a prediction, and we are simply checking that there is no positive $y_{t,\ell}$ that is on the ``wrong'' side of any left or right split that we take. If we reach a leaf node, we can conclude that the current solution $(\xb,\yb)$ does not violate any of the split constraints of tree $t$. 

\begin{algorithm}
\small
\caption{Procedure for verifying feasibility of candidate solution $(\xb, \yb)$. }
\label{algorithm:ConstraintGeneration}
\begin{algorithmic}
\REQUIRE Candidate solution $(\xb, \yb)$, satisfying constraint~\eqref{prob:TEOMIO_ysumtoone}, \eqref{prob:TEOMIO_categorical} - \eqref{prob:TEOMIO_ycontinuous}
\STATE Initialize $\nu \leftarrow \root(t)$ 
\WHILE{$\nu \notin \leaves(t)$}
	\IF {$\sum_{j \in \cond(\nu)} x_{\var(\nu),j} =1$} %
		\IF{ $\sum_{\ell \in \rightleaves(\nu)} y_{t,\ell} > 1 - \sum_{j \in \cond(\nu)} x_{\var(\nu),j} $ }
			\RETURN Violated right constraint~\eqref{prob:TEOMIO_right} at split $\nu$
		\ELSE
			\STATE Set $\nu \leftarrow \leftchild(\nu)$
		\ENDIF
	\ELSE %
		\IF{ $\sum_{\ell \in \leftleaves(\nu)} y_{t,\ell} > \sum_{j \in \cond(\nu)} x_{\var(\nu),j} $ } 
			\RETURN Violated left constraint~\eqref{prob:TEOMIO_left} at split $\nu$ 
		\ELSE
			\STATE Set $\nu \leftarrow \rightchild(\nu)$
		\ENDIF
	\ENDIF
\ENDWHILE
\end{algorithmic}
\end{algorithm}

With regard to running time, Algorithm~\ref{algorithm:ConstraintGeneration}, like \getLeaf, has a complexity of $O(d_{\max})$: for a given tree, Algorithm~\ref{algorithm:ConstraintGeneration} will check at most $d_{\max}$ left/right split constraints in the worst case. In contrast, the total number of split constraints for a given tree $t$ in \problemeqref{prob:TEOMIO} could be much larger than $O(d_{\max})$. For example, in a complete binary tree where every leaf is at depth $d_{\max}+1$, the total number of split constraints will be $O(2^{d_{\max}})$. Algorithm~\ref{algorithm:ConstraintGeneration} thus provides a computationally tractable path to solving instances of \problemeqref{prob:TEOMIO} with deep trees and large numbers of split constraints.

As with our Benders approach, we comment on two important aspects of this solution approach. First, the separation procedure provided in Algorithm~\ref{algorithm:ConstraintGeneration} can only be used for candidate \emph{integer} solutions to \problemeqref{prob:TEOMIO}; it cannot be used to separate fractional solutions, and thus cannot be used to solve the LO relaxation of \problemeqref{prob:TEOMIO}. Extending Algorithm~\ref{algorithm:ConstraintGeneration} and Proposition~\ref{prop:TreeTraversalFeasibility} to fractional solutions $(\xb,\yb)$ is an interesting direction for future research.

Second, we employ Algorithm~\ref{algorithm:ConstraintGeneration} in a manner similar to our Benders approach in Section~\ref{subsec:solutionmethods_bendersdecomposition}. In particular, we solve \problemeqref{prob:TEOMIO} using branch-and-bound, without any of the left/right split constraints enforced. Whenever a new candidate integer solution $(\xb,\yb)$ is found in the branch-and-bound process, we execute Algorithm~\ref{algorithm:ConstraintGeneration} to determine if $(\xb,\yb)$ violates any left/right split constraints for any of the $T$ trees; if any violated constraints are found, they are added to every node in the branch-and bound tree. (This may result in more than one violated constraint being added, where each violated constraint will correspond to one tree.) This computational approach is also implemented using lazy constraint generation.

\section{Computational experiments}
\label{sec:computationalexperiments}

In this section, we describe our first set of computational results. Section~\ref{subsec:computationalexperiments_background} provides the background of our experiments. Section~\ref{subsec:computationalexperiments_fullmio} provides initial results on the full MIO formulation~\eqref{prob:TEOMIO}, while Section~\ref{subsec:computationalexperiments_depth} provides results on the depth approximation scheme of Section~\ref{subsec:model_depthapproximation}. Finally, Section~\ref{subsec:computationalexperiments_solutionmethod} compares the Benders and split generation solution methods against directly solving \problemeqref{prob:TEOMIO}.

\subsection{Background}
\label{subsec:computationalexperiments_background}

We test our optimization formulation~\eqref{prob:TEOMIO} and the associated solution methods from Section~\ref{sec:solutionmethods} using tree ensemble models estimated from real data sets, whose details are provided in Table~\ref{table:datasets}. We wish to point out that in these data sets, the independent variables may in reality not be controllable. However, they are still useful in that they furnish us with real tree ensemble models for evaluating our optimization methodology.

\begin{table}
\centering
\small
\begin{tabular}{llrrl} \toprule
& & Num. & Num. & \\
Data set & Source & Vars. & Obs. &   Description \\  \midrule
\texttt{winequalityred}\ *& \cite{cortez2009modeling} & 11 & 1599 & Predict quality of (red) wine  \\
& & & &  \\[0.25em]
\texttt{concrete} **& \cite{yeh1998modeling} & 8 & 1030 & Predict strength of concrete   \\
& &  & &  \\[0.25em]
\texttt{permeability} **& \cite{kansy1998physicochemical} & 1069 & 165  & Predict permeability of compound \\ 
& &  & &   \\[0.25em]
\texttt{solubility} **& \cite{tetko2001estimation}, & 228 & 951 & Predict solubility of compound \\
& \cite{huuskonen2000estimation}   &   & &  \\ \bottomrule %
\end{tabular}

\caption{Summary of real data sets used in numerical experiments. Note: * = accessed via UCI Machine Learning Repository  \citep{Lichman2013};  ** = accessed via \texttt{AppliedPredictiveModeling} package in R \citep{kuhn2014applied}. \label{table:datasets}}
\end{table}

We specifically focus on random forest models. Unless otherwise stated, all random forests are estimated in R using the \texttt{randomForest} package \citep{liaw2002classification}, using the default parameters. All linear and mixed-integer optimization models are formulated in the Julia programming language \citep{bezanson2012julia}, using the \texttt{JuMP} package \citep[Julia for Mathematical Programming; see][]{lubin2015computing}, and solved using Gurobi 6.5 \citep{gurobi}. All experiments were executed on a late 2013 Apple Macbook Pro Retina laptop, with a quad-core 2.6GHz Intel i7 processor and 16GB of memory.

\subsection{Full MIO formulation experiments} %
\label{subsec:computationalexperiments_fullmio}

As part of our first experiment, we consider solving the unconstrained tree ensemble problem for each data set. For each data set, we consider optimizing the default random forest model estimated in R which uses 500 trees (the parameter \texttt{ntree} in \texttt{randomForest} is set to 500). For each data set, we also consider solving the tree ensemble problem using only the first $T$ trees of the complete forest, where $T$ ranges in $\{10, 50, 100, 200\}$. For each data set and each value of $T$, we solve the MIO formulation~\eqref{prob:TEOMIO}, as well as its LO relaxation. 

We compare our MIO formulation against two other approaches:
\begin{enumerate}
\item \textbf{Local search}: We solve the tree ensemble problem~\eqref{prob:TreeEnsembleOpt_abstract} using a local search heuristic. The details of this local search are provided in Section~\ref{appendix:localsearch} of the e-companion; at a high level, it starts from a randomly chosen initial solution and iteratively improves the solution, one independent variable at a time, until a local optimum is reached. The heuristic is repeated from ten starting points, out of which we only retain the best (highest objective value) solution. We test such an approach to establish the value of our MIO-based approach, which obtains a globally optimal solution, as opposed to a locally optimal solution.
\item \textbf{Standard linearization MIO}: We solve the standard linearization MIO~\eqref{prob:TEOLinearized} and its relaxation, in order to obtain a relative sense of the strength of formulation~\eqref{prob:TEOMIO}. Due to this formulation being much harder to solve, we impose a 30 minute time limit on the solution time of the integer formulation. 
\end{enumerate}

We consider several metrics:
\begin{itemize}
\item $N_{Levels}$: the number of levels (\ie, dimension of $\xb$), defined as $N_{Levels} = \sum_{i=1}^n K_i$.\\[-0.8em]
\item $N_{Leaves}$: the number of leaves (\ie, dimension of $\yb$), defined as $N_{Leaves} = \sum_{t=1}^T | \leaves(t) |$. 
\item $\Tcal_{MIO}$: the time (in seconds) to solve our MIO~\eqref{prob:TEOMIO}.
\item $\Tcal_{StdLin,MIO}$: the time (in seconds) to solve the standard linearization MIO~\eqref{prob:TEOLinearized}.
\item $\Tcal_{LS}$: the time (in seconds) to run the local search procedure (value reported is the total for ten starting points).
\item $G_{LS}$: the gap of the local search solution; if $Z_{LS}$ is the objective value of the local search solution and $Z^*$ is the optimal objective value of \problemeqref{prob:TEOMIO}, then
$$ G_{LS} = 100\% \times (Z^* - Z_{LS}) / Z^*.$$
\item $G_{LO}$: the gap of the LO relaxation of \problemeqref{prob:TEOMIO}; if $Z_{LO}$ is the objective value of the LO relaxation and $Z^*$ is the optimal integer objective as before, then 
$$ G_{LO} = 100\% \times (Z_{LO} - Z^*) / Z^*.$$
\item $G_{StdLin,LO}$: the gap of the LO relaxation of the standard linearization MIO~\eqref{prob:TEOLinearized}; if $Z_{StdLin,LO}$ is the optimal value of the relaxation, then
$$ G_{StdLin,LO} = 100\% \times (Z_{StdLin,LO} - Z^*) / Z^*.$$
\item $G_{StdLin,MIO}$: the optimality gap of the standard linearization MIO~\eqref{prob:TEOLinearized}; if $Z_{StdLin,UB}$ and $Z_{StdLin,LB}$ are the best upper and lower bounds, respectively, of problem~\eqref{prob:TEOLinearized} upon termination, then
$$ G_{StdLin,MIO} = 100\% \times (Z_{StdLin,UB} - Z_{StdLin,LB}) / Z_{StdLin,UB}.$$
\end{itemize}
Note that all $\Tcal$ metrics do \emph{not} include the random forest estimation time in R.

Table~\ref{table:unconstrained_results_times} shows solution times and problem size metrics, while Table~\ref{table:unconstrained_results_gaps} shows the gap metrics. From these two tables, we can draw several conclusions. First, the time required to solve \problemeqref{prob:TEOMIO} is very reasonable; in the most extreme case (\texttt{winequalityred}, $T = 500$), \problemeqref{prob:TEOMIO} can be solved to full optimality in about 20 minutes. (Note that no time limit was imposed on \problemeqref{prob:TEOMIO}; all values of $\Tcal_{MIO}$ correspond to the time required to solve \problemeqref{prob:TEOMIO} to full optimality.)  In contrast, the standard linearization problem~\eqref{prob:TEOLinearized} was only solved to full optimality in two out of twenty cases within the 30 minute time limit. In addition, for those instances where the solver reached the time limit, the optimality gap of the final integer solution, $G_{StdLin,MIO}$, was quite poor, ranging from 50 to over 100\%. 

Second, the integrality gap $G_{LO}$ is quite small -- on the order of a few percent in most cases. This suggests that the LO relaxation of \problemeqref{prob:TEOMIO} is quite tight. In contrast, the LO relaxation bound from \problemeqref{prob:TEOLinearized} is weaker than that of \problemeqref{prob:TEOMIO}, as predicted by Proposition~\ref{prop:relaxationstrength}, and strikingly so. The weakness of the relaxation explains why the corresponding integer problem cannot be solved to a low optimality gap within the 30 minute time limit. These results, together with the results above on the MIO solution times and the final optimality gaps of \problemeqref{prob:TEOLinearized}, show the advantages of our formulation~\eqref{prob:TEOMIO} over the standard linearization formulation~\eqref{prob:TEOLinearized}.

Third, although there are many cases where the local search solution performs quite well, there are many where it can be quite suboptimal, even when repeated with ten starting points. Moreover, while the local search time $T_{LS}$ is generally smaller than the MIO time $T_{MIO}$, in some cases it is not substantially lower (for example, \texttt{solubility} for $T = 500$), and the additional time required by the MIO formulation~\eqref{prob:TEOMIO} may therefore be justified for the guarantee of provable optimality. 

\begin{table}
\centering
\begin{tabular}{llllccccccccc} \toprule
Data set & $T$ & $N_{Levels}$ & $N_{Leaves}$ & $\Tcal_{MIO}$ & $\Tcal_{StdLin,MIO}$ &  $\Tcal_{LS}$ \\ \midrule
\texttt{solubility} & 10 & 1253 & 3157 & 0.1 & 215.2 & 0.2 &\\  
 & 50 & 2844 & 15933 & 0.8 & 1800.3 & 1.8 &\\  
 & 100 & 4129 & 31720 & 1.7 & 1801.8 & 8.8 &\\  
 & 200 & 6016 & 63704 & 4.5 & 1877.8 & 33.7 &\\  
 & 500 & 9646 & 159639 & 177.9 & 1800.3 & 147.0 &\\[0.25em]  
\texttt{permeability} & 10 & 2138 & 604 & 0.0 & 122.6 & 1.0 &\\  
 & 50 & 2138 & 3056 & 0.1 & 1800.0 & 1.9 &\\  
 & 100 & 2138 & 6108 & 0.2 & 1800.3 & 3.1 &\\  
 & 200 & 2138 & 12214 & 0.5 & 1800.0 & 6.1 &\\  
 & 500 & 2138 & 30443 & 2.7 & 1800.0 & 19.1 &\\[0.25em]  
\texttt{winequalityred} & 10 & 1370 & 3246 & 1.8 & 1800.1 & 0.0 &\\  
 & 50 & 2490 & 16296 & 18.5 & 1800.1 & 0.6 &\\  
 & 100 & 3000 & 32659 & 51.6 & 1800.1 & 2.5 &\\  
 & 200 & 3495 & 65199 & 216.0 & 1800.2 & 11.4 &\\  
 & 500 & 3981 & 162936 & 1159.7 & 1971.8 & 34.6 &\\[0.25em]  
\texttt{concrete} & 10 & 1924 & 2843 & 0.2 & 1800.8 & 0.1 &\\  
 & 50 & 5614 & 14547 & 22.7 & 1800.1 & 1.3 &\\  
 & 100 & 7851 & 29120 & 67.8 & 1800.1 & 4.3 &\\  
 & 200 & 10459 & 58242 & 183.8 & 1800.2 & 20.2 &\\  
 & 500 & 13988 & 145262 & 846.9 & 1809.4 & 81.6 &\\
 \bottomrule
\end{tabular}
\caption{Solution times for tree ensemble optimization experiment. \label{table:unconstrained_results_times}} \vspace{-1em}
\end{table}

\begin{table}
\centering
\begin{tabular}{llcccc} \toprule
Data set & $T$ & $G_{LO}$ & $G_{StdLin,LO}$ & $G_{StdLin,MIO}$ & $G_{LS}$  \\ \midrule
\texttt{solubility} & 10 & 0.0 & 485.5 & 0.0 & 18.6 \\  
 & 50 & 0.0 & 498.0 & 50.1 & 9.5 \\  
 & 100 & 0.0 & 481.2 & 70.5 & 0.3 \\  
 & 200 & 0.0 & 477.5 & 77.7 & 0.2 \\  
 & 500 & 0.0 & 501.3 & 103.2 & 0.2 \\[0.25em]  
\texttt{permeability} & 10 & 0.0 & 589.5 & 0.0 & 6.1 \\  
 & 50 & 0.0 & 619.4 & 71.9 & 3.5 \\  
 & 100 & 0.0 & 614.1 & 75.0 & 1.8 \\  
 & 200 & 0.0 & 613.0 & 80.0 & 0.1 \\  
 & 500 & 0.0 & 610.4 & 85.9 & 0.0 \\[0.25em]  
\texttt{winequalityred} & 10 & 1.5 & 11581.3 & 89.8 & 1.2 \\  
 & 50 & 3.4 & 11873.6 & 98.3 & 2.3 \\  
 & 100 & 4.3 & 12014.9 & 98.8 & 0.6 \\  
 & 200 & 4.3 & 12000.6 & 99.0 & 1.2 \\  
 & 500 & 4.5 & 12031.8 & 99.2 & 1.4 \\[0.25em]  
\texttt{concrete} & 10 & 0.0 & 6210.6 & 72.5 & 0.0 \\  
 & 50 & 1.8 & 6657.1 & 95.0 & 0.0 \\  
 & 100 & 2.6 & 6706.6 & 98.3 & 0.0 \\  
 & 200 & 1.6 & 6622.2 & 98.5 & 0.0 \\  
 & 500 & 2.2 & 6652.6 & 98.8 & 0.0 \\[0.25em]  
 \bottomrule
\end{tabular}
\caption{Gaps for tree ensemble optimization experiment. \label{table:unconstrained_results_gaps}}
\end{table}

\subsection{Depth $d$ approximation experiments}
\label{subsec:computationalexperiments_depth}

In this section, we investigate the use of the depth $d$ tree problem (formulation~\eqref{prob:TEOMIO_barOmega} with $\bar{\Omega}_d$) for approximating the full depth problem~\eqref{prob:TEOMIO}. We focus on the same data sets as before with $T = 100$. %
We solve \problemeqref{prob:TEOMIO_barOmega} with $\bar{\Omega}_d$ and vary the depth $d$ of the approximation. We consider the upper bound $Z^*_{MIO,d}$ (denoted by ``UB''), the actual value of the solution $Z_d$ (denoted by ``Actual'') and the lower bound $Z^*_{MIO,d} - \sum_{t=1}^T \lambda_t \Delta_t$ (denoted by ``LB'').

\begin{figure}
\centering 
\includegraphics[width=0.6\textwidth]{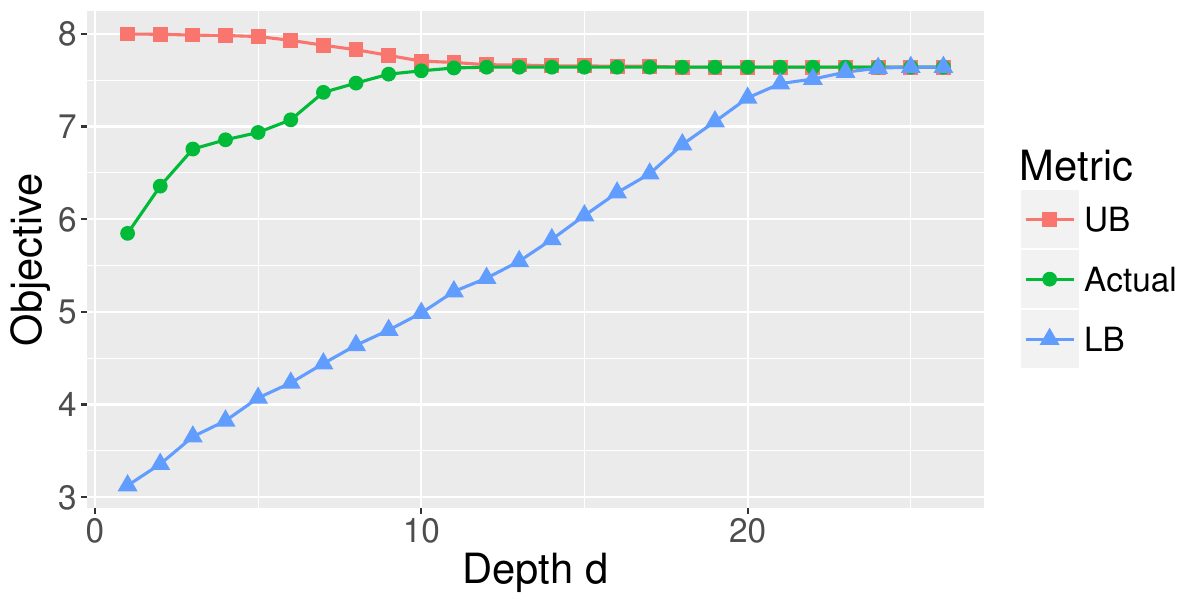} 

\caption{Plot of UB, Actual and LB versus depth for \texttt{winequalityred} with $T = 100$. \label{plot:UBActualLB_vs_depth_wqr_ntree100}}
\end{figure}

\begin{figure}
\centering
\includegraphics[width=0.6\textwidth]{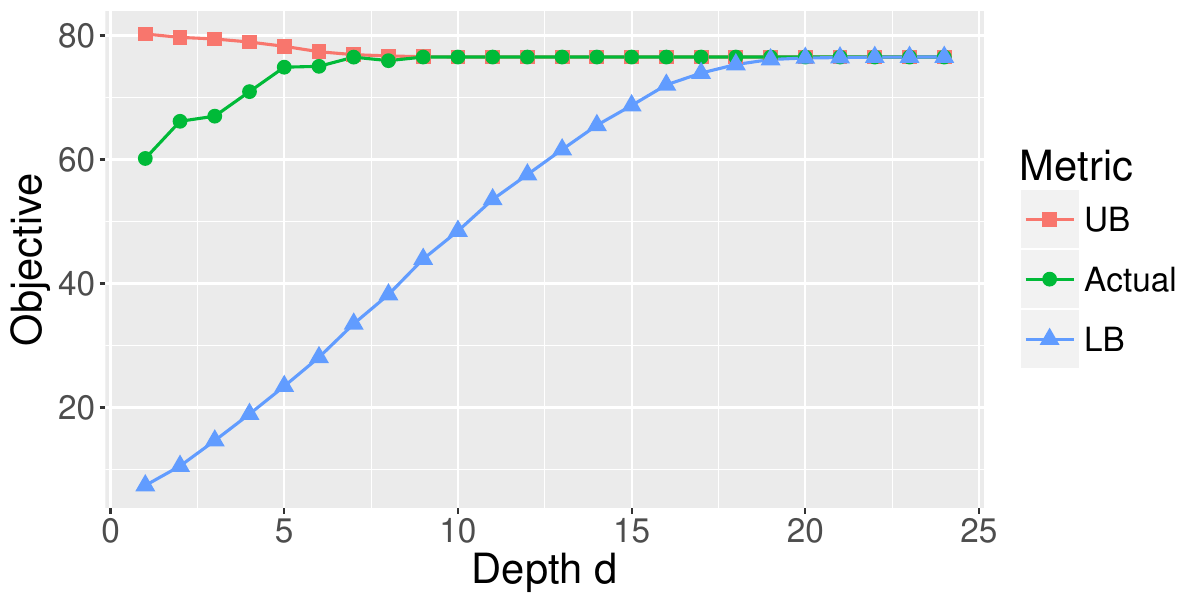} 

\caption{Plot of UB, Actual and LB versus depth for \texttt{concrete} with $T = 100$. \label{plot:UBActualLB_vs_depth_concrete_ntree100}}
\end{figure}

Figures~\ref{plot:UBActualLB_vs_depth_wqr_ntree100} and \ref{plot:UBActualLB_vs_depth_concrete_ntree100} plot the above three metrics for the \texttt{winequalityred} and \texttt{concrete} data sets, respectively. From these plots, we can see that the upper bound is decreasing, while the lower bound and the actual objective are increasing. We can also see that the lower bound is quite loose, and the depth $d$ needs to increase significantly in order for the lower bound to be close to the upper bound. However, even when the depth $d$ is small and the lower bound is loose, the actual objective of the solution produced by the approximation is very good. In the case of \texttt{winequalityred}, the solution is essentially optimal ($Z*_{MIO,d}$ and $Z_d$ are close to or equal) after a depth of $d = 15$ (compared to a maximum depth of 26); for \texttt{concrete}, this occurs for a depth of $d = 9$ (compared to a maximum depth of 24). 

\begin{figure}
\centering 
\includegraphics[width=0.6\textwidth]{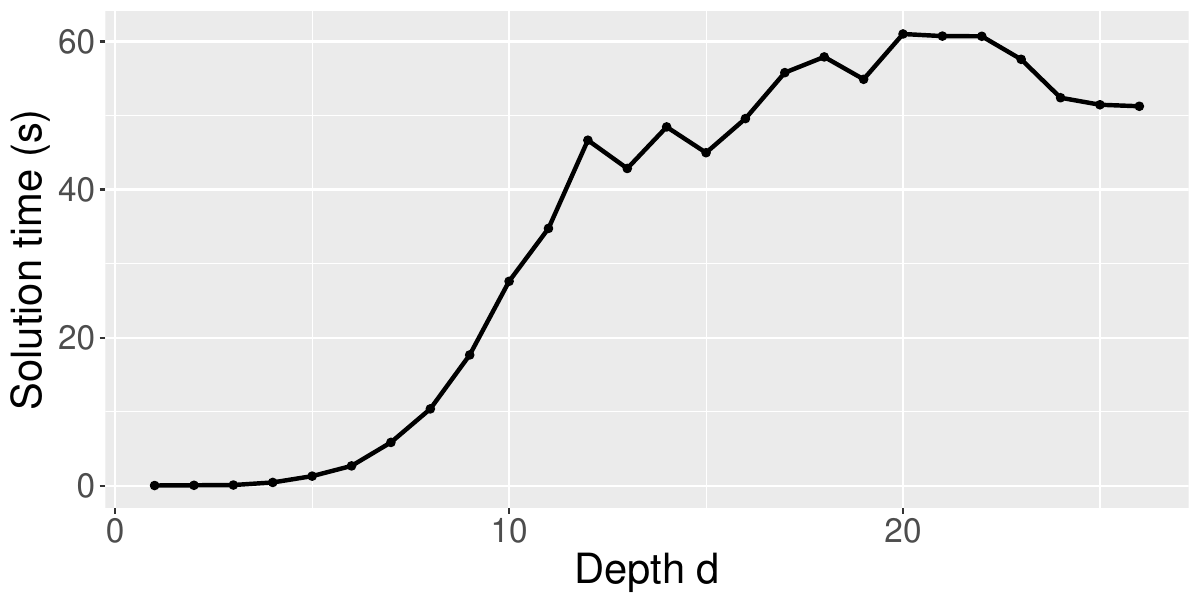} 

\caption{Plot of solution time for problem~\eqref{prob:TEOMIO_barOmega} with $\bar{\Omega}_d$ versus depth $d$ for \texttt{winequalityred} with $T = 100$. \label{plot:time_vs_depth_wqr_ntree100}}
\end{figure}

\begin{figure}
\centering
\includegraphics[width=0.6\textwidth]{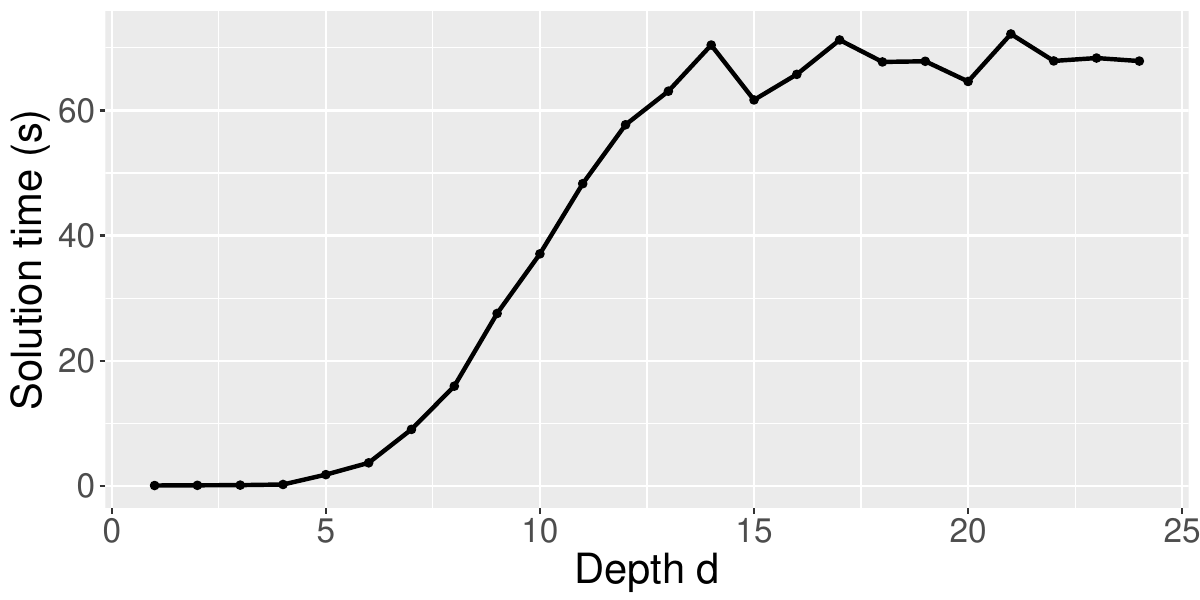} 

\caption{Plot of solution time for problem~\eqref{prob:TEOMIO_barOmega} with $\bar{\Omega}_d$ versus depth $d$ for \texttt{concrete} with $T = 100$. \label{plot:time_vs_depth_concrete_ntree100}}
\end{figure}

To complement these results on the objective values, Figures~\ref{plot:time_vs_depth_wqr_ntree100} and \ref{plot:time_vs_depth_concrete_ntree100} show the computation time of the depth approximation formulation as $d$ varies for the same two data sets. Here we can see that the solution time required to solve the depth approximation formulation initially increases in an exponential fashion as $d$ increases; this is to be expected, because with each additional layer of splits, the number of splits roughly doubles. Interestingly, though, the solution time seems to plateau after a certain depth, and no longer continues to increase. Together with Figures~\ref{plot:UBActualLB_vs_depth_wqr_ntree100} and \ref{plot:UBActualLB_vs_depth_concrete_ntree100}, these plots suggest the potential of the depth approximation approach to obtain near-optimal and optimal solutions with significantly reduced computation time relative to the full depth problem. 

\subsection{Solution method experiments}
\label{subsec:computationalexperiments_solutionmethod}

In this final set of experiments, we evaluate the effectiveness of the two solution methods from Section~\ref{sec:solutionmethods} -- Benders decomposition and split constraint generation -- on solving large instances of \problemeqref{prob:TEOMIO}. We use the same data sets as before with $T = 500$. For each instance, we consider $\Tcal_{Benders}$, $\Tcal_{SplitGen}$ and $\Tcal_{Direct}$, which are the times to solve \problemeqref{prob:TEOMIO} to full optimality using the Benders, split constraint generation and direct solution approaches, respectively. 

Table~\ref{table:solutionmethod_results} shows the results from this comparison. Both approaches can lead to dramatic reductions in the solution time relative to solving \problemeqref{prob:TEOMIO} directly with all split constraints enforced at the start. In the most extreme case (\texttt{concrete}), we observe a reduction from about 800 seconds for the standard solution method to about 32 seconds for split constraint generation and about 37 seconds for the Benders approach -- a reduction in solution time of over 95\%. In some cases, Benders decomposition is slightly faster than split generation;  for example, for \texttt{winequalityred} with $T = 500$, the Benders approach requires just under 11 minutes whereas split generation requires just over 13 minutes. In other cases, split generation is faster (for example, \texttt{solubility} with $T = 500$). 

\begin{table}
\centering
\begin{tabular}{llllllll} \toprule
Data set & $T$ & $N_{Levels}$ & $N_{Leaves}$ & $\Tcal_{Direct}$ & $\Tcal_{Benders}$ & $\Tcal_{SplitGen}$ \\ \midrule
 \texttt{solubility} & 100 & 4129 & 31720 & 1.7 & 1.1 & 0.8\\  
 & 200 & 6016 & 63704 & 4.5 & 2.6 & 0.8\\  
 & 500 & 9646 & 159639 & 177.9 & 28.6 & 13.7\\[0.25em]  
\texttt{permeability} & 100 & 2138 & 6108 & 0.2 & 0.0 & 0.0\\  
 & 200 & 2138 & 12214 & 0.5 & 0.2 & 0.1\\  
 & 500 & 2138 & 30443 & 2.7 & 0.4 & 0.7\\[0.25em]  
\texttt{winequalityred} & 100 & 3000 & 32659 & 51.6 & 41.1 & 56.5\\  
 & 200 & 3495 & 65199 & 216.0 & 152.3 & 57.2\\  
 & 500 & 3981 & 162936 & 1159.7 & 641.8 & 787.8\\[0.25em]  
\texttt{concrete} & 100 & 7851 & 29120 & 67.8 & 8.8 & 8.6\\  
 & 200 & 10459 & 58242 & 183.8 & 12.5 & 15.4\\  
 & 500 & 13988 & 145262 & 846.9 & 37.5 & 32.3\\ %
 \bottomrule 
\end{tabular}
\caption{Results of experiments comparing solution methods. \label{table:solutionmethod_results}}
\end{table}

\section{Case study 1: drug design}
\label{sec:drugdesign}

In this section, we describe our case study in drug design. Section~\ref{subsec:drugdesign_background} provides the background on the problem and the data. Section~\ref{subsec:drugdesign_optimization} shows results on the unconstrained optimization problem, while Section~\ref{subsec:drugdesign_proximity} shows results for when the similarity to existing compounds is constrained to be small.

\subsection{Background}
\label{subsec:drugdesign_background}

For this set of experiments, we use the data sets from \cite{ma2015deep}. These data sets were created for a competition sponsored by Merck and hosted by Kaggle. There are fifteen different data sets. In each data set, each observation corresponds to a compound/molecule. Each data set has a single dependent variable, which is different in each data set and represents some measure of ``activity'' (a property of the molecule or the performance of the molecule for some task). The independent variables in each data set are the so-called ``atom pair'' and ``donor-acceptor pair'' features, which describe the substructure of each molecule (see \citealt{ma2015deep} for further details). The goal of the competition was to develop a model to predict activity using the molecular substructure; such models are known as quantitative structure-activity relationship (QSAR) models.

The optimization problem that we will consider is to find the molecule that maximizes activity as predicted by a random forest model. Our interest in this problem is two-fold. First, this is a problem of significant practical interest, as new drugs are extremely costly to develop. Moreover, these costs are rising: the number of drugs approved by the FDA per billion dollars of pharmaceutical R\&D spending has been decreasing by about 50\% every 10 years since 1950 \citep[a phenomenon known as ``Eroom's Law'' -- Moore's Law backwards; see][]{scannell2012diagnosing}. As a result, there has been growing interest in using analytics to identify promising drug candidates in academia as well as industry (see for example \citealt{atomwise}). We note that random forests are widely used in this domain: the QSAR community was one of the first to adopt them \citep{svetnik2003random} and they have been considered a ``gold standard'' in QSAR modeling \citep{ma2015deep}. %

Second, the problem is of a very large scale. The number of independent variables ranges from about 4000 to just under 10,000, while the number of observations ranges from about 1500 to just over 37,000; in terms of file size, the smallest data set is approximately 15MB, while the largest is just over 700MB. Table~\ref{table:merck_data_summary} summarizes the data sets. Estimating a random forest model using the conventional \texttt{randomForest} package in R on any one of these data sets is a daunting task; to give an example, a single tree on the largest data set requires more than 5 minutes to estimate, which extrapolates to a total computation time of over 8 hours for a realistic forest of 100 trees. 

\begin{table}
\small 
\centering
\begin{tabular}{llrr} \toprule
Data set & Name & Num. Obs. & Num. Variables \\ \midrule
1 & 3A4 & 37241 & 9491 \\ 
2 & CB1 & 8716 & 5877 \\ 
3 & DPP4 & 6148 & 5203 \\ 
4 & HIVINT & 1815 & 4306 \\ 
5 & HIVPROT & 3212 & 6274 \\ 
6 & LOGD & 37388 & 8921 \\ 
7 & METAB & 1569 & 4505 \\ 
8 & NK1 & 9965 & 5803 \\ 
9 & OX1 & 5351 & 4730 \\ 
10 & OX2 & 11151 & 5790 \\ 
11 & PGP & 6399 & 5135 \\ 
12 & PPB & 8651 & 5470 \\ 
13 & RAT\_F & 6105 & 5698 \\ 
14 & TDI & 4165 & 5945 \\ 
15 & THROMBIN & 5059 & 5552 \\  \bottomrule
\end{tabular}
\caption{Summary of drug design data sets (see \citealt{ma2015deep} for further details). \label{table:merck_data_summary}}
\end{table}

\subsection{Unconstrained optimization results}
\label{subsec:drugdesign_optimization}

In our first set of experiments, we proceed as follows. For each data set, we estimate a random forest model to predict the activity variable using all available independent variables. To reduce the computational burden posed by estimating random forest models from such large data sets, we deviate from our previous experiments by using the \texttt{ranger} package in R \citep{wright2017ranger}, which is a faster implementation of the random forest algorithm suited for high dimensional data sets. In addition, we follow \cite{ma2015deep} in restricting the number of trees to 100. For each such random forest model, we solve the (unconstrained) tree ensemble optimization problem~\eqref{prob:TEOMIO} using the Benders approach of Section~\ref{subsec:solutionmethods_bendersdecomposition} and the split generation approach of Section~\ref{subsec:solutionmethods_primalconstraintgeneration}, as well as directly using Gurobi. We impose a time limit of two hours. We also solve each tree ensemble optimization problem using local search with ten repetitions. We consider the following metrics:
\begin{itemize}
\item $G_{Direct}, G_{SplitGen}, G_{Benders}$: the optimality gap of the solution produced by solving the problem directly using Gurobi, solving it using the split generation method and solving it using the Benders methods, respectively. If $Z_{m,LB}$ is the lower bound of a method $m$ and $Z_{m,UB}$ is the upper bound, then $G_{m}$ are defined as
\begin{align*}
G_{m} = 100\% \times (Z_{m,UB} - Z_{m,LB}) / Z_{m,UB}.
\end{align*}
\item $G_{LS}$: the optimality gap of the local search solution, relative to the best split generation solution. If $Z_{LS}$ is the local search objective, it is defined as
\begin{equation*}
G_{LS} = 100\% \times (Z_{SplitGen,LB} - Z_{LS}) / Z_{SplitGen,LB}.
\end{equation*}
\item $\Tcal_{Direct}$ $\Tcal_{SplitGen}, \Tcal_{Benders}$: the time (in seconds) to solve \problemeqref{prob:TEOMIO} directly using Gurobi, using the split generation method and using the Benders method, respectively. (A time that is below 7200 indicates that the problem was solved to full optimality.)
\item $\Tcal_{LS}$: the time (in seconds) to execute the local search procedure. (The time reported is the total of ten repetitions.)
\item $N_{Levels}$ and $N_{Leaves}$: the number of levels and the number of leaves in the ensemble, respectively, defined as in Section~\ref{subsec:computationalexperiments_fullmio}. 
\end{itemize}
Table~\ref{table:merck_solution_methods} displays the results of this experiment. We first discuss the direct solution approach. For this approach, we can see that out of the fifteen data sets, nine were solved to full optimality within one hour, and one more data set was solved within the next hour. For the remaining five data sets (2, 3, 6, 8 and 13), the solver terminated after two hours with very low optimality gaps (four of the data sets having an optimality gap of below 0.5\%, and one with an optimality gap of 2.4\%). The high optimality gap for set 6 is to be expected, as this data set is among the two largest data sets in terms of the number of levels and the total number of leaves (which stems from the number of variables and the number of observations in that data set; see Table~\ref{table:merck_data_summary}). With regard to the split generation approach, we can see that it improves on the direct solution method; with split generation, ten data sets are solved to full optimality within one hour, and another three are solved within the next hour. The greatest improvement is for data set 10, where the direct approach required over an hour, but split generation terminates in just over 20 minutes.

\begin{table}
\centering
\begin{tabular}{crrrrrrrrrr} \toprule
Data set & $N_{Levels}$ & $N_{Leaves}$ & $\Tcal_{Direct}$ & $\Tcal_{SplitGen}$ & $\Tcal_{Benders}$ & $\Tcal_{LS}$ & $G_{Direct}$ & $G_{SplitGen}$ & $G_{Benders}$ & $G_{LS}$ \\ \midrule
 1 & 27145 & 852533 & 151.2 & 97.3 & 7200.0 & 390.8 & 0.00 & 0.00 & 0.07 & 9.26 \\ 
 2 & 16480 & 289800 & 7201.6 & 6533.8 & 7200.0 & 132.6 & 0.13 & 0.00 & 1.30 & 8.64 \\ 
 3 & 13697 & 201265 & 7200.7 & 6252.2 & 7200.0 & 84.9 & 0.22 & 0.00 & 1.41 & 11.58 \\ 
 4 & 11790 & 59552 & 2.9 & 2.1 & 6.3 & 55.4 & 0.00 & 0.00 & 0.00 & 5.92 \\ 
 5 & 16426 & 109378 & 62.5 & 23.9 & 7200.0 & 108.9 & 0.00 & 0.00 & 0.07 & 6.99 \\ 
 6 & 26962 & 1307848 & 7203.5 & 7219.8 & 7200.1 & 409.3 & 2.40 & 2.15 & 12.82 & 23.71 \\ 
 7 & 12523 & 53934 & 16.5 & 12.3 & 1743.0 & 60.8 & 0.00 & 0.00 & 0.00 & 11.17 \\ 
 8 & 17319 & 328705 & 7202.5 & 7200.5 & 7200.1 & 146.9 & 0.37 & 0.27 & 2.67 & 5.22 \\ 
 9 & 12595 & 184841 & 370.2 & 55.0 & 7200.1 & 73.8 & 0.00 & 0.00 & 0.46 & 15.76 \\ 
 10 & 15780 & 379583 & 4101.1 & 1339.5 & 7200.2 & 124.7 & 0.00 & 0.00 & 4.71 & 12.03 \\ 
 11 & 15111 & 217395 & 281.1 & 81.0 & 7200.1 & 94.0 & 0.00 & 0.00 & 0.55 & 12.87 \\ 
 12 & 15737 & 291709 & 32.6 & 40.0 & 7200.0 & 94.4 & 0.00 & 0.00 & 0.02 & 12.17 \\ 
 13 & 17841 & 212926 & 7202.0 & 6731.5 & 7200.0 & 137.3 & 0.43 & 0.00 & 3.73 & 26.38 \\ 
 14 & 16272 & 145476 & 13.1 & 11.5 & 41.8 & 110.9 & 0.00 & 0.00 & 0.00 & 17.10 \\ 
 15 & 14863 & 169638 & 388.1 & 223.1 & 7200.0 & 111.7 & 0.00 & 0.00 & 0.89 & 13.10 \\ \bottomrule
\end{tabular}
\caption{Comparison of split generation and Benders decomposition for drug design data sets. \label{table:merck_solution_methods}}
\end{table}

For the Benders approach, we can see that the performance is quite different. The optimality gap is substantially higher than that achieved by both the direct approach and split generation after two hours. The Benders approach is only able to solve two instances to full optimality within the two hour time limit and in both instances, the split generation approach is able to solve the same instance to full optimality more quickly. %

The last important insight from Table~\ref{table:merck_solution_methods} concerns the performance of the local search procedure. With regard to solution times, we can see that in some cases the total time required for the ten repetitions of the local search \emph{exceeds} the time required for split generation (see data set 1 for example). In addition, and more importantly, the best solution obtained by local search in each data set is highly suboptimal, as evidenced by the high values of $G_{LS}$. In the best case (data set 8), $G_{LS}$ is about 5\%, whereas in the worst case (data set 13), $G_{LS}$ is as high as 26\%. These results suggest that local search is not adequate for this problem: our approaches, which are provably optimal and based on mixed-integer optimization, deliver significantly better solutions. 

\subsection{Controlling proximity}
\label{subsec:drugdesign_proximity}

In the random forest literature, one concept that is useful for analyzing random forest models is that of \emph{proximity}. The proximity of two observations $\Xb, \Xb' \in \Xcal$ is defined as the proportion of trees for which $\Xb$ and $\Xb'$ fall in the same leaf:
\begin{equation*}
\pi(\Xb, \Xb') = \frac{1}{T} \sum_{t=1}^T  \Ibb \{ \ell_t(\Xb) = \ell_t(\Xb') \},
\end{equation*}
where $\ell_t(\Xb)$ is the leaf to which tree $t$ maps the observation $\Xb$. Proximity is valuable to consider because it allows one to use the random forest model to determine how similar two observations are, based on how each tree categorizes those two observations. In \cite{svetnik2003random}, proximity was used as a similarity metric for clustering molecules, and was shown to be competitive with Tanimoto similarity (a different metric, commonly used in QSAR modeling; see \citealt{willett1998chemical}) in producing meaningful molecule clusters. In the remainder of this section, we first analyze the proximity of the solutions from Section~\ref{subsec:drugdesign_optimization}, and verify that these solutions are distinct from the molecules in the data sets. Then, we consider solving \problemeqref{prob:TEOMIO} with constraints on the proximity. By imposing upper bounds on proximity, one can find molecules that maximize predicted activity while limiting how similar they are to existing molecules or equivalently, requiring them to be dissimilar to existing molecules; such a constrained optimization approach may potentially be useful in identifying novel molecules.

Table~\ref{table:merck_proximity} displays, for each data set, the average proximity $\pi_{\avg}$ and maximum proximity $\pi_{\max}$ between the split generation solution and all of the molecules in the data set. We can see that the maximum proximity -- the highest proximity between the solution and any point in the training data -- is in general low. For example, for data set 3, $\pi_{\max}$ is 0.11: for any point in the data set that we consider, at most eleven trees out of the 100 trees in the forest will place both the training point and our solution in the same leaf. In addition, the average proximity is much lower than the maximum proximity, which suggests that for most training set points, the actual proximity to the solution is close to zero (\ie, there is very little similarity between the solution and the training point). %
Most importantly, all solutions have a $\pi_{\max}$ strictly lower than 1, indicating that all fifteen solutions are in fact different from all of the molecules of their respective training data sets.

\begin{table}
\centering
\begin{tabular}{ccc}  \toprule
Data set & $\pi_{\avg}$ & $\pi_{\max}$ \\  \midrule
1 & 0.00006 & 0.45 \\ 
2 & 0.00033 & 0.37 \\ 
3 & 0.00042 & 0.11 \\ 
4 & 0.00157 & 0.71 \\ 
5 & 0.00071 & 0.40 \\ 
6 & 0.00005 & 0.15 \\ 
7 & 0.00115 & 0.43 \\ 
8 & 0.00029 & 0.42 \\ \bottomrule
\end{tabular}
\quad 
\begin{tabular}{ccc}  \toprule
Data set & $\pi_{\avg}$ & $\pi_{\max}$ \\  \midrule
9 & 0.00040 & 0.51 \\ 
10 & 0.00027 & 0.48 \\ 
11 & 0.00028 & 0.36 \\ 
12 & 0.00032 & 0.54 \\ 
13 & 0.00030 & 0.36 \\ 
14 & 0.00053 & 0.41 \\ 
15 & 0.00039 & 0.21 \\  \bottomrule
\end{tabular}

\caption{Average and maximum proximity of split generation solutions for drug design data sets. \label{table:merck_proximity}}

\end{table}

We now turn our attention to solving \problemeqref{prob:TEOMIO} with an added constraint on the maximum proximity of the solution to the training points. Such a constraint is defined as follows: let $\Xb^{(1)}, \dots, \Xb^{(M)}$ be the set of observations in the training data, and define for each observation $m$ the vector $\yb^{(m)} \in \mathbb{R}^{ \sum_{t=1}^T | \leaves(t)|}$ as 
\begin{equation*}
y^{(m)}_{t,\ell} = \Ibb\{ \ell_t(\Xb^{(m)}) = \ell \}.
\end{equation*}
The proximity between the fixed observation $\Xb^{(m)}$ and the solution encoded by $(\xb,\yb)$ can then be written as an affine function of $\yb$: 
\begin{equation*}
\frac{1}{T} \sum_{t=1}^T \sum_{\ell \in \leaves(t)}  y^{(m)}_{t,\ell} \cdot y_{t,\ell}.
\end{equation*}
(Note that $\yb^{(1)}, \dots \yb^{(M)}$ are data, and not decision variables.) We can thus enforce a constraint on the proximity of the solution encoded by $\xb$ to each observation $\Xb^{(1)}, \dots, \Xb^{(M)}$ to be at most $c \in [0,1]$ through the following family of linear constraints on $\yb$:
\begin{equation}
\frac{1}{T} \sum_{t=1}^T \sum_{\ell \in \leaves(t)}  y^{(m)}_{t,\ell} \cdot y_{t,\ell} \leq c, \quad \forall \ m \in \{1,\dots, M\}. \label{eq:maxproximity}
\end{equation}
We solve \problemeqref{prob:TEOMIO} with constraint~\eqref{eq:maxproximity} and vary the parameter $c$ to generate a Pareto efficient frontier of solutions that optimally trade-off their maximum proximity to the training data and their predicted value under the random forest model. We impose a time limit of one hour on each solve of \problemeqref{prob:TEOMIO}. We solve each constrained instance using the split generation approach.

\begin{figure}
\centering
\begin{tabular}{cc}
\includegraphics[width=0.5\textwidth]{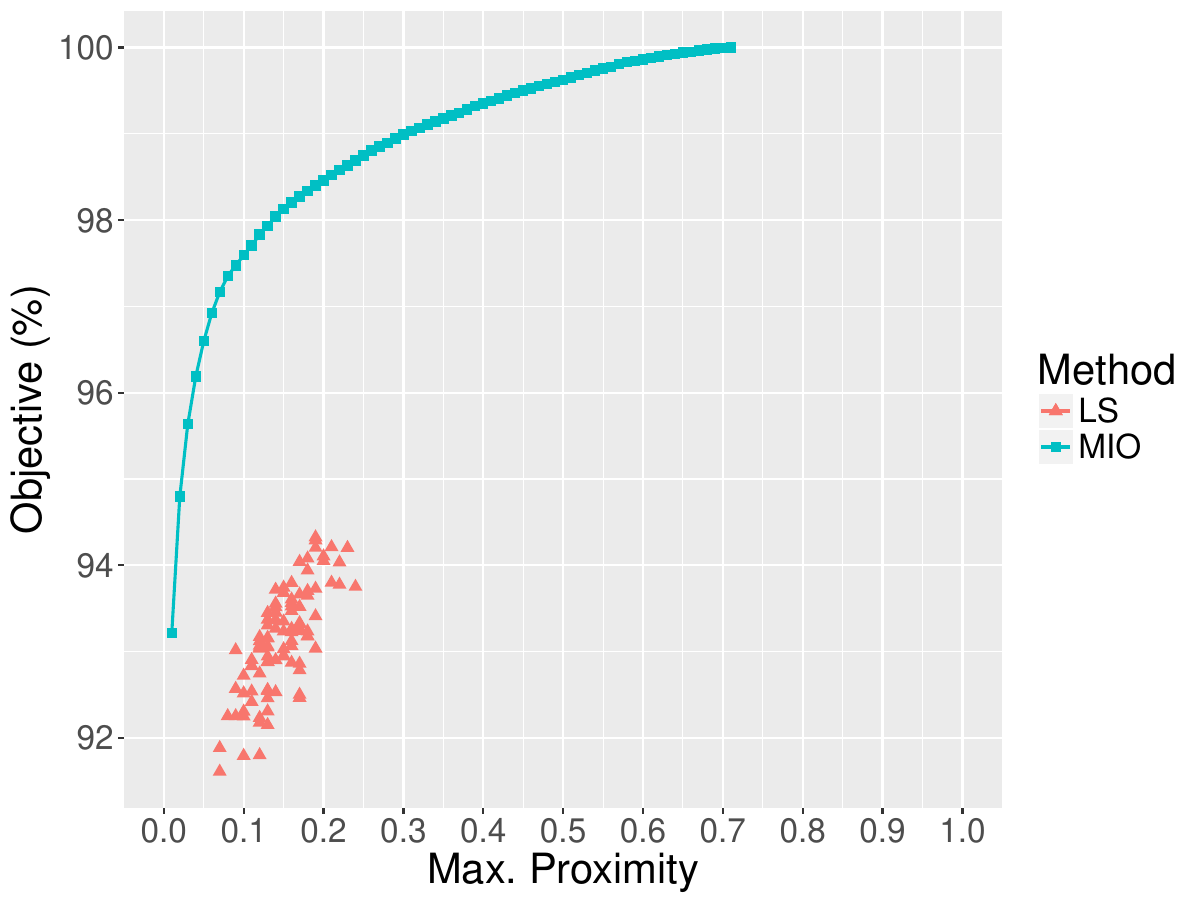} & \includegraphics[width=0.5\textwidth]{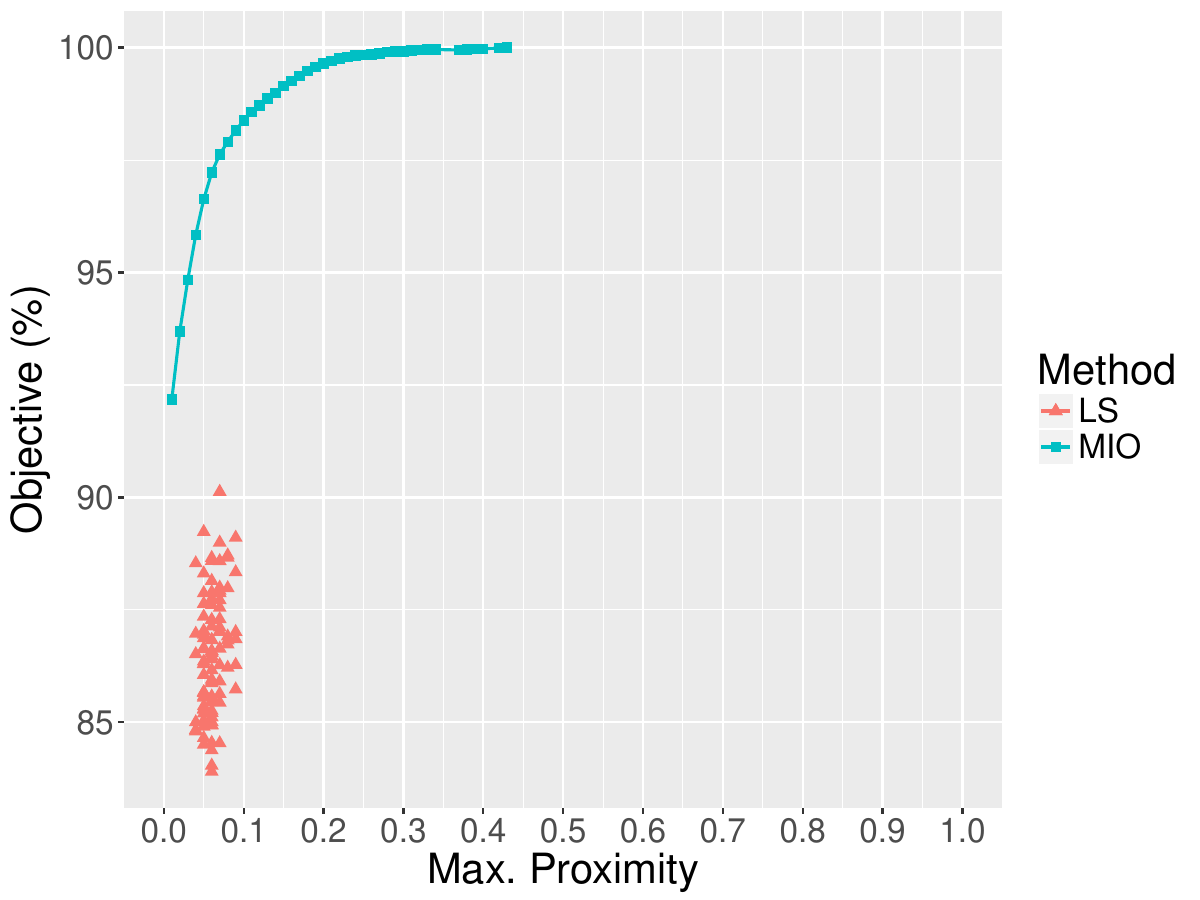}
\end{tabular}
\caption{Plot of objective-maximum proximity frontier for drug design data set 4 (left) and data set 7 (right). \label{plot:proxwhile_4and7}}  %
\end{figure}

To demonstrate, we focus on data sets 4 and 7. Figure~\ref{plot:proxwhile_4and7} shows the proximity-objective value frontier (the points labeled ``MIO'') for data sets 4 and 7.  (Note that the objective value is expressed in terms of the maximum unconstrained objective value, \ie, the objective value attained when the proximity constraint is omitted.) We note that the right-most point in each frontier (maximum proximities of 0.71 for data set 4 and 0.43 for data set 7) corresponds to the original unconstrained solution. As a comparison, we also solve the (unconstrained) problem for each data set using local search with 100 repetitions, and plot the proximity and objective value of each solution from each repetition (the points labeled ``LS''). We remark that both frontiers are approximate because the constrained optimization problem associated with each point is not necessarily solved to full optimality, due to the one hour time limit. However, upon closer examination, only six points out of a total of 112 across both plots did not solve to full optimality, with a maximum gap of only about 0.12\%. Thus, although constraint~\eqref{eq:maxproximity} adds to the problem size, it does not drastically impact our ability to solve \problemeqref{prob:TEOMIO}.

From both of these figures, we obtain several insights. First, in these two cases, we are able to push the proximity to its lowest; for both data sets, we can find solutions with maximum proximities of 0.01 (\ie, one tree out of the ensemble of 100). (A solution with maximum proximity of zero does not exist, because for any tree, the leaf $\ell_t(\Xb)$ that the candidate solution $\Xb$ is mapped to will also be shared by at least one point in the training set.) 
Second, we can can see that the price of dissimilarity is low: as we decrease the maximum proximity, we can still obtain solutions with very good predicted performance. For example, for data set 4, we can see that if we lower the proximity to 0.01, the relative objective value decreases by only about 7\%. Third, although the solutions obtained by local search have smaller maximum proximities than the unconstrained MIO solution, they are highly suboptimal with respect to objective value (best relative objectives for data sets 4 and 7 are roughly 94\% and 90\%, respectively) and are dominated in both maximum proximity and objective value by the MIO solutions. Overall, these results suggest that our MIO formulation can be used to systematically identify promising molecules that have good predicted performance and are sufficiently different from existing molecules in the training data.

\section{Case study 2: customized pricing}

\label{sec:pricing_summary}

In addition to our drug design case study, we also consider the problem of customized pricing. The results of this case study are fully reported in Section~\ref{sec:pricing} of the e-companion. At a high level, we estimate random forest models of profits for eleven orange juice products as a function of their prices and store-level attributes, and then optimize these models to derive prices for each store. We also estimate log-linear and semi-log demand models in a hierarchical Bayesian (HB) framework, and heuristically optimize these models to obtain prices as well. We find that our random forest models can be optimized quickly (on the order of seconds). Moreover, the random forest-optimal prices perform well under the HB models, and are less extreme, in that fewer products are set to their highest or lowest allowable prices across the whole store chain. From a predictive standpoint, our random forest models also achieve significantly higher out-of-sample $R^2$ values than the HB models, which may be of independent interest.

\section{Conclusion}
\label{sec:conclusion}

In this paper, we developed a modern optimization approach to the problem of finding the decision that optimizes the prediction of a tree ensemble model. At the heart of our approach is a mixed-integer optimization formulation that models the action of each tree in the ensemble. We showed that this formulation is better than a general alternate formulation, that one can construct an hierarchy of approximations to the formulation with bounded approximation quality through depth-based truncation and that one can exploit the structure of the formulation to derive efficient solution methods, based on Benders decomposition and split constraint generation. We demonstrated the utility of our approach using real data sets, including two case studies in drug design and customized pricing. 
Given the prevalence of tree ensemble models, we believe that this methodology will become an important asset in the modern business analytics toolbox and is an exciting starting point for future research at the intersection of optimization and machine learning.

\ACKNOWLEDGMENT{The author sincerely thanks the area editor Marina Epelman, the associate editor and the two anonymous referees for their careful reading of the paper and thoughtful comments that have helped to greatly improve the quality of the paper. The author thanks Fernanda Bravo, Vishal Gupta and Auyon Siddiq for helpful conversations. The authors of \cite{ma2015deep} and Merck Research Laboratories are gratefully acknowledged for the data set used in the drug design case study of Section~\ref{sec:drugdesign}. Dominick's Finer Foods and Peter Rossi are gratefully acknowledged for the data set used in the customized pricing case study of Section~\ref{sec:pricing}.
}

\bibliographystyle{plainnat}
\bibliography{teo_literature}

\ECSwitch

\ECHead{Electronic companion for ``Optimization of Tree Ensembles''} %

\section{Proofs}

\subsection{Auxiliary results}

Before proceeding to the proofs, we first state two auxiliary results. 

\begin{lemma}
Fix a tree $t$ and a leaf $\ell \in \leaves(t)$. Then 
\begin{equation}
\{ \ell' \in \leaves(t) \, | \, \ell' \neq \ell \} = \bigcup_{s \in \leftsplits(\ell)} \rightleaves(s) \cup \bigcup_{s \in \rightsplits(\ell)} \leftleaves(s). \label{eq:leafdecomp}
\end{equation}
Furthermore, define the collection $ \Scal = \{ \rightleaves(s) \, | \, s \in \leftsplits(\ell)\} \cup \{ \leftleaves(s) \, | \, s \in \rightsplits(\ell) \}$. Then each pair of distinct sets $A, B \in \Scal$, $A \neq B$, is disjoint.
\label{lemma:leafdecomp}
\end{lemma}
To gain some intuition for the right hand set in equation~\eqref{eq:leafdecomp}, recall that $\leftsplits(\ell)$ is the set of splits for which we follow the left branch in order to reach $\ell$, and $\rightsplits(\ell)$ is the set of splits for which we follow the right branch to reach $\ell$. For each $s \in \leftsplits(\ell)$, $\rightleaves(s)$ is the set of leaves that is on the ``wrong side'' of split $s$ (we take the left branch to reach $\ell$, but each leaf in $\rightleaves(s)$ is only reachable by taking the right branch). Similarly, for each $s \in \rightsplits(\ell)$, $\leftleaves(s)$ is the set of leaves that is on the wrong side of split $s$. The union of all leaves $\ell'$ on the wrong side of each split $s \in \rightsplits(\ell) \cup \leftsplits(\ell)$ covers all leaves except $\ell$. 
\proof{Proof of Lemma~\ref{lemma:leafdecomp}:}
We prove this by establishing two inclusions: %
\begin{align*}
\{ \ell' \in \leaves(t) \, | \, \ell' \neq \ell \}\  \subseteq \ \bigcup_{s \in \leftsplits(\ell)} \rightleaves(s) \cup \bigcup_{s \in \rightsplits(\ell)} \leftleaves(s),\\
\{ \ell' \in \leaves(t) \, | \, \ell' \neq \ell \}\  \supseteq \ \bigcup_{s \in \leftsplits(\ell)} \rightleaves(s) \cup \bigcup_{s \in \rightsplits(\ell)} \leftleaves(s).
\end{align*}
For the first inclusion, let us fix $\ell' \in \leaves(t)$ such that $\ell' \neq \ell$. Then there must exist a split $\bar{s} \in \splits(t)$ such that either $\ell' \in \leftleaves(\bar{s})$, $\ell \in \rightleaves(\bar{s})$, or $\ell' \in \rightleaves(\bar{s})$, $\ell \in \leftleaves(\bar{s})$. In the former case, we have that $\bar{s} \in \rightsplits(\ell)$, so that 
$$ \ell' \in \leftleaves(\bar{s}) \subseteq   \bigcup_{s \in \leftsplits(\ell)} \rightleaves(s) \cup \bigcup_{s \in \rightsplits(\ell)} \leftleaves(s).$$
In the latter case, we have that $\bar{s} \in \leftsplits(\ell)$, so that 
$$ \ell' \in \rightleaves(\bar{s}) \subseteq \bigcup_{s \in \leftsplits(\ell)} \rightleaves(s) \cup \bigcup_{s \in \rightsplits(\ell)} \leftleaves(s).$$
This proves the first inclusion.

For the second inclusion, we will argue the contrapositive. We have that 
\begin{align*}
\left[ \bigcup_{s \in \leftsplits(\ell)} \rightleaves(s) \cup \bigcup_{s \in \rightsplits(\ell)} \leftleaves(s) \right]^C & = \bigcap_{s \in \leftsplits(\ell)} (\rightleaves(s))^C \cap \bigcap_{s \in \rightsplits(\ell)} (\leftleaves(s))^C \\
& \supseteq \bigcap_{s \in \leftsplits(\ell)} \leftleaves(s) \cap \bigcap_{s \in \rightsplits(\ell)} \rightleaves(s) \\
& \supseteq \{ \ell \}
\end{align*}
where the first step follows by De Morgan's law; the second follows by the fact that for any split $s \in \splits(t)$, $\leftleaves(s)$ and $\rightleaves(s)$ are disjoint; and the last by the definition of $\leftsplits(\ell)$ and $\rightsplits(\ell)$. This proves the second inclusion, and thus proves the equivalence.

Finally, to show that the sets in $\Scal$ are pairwise disjoint, order the splits in $\leftsplits(\ell) \cup \rightsplits(\ell)$ according to their depth:
$$ \leftsplits(\ell) \cup \rightsplits(\ell) = \{s_1, s_2, \dots, s_K\},$$
where $K$ is the total number of splits in $\leftsplits(\ell) \cup \rightsplits(\ell)$. Let us also define the sequence of sets $A_1, A_2, \dots, A_K$ as 
$$ A_i = \left\{ \begin{array}{ll} \leftleaves(s_i) & \text{if}\ s_i \in \leftsplits(\ell), \\
\rightleaves(s_i) & \text{if}\ s_i \in \rightsplits(\ell), \end{array} \right. $$
and the sequence of sets $B_1, B_2, \dots, B_K$ as
$$ B_i = \left\{ \begin{array}{ll} \rightleaves(s_i) & \text{if}\ s_i \in \leftsplits(\ell), \\
\leftleaves(s_i) & \text{if}\ s_i \in \rightsplits(\ell). \end{array} \right. $$
We need to show that the collection $\{B_1,\dots, B_K\}$ (this is the collection $\Scal$) is disjoint. Observe that by the definition of $\leftsplits(\ell)$ and $\rightsplits(\ell)$, the sets $A_1, A_2, \dots, A_K$ form a nested sequence, \ie, 
$$ A_1 \supseteq A_2 \supseteq \dots \supseteq A_K.$$
Notice also that for each $i \in \{2,\dots,K\}$, 
$$ B_i \subseteq A_{i-1},$$
and for each $i \in \{1,\dots, K\}$,
$$ B_i \cap A_i = \emptyset.$$
It therefore follows that given $i, j \in \{1,\dots, K\}$ with $i < j$, that 
$$ B_i \cap B_j \subseteq B_i \cap A_i = \emptyset,$$
which establishes that $\{B_1,\dots, B_K\}$ are pairwise disjoint. \Halmos
\endproof

In addition to Lemma~\ref{lemma:leafdecomp}, it will also be useful to state an analogous lemma for splits. With a slight abuse of notation, let us define $\leftsplits(s)$ for a split $s \in \splits(t)$ as the sets of splits $s'$ such that $s$ is on the left subtree of $s'$; similarly, we define $\rightsplits(s)$ for a split $s \in \splits(t)$ as the set of splits $s'$ such that $s$ is on the right subtree of $s'$. We then have the following lemma; the proof follows along similar lines to Lemma~\ref{lemma:leafdecomp} and is omitted.

\begin{lemma}
For a given tree $t$, let $s \in \splits(t)$. We then have 
\begin{equation*}
[ \leftleaves(s) \cup \rightleaves(s) ]^C = \bigcup_{s' \in \rightsplits(s)} \leftleaves(s') \cup \bigcup_{s' \in \leftsplits(s)} \rightleaves(s').
\end{equation*}
\label{lemma:leafdecomp_splits}
\end{lemma}

\subsection{Proof of Proposition~\ref{prop:TEO_NPComplete}}
\label{appendix:proof_TEO_NPComplete}

To prove that problem~\eqref{prob:TreeEnsembleOpt_abstract} is NP-Hard, we will show that it can be used to solve the minimum vertex cover problem. An instance of the minimum vertex cover problem is defined by a graph $(V,E)$, where $V$ is a set of vertices and $E$ is the set of edges between these vertices. The minimum vertex cover problem is to find the smallest set of vertices $S$ from $V$ such that each edge in $E$ is incident to at least one vertex from  $S$. 

\begin{figure}
\begin{center}
\includegraphics[valign=t]{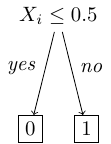}
\qquad
\includegraphics[valign=t] {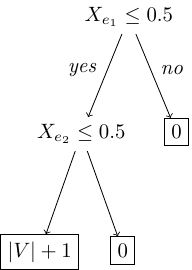}
\end{center}
\caption{Type 1 tree (left) and type 2 tree (right) for vertex cover reduction. \label{figure:trees_NPHard} }
\end{figure}

We now show how to cast this problem as a tree ensemble optimization problem. For convenience, let us index the vertices from $1$ to $|V|$ and the edges from $1$ to $|E|$. Also, let $e_1$ and $e_2$ be the nodes to which edge $e \in E$ is incident to. Suppose that our independent variable $\Xb$ is given by $\Xb = (X_1, \dots, X_{|V|})$, where $X_i$ is a numeric variable that is 1 or 0. The tree ensemble we will consider will consist of the following two types of trees:
\begin{enumerate}
\item \textbf{Type 1}: Trees $1$ to $|V|$, where for $i \in \{1,\dots,|V|\}$,
\begin{equation*}
f_i( \Xb) = \left\{ \begin{array}{ll} 0 & \text{if} \ X_i \leq 0.5, \\
1 & \text{if} \ X_i > 0.5. \\ 
\end{array} \right.
\end{equation*}
\item \textbf{Type 2}: Trees $|V|+1$ to $|V| + |E|$, where for each edge $e \in E$, 
\begin{equation*}
f_{|V| + e} ( \Xb) = \left\{ \begin{array}{ll} 0 & \text{if} \ X_{e_1} > 0.5, \\
0 & \text{if} \ X_{e_1} \leq 0.5, \ X_{e_2} > 0.5,  \\
+|V|+1 & \text{if} \ X_{e_1} \leq 0.5, \ X_{e_2} \leq 0.5. \end{array} \right.
\end{equation*}
\end{enumerate}
These two types of trees are visualized in Figure~\ref{figure:trees_NPHard}.

We let the weight $\lambda_t$ of each tree $t$ be $-1$. The corresponding tree ensemble optimization problem is
\begin{equation}
\underset{\Xb \in \{0,1\}^{|V|}}{\text{maximize}} \ - \sum_{t=1}^{|V|} f_t( \Xb) - \sum_{t=|V|+1}^{|V| + |E|} f_t(\Xb). \label{prob:VertexCover_as_TEO}
\end{equation}
The above problem is identical to the minimum vertex cover problem. In particular, the independent variable $\Xb$ encodes the cover; $X_i = 1$ indicates that vertex $i$ is part of the cover, and $X_i = 0$ indicates that vertex $i$ is not in the cover. The type 1 trees count the size of the cover, while the type 2 trees penalize the solution if an edge is not covered by the set. More precisely, to understand the role of the type 2 trees, observe that:
\begin{itemize}
\item If the set of vertices encoded by $\Xb$ is a feasible cover, then $\sum_{t=|V|+1}^{|V| + |E|} f_t(\Xb) = 0$, and the objective only consists $- \sum_{t=1}^{|V|} f_t( \Xb)$, which counts the number of vertices in the set encoded by $\Xb$. 
\item If the set of vertices encoded by $\Xb$ is not a feasible cover, then $f_t(\Xb) = |V|+1$ for at least one $t \in \{|V| + 1, \dots, |V| + |E|\}$, and therefore the objective satisfies the bound
\begin{equation}
- \sum_{t=1}^{|V|} f_t( \Xb) - \sum_{t=|V|+1}^{|V| + |E|} f_t(\Xb) \leq - (|V| + 1).  \label{eq:vertexcoverinfeasiblebound}
\end{equation}
\end{itemize}
Observe that the bound in inequality~\eqref{eq:vertexcoverinfeasiblebound} is strictly worse than selecting all of the vertices, that is, setting $X_1 = X_2 = \dots = X_{|V|} = 1$; using all of the vertices corresponds to an objective value of $- |V|$. Therefore, at optimality, the set of vertices encoded by $\Xb$ must be a feasible cover. As stated above, the objective value of $\Xb$ when it corresponds to a feasible cover reduces to 
\begin{equation*}
- \sum_{t=1}^{|V|} f_t(\Xb) = - | \{ i \, | \, X_i = 1 \} |,
\end{equation*}
which is (the negative of) the size of the set of vertices encoded by $\Xb$. Maximizing this quantity is equivalent to minimizing its negative, which is the same as minimizing the size of the set of vertices that covers $E$. Therefore, solving \problemeqref{prob:VertexCover_as_TEO} is equivalent to solving the minimum vertex cover problem for $(V,E)$. 

Since the minimum vertex cover problem is NP-Complete \citep{garey2002computers}, it follows that the tree ensemble optimization problem~\eqref{prob:TreeEnsembleOpt_abstract} is NP-Hard. \Halmos

\subsection{Proof of Proposition~\ref{prop:relaxationstrength}} 
\label{appendix:proof_relaxationstrength}

To prove the proposition, we will show that any optimal solution $(\xb,\yb)$ of the relaxation of \problemeqref{prob:TEOMIO} is a feasible solution of the relaxation of the standard linearization problem~\eqref{prob:TEOLinearized}. Since the objective functions of problems~\eqref{prob:TEOMIO} and \eqref{prob:TEOLinearized} are the same, it will follow that the objective of $(\xb,\yb)$, which is $Z^*_{LO}$, is less than or equal to $Z^*_{LO,StdLin}$, which is the optimal objective value of the relaxation of \problemeqref{prob:TEOLinearized}. 

Let $(\xb,\yb)$ be an optimal solution of the relaxation of problem~\eqref{prob:TEOMIO}. To show it is feasible for the relaxation of \problemeqref{prob:TEOLinearized}, we need to show that it satisfies the constraints of that formulation. We only need to show that constraints~\eqref{prob:TEOLinearized_left} -- \eqref{prob:TEOLinearized_lower} are satisfied, since the other constraints of \problemeqref{prob:TEOLinearized} are the same as in \problemeqref{prob:TEOMIO}. 

To verify constraint~\eqref{prob:TEOLinearized_left}, observe that $(\xb,\yb)$ satisfies constraint~\eqref{prob:TEOMIO_left}. For any $t \in \{1,\dots,T\}$, $\ell \in \leaves(t)$ and $s \in \leftsplits(\ell)$, we have
\begin{align*}
\sum_{j \in \cond(s)} x_{\var(s),j} & \geq \sum_{\ell' \in \leftleaves(s) } y_{t,\ell'} \\ %
& \geq y_{t,\ell} 
\end{align*}
where the first inequality is exactly constraint~\eqref{prob:TEOMIO_left}, and the second inequality follows because $\ell \in \leftleaves(s)$ (this is because $s \in \leftsplits(\ell)$) and all $y_{t,\ell'}$'s are nonnegative (by constraint~\eqref{prob:TEOMIO_ycontinuous}). Therefore, $(\xb,\yb)$ satisfies constraint~\eqref{prob:TEOLinearized_left}. Similar reasoning can be used to establish constraint~\eqref{prob:TEOLinearized_right}.

To verify constraint~\eqref{prob:TEOLinearized_lower}, observe that $(\xb,\yb)$ satisfies 
\begin{equation*}
\sum_{\ell' \in \leaves(t)} y_{t,\ell'} \geq 1,
\end{equation*}
for any tree $t$, by virtue of constraint~\eqref{prob:TEOMIO_ysumtoone}. Fix a tree $t$ and a leaf $\ell \in \leaves(t)$, and re-arrange the above to obtain
\begin{equation}
y_{t,\ell} \geq 1 - \sum_{\ell' \neq \ell} y_{t,\ell'}. \label{eq:ytellgreaterthannonell}
\end{equation}
We then have
{\allowdisplaybreaks
\begin{align*}
y_{t,\ell} & \geq 1 - \sum_{\ell' \neq \ell} y_{t,\ell'} \\
& = 1 - \sum_{s \in \leftsplits(\ell)} \sum_{\ell' \in \rightleaves(s)} y_{t,\ell'}  - \sum_{s \in \rightsplits(\ell)} \sum_{\ell' \in \leftleaves(s)} y_{t,\ell'} \\
& \geq 1 - \sum_{s \in \leftsplits(\ell)} ( 1 - \sum_{j \in \cond(s)} x_{\var(s),j}) - \sum_{s \in \rightsplits(\ell)} ( \sum_{j \in \cond(s)} x_{\var(s),j}) \\
& = 1 - \sum_{s \in \leftsplits(\ell)} ( 1 -   \sum_{j \in \cond(s)} x_{\var(s),j} ) - \sum_{s \in \rightsplits(\ell)} ( 1 - (1 -  \sum_{j \in \cond(s)} x_{\var(s),j}) ) \\
& = \sum_{s \in \leftsplits(\ell)} \sum_{j \in \cond(s)} x_{\var(s),j}  + \sum_{s \in \rightsplits(\ell)} (1 - \sum_{j \in \cond(s)} x_{\var(s),j}) - (|\leftsplits(\ell)| + |\rightsplits(\ell)| - 1) 
\end{align*}}%
where the first inequality is just inequality~\eqref{eq:ytellgreaterthannonell} from earlier, and the first equality follows from Lemma~\ref{lemma:leafdecomp}; the second inequality follows by constraints~\eqref{prob:TEOMIO_left} and \eqref{prob:TEOMIO_right}; and the last two equalities follow by simple algebra. This establishes that $(\xb,\yb)$ satisfies constraint~\eqref{prob:TEOLinearized_lower}. This establishes that $(\xb,\yb)$ is feasible for \problemeqref{prob:TEOLinearized}, which concludes the proof. \Halmos

\subsection{Proof of Proposition~\ref{prop:depth_ordering}}
\label{appendix:proof_depth_ordering}

Observe that the sets of tree-split pairs are nested in the following way:
\begin{equation*}
\bar{\Omega}_1 \subseteq \bar{\Omega}_2 \subseteq \dots \subseteq \bar{\Omega}_{d_{\max}}.
\end{equation*}
As a consequence, the feasible region of \problemeqref{prob:TEOMIO_barOmega} at depth $d$ is a superset of the feasible region of \problemeqref{prob:TEOMIO_barOmega} at depth $d+1$ and so we have
\begin{equation*}
Z^*_{MIO,1} \geq Z^*_{MIO,2} \geq \dots \geq Z^*_{MIO, d_{\max}}.
\end{equation*}
The equality $Z^*_{MIO, d_{\max}} = Z^*_{MIO}$ follows by the definition of $d_{\max}$ as the maximum depth of any tree in the ensemble. Combining this equality with the above sequence of inequalities establishes the result. \Halmos

\subsection{Proof of Theorem~\ref{theorem:depth_approximation}}
\label{appendix:proof_depth_approximation}

Let $(\xb, \yb)$ be an optimal solution of the depth $d$ problem (\ie, problem~\eqref{prob:TEOMIO_barOmega} with $\bar{\Omega}_d$). Let $Z^*_{MIO,d}$ be the objective value of $(\xb, \yb)$ within problem~\eqref{prob:TEOMIO_barOmega} with $\bar{\Omega}_d$.
 
Let $(\xb, \tilde{\yb})$ be the solution for \problemeqref{prob:TEOMIO} (the full-depth problem), obtained by finding the unique value of $\tilde{\yb}$ such that $(\xb, \tilde{\yb})$ is feasible for \problemeqref{prob:TEOMIO}. The existence and uniqueness of such a $\tilde{\yb}$ is guaranteed by Proposition~\ref{prop:SubPrimalBinaryOptimal}. Let $Z_d$ be the objective value of $(\bar{\xb}, \tilde{\yb})$ within \problemeqref{prob:TEOMIO}.

For a given tree $t \in \{1,\dots, T\}$, let us consider the difference of the prediction of tree $t$ for $(\xb, \yb)$ and $(\xb, \tilde{\yb})$: 
\begin{equation}
\sum_{\ell \in \leaves(t)} p_{t,\ell} y_{t,\ell} -  \sum_{\ell \in \leaves(t)} p_{t,\ell} \tilde{y}_{t,\ell}. \label{eq:prediction_difference}
\end{equation}

In order to understand this quantity, we need to understand which part of the tree $\bar{\xb}$ will get mapped to. We will do this through the following procedure:
\begin{enumerate}
\item Initialize $\nu$ to the root node of the tree. 
\item If $\nu \in \leaves(t)$ or $\nu \in \splits(t,d)$, stop; otherwise:
	\begin{itemize}
	\item If $\sum_{j \in \cond(\nu)} x_{\var(\nu),j} = 1$, set $\nu$ to its left child node;
	\item If $\sum_{j \in \cond(\nu)} x_{\var(\nu),j} = 0$, set $\nu$ to its right child node;
	\end{itemize}
	and repeat step 2.
\end{enumerate}
Upon termination, $\nu$ is some node in the tree -- either a leaf or a depth $d$ split node. Regardless of the type of node, we know that for any $s \in \leftsplits(\nu)$, the depth of $s$ is in $\{1,\dots, d\}$, and so $(\xb, \yb)$ satisfies the right split constraint~\eqref{prob:TEOMIO_barOmega_right} of the depth $d$ problem for $s$. By the definition of the procedure, it also must be that $\sum_{j \in \cond(s)} x_{\var(s),j} = 1$, which implies that
\begin{align*}
& \sum_{\ell \in \rightleaves(s)} y_{t,\ell} \leq 1 - \sum_{j \in \cond(s)} x_{\var{s},j} = 1 - 1 = 0, \\
& \Rightarrow \sum_{\ell \in \rightleaves(s)} y_{t,\ell} = 0, \\
& \Rightarrow y_{t,\ell} = 0, \quad \forall \ell \in \rightleaves(s). 
\end{align*}
For $s \in \rightsplits(\nu)$, similar reasoning using the left split constraint~\eqref{prob:TEOMIO_barOmega_left} allows us to assert that for any $s \in \rightsplits(\nu)$,
\begin{align*}
& \sum_{\ell \in \leftleaves(s)} y_{t,\ell} \leq \sum_{j \in \cond(s)} x_{\var{s},j} = 0, \\
& \Rightarrow \sum_{\ell \in \leftleaves(s)} y_{t,\ell} = 0, \\
& \Rightarrow y_{t,\ell} = 0, \quad \forall \ell \in \leftleaves(s). 
\end{align*}
We thus know the following about $\yb$:
\begin{equation*}
y_{t,\ell} = 0, \quad \forall \ \ell \in \bigcup_{s \in \leftsplits(\nu)} \rightleaves(s) \cup \bigcup_{s \in \rightsplits(\nu)} \leftleaves(s).
\end{equation*}
We can also assert the same about $\tilde{\yb}$, since $(\xb, \tilde{\yb})$ satisfies the left and right split constraints~\eqref{prob:TEOMIO_left} and \eqref{prob:TEOMIO_right} at all depths:
\begin{equation*}
\tilde{y}_{t,\ell} = 0, \quad \forall \ \ell \in \bigcup_{s \in \leftsplits(\nu)} \rightleaves(s) \cup \bigcup_{s \in \rightsplits(\nu)} \leftleaves(s).
\end{equation*}
We now consider three possible cases for the type of node $\nu$ is:
\begin{enumerate}
\item \textbf{Case 1}: $\nu \in \leaves(t)$. In this case, by Lemma~\ref{lemma:leafdecomp}, we can assert that 
\begin{equation*}
y_{t,\ell} = 0, \quad \forall \ \ell \neq \nu, 
\end{equation*}
\begin{equation*}
\tilde{y}_{t,\ell} = 0, \quad \forall \ \ell \neq \nu.
\end{equation*}
Since both $\yb$ and $\tilde{\yb}$ are nonnegative and sum to one, it follows that $y_{t,\nu} = \tilde{y}_{t,\nu} = 1$. We therefore have that the prediction difference~\eqref{eq:prediction_difference} is simply
\begin{equation}
\sum_{\ell \in \leaves(t)} p_{t,\ell} y_{t,\ell} -  \sum_{\ell \in \leaves(t)} p_{t,\ell} \tilde{y}_{t,\ell} = p_{t, \nu} - p_{t,\nu} = 0.
\end{equation}
\item \textbf{Case 2}: $\nu \in \splits(t,d)$ and $\sum_{j \in \cond(\nu)} x_{\var(\nu),j} = 1$. In this case, by Lemma~\ref{lemma:leafdecomp_splits}, we have that 
\begin{align*}
y_{t,\ell} = 0, \quad \forall \ \ell \notin \leftleaves(\nu) \cup \rightleaves(\nu), \\
\tilde{y}_{t,\ell} = 0, \quad \forall \ \ell \notin \leftleaves(\nu) \cup \rightleaves(\nu).
\end{align*}
In addition, since $\sum_{j \in \cond(\nu)} x_{\var(\nu),j} = 1$, then by the right split constraints, we additionally have $y_{t,\ell} = 0$ and $\tilde{y}_{t,\ell} = 0$ for all $\ell \in \rightleaves(\nu)$, which implies that 
\begin{align*}
y_{t,\ell} = 0, \quad \forall \ \ell \notin \leftleaves(\nu), \\
\tilde{y}_{t,\ell} = 0, \quad \forall \ \ell \notin \leftleaves(\nu).
\end{align*}
We can use the above properties of $\yb_t$ and $\tilde{\yb}_t$ to bound the prediction difference~\eqref{eq:prediction_difference} as follows:
\begin{align*}
\sum_{\ell \in \leaves(t)} p_{t,\ell} y_{t,\ell} - \sum_{\ell \in \leaves(t)} p_{t,\ell} \tilde{y}_{t,\ell} & \leq \max\left\{ \sum_{\ell } p_{t,\ell} y'_{t,\ell} \ \vline \ \sum_{\ell} y'_{t,\ell} = 1;\  y'_{t,\ell} \geq 0, \ \forall \ell;\  y'_{t,\ell} = 0, \ \forall \ell \notin \leftleaves(\nu) \right\}  \\
& \phantom{\leq} - \min\left\{ \sum_{\ell } p_{t,\ell} y'_{t,\ell} \ \vline \ \sum_{\ell} y'_{t,\ell} = 1;\  y'_{t,\ell} \geq 0, \ \forall \ell;\  y'_{t,\ell} = 0, \ \forall \ell \notin \leftleaves(\nu) \right\} \\
& = \max_{\ell \in \leftleaves(\nu)} p_{t,\ell} - \min_{\ell \in \leftleaves(\nu)} p_{t,\ell}.
\end{align*}
\item \textbf{Case 3}: $\nu \in \splits(t,d)$ and $\sum_{j \in \cond(\nu)} x_{\var(\nu),j} = 0$. Similar reasoning as in case 2 can be used to establish that 
\begin{align*}
y_{t,\ell} = 0, \quad \forall \ \ell \notin \rightleaves(\nu), \\
\tilde{y}_{t,\ell} = 0, \quad \forall \ \ell \notin \rightleaves(\nu),
\end{align*}
and to bound the prediction difference~\eqref{eq:prediction_difference} as
\begin{align*}
\sum_{\ell \in \leaves(t)} p_{t,\ell} y_{t,\ell} - \sum_{\ell \in \leaves(t)} p_{t,\ell} \tilde{y}_{t,\ell} & \leq \max\left\{ \sum_{\ell } p_{t,\ell} y'_{t,\ell} \ \vline \ \sum_{\ell} y'_{t,\ell} = 1;\  y'_{t,\ell} \geq 0, \ \forall \ell;\  y'_{t,\ell} = 0, \ \forall \ell \notin \rightleaves(\nu) \right\}  \\
& \phantom{\leq} - \min\left\{ \sum_{\ell } p_{t,\ell} y'_{t,\ell} \ \vline \ \sum_{\ell} y'_{t,\ell} = 1;\  y'_{t,\ell} \geq 0, \ \forall \ell;\  y'_{t,\ell} = 0, \ \forall \ell \notin \rightleaves(\nu) \right\} \\
& = \max_{\ell \in \rightleaves(\nu)} p_{t,\ell} - \min_{\ell \in \rightleaves(\nu)} p_{t,\ell}.
\end{align*}
\end{enumerate}

Given cases 2 and 3, observe that if we know that $\nu \in \splits(t,d)$, then a valid upper bound on the prediction difference is simply
\begin{align*}
\sum_{\ell \in \leaves(t)} p_{t,\ell} y_{t,\ell} - \sum_{\ell \in \leaves(t)} p_{t,\ell} \tilde{y}_{t,\ell} & \leq \max\left\{ \max_{\ell \in \leftleaves(\nu)} p_{t,\ell} - \min_{\ell \in \leftleaves(\nu)} p_{t,\ell},  \max_{\ell \in \rightleaves(\nu)} p_{t,\ell} - \min_{\ell \in \rightleaves(\nu)} p_{t,\ell} \} \right\}  \\
& =  \delta_{t,\nu}.
\end{align*}

Now, if we do not know what type of node $\nu$ is -- whether it is a leaf, or which split in $\splits(t,d)$ it is -- then we can construct an upper bound, based on cases 1, 2 and 3 above, for the prediction difference as
\begin{align*}
\sum_{\ell \in \leaves(t)} p_{t,\ell} y_{t,\ell} - \sum_{\ell \in \leaves(t)} p_{t,\ell} \tilde{y}_{t,\ell} & \leq \max_{\nu \in \splits(t,d)} \delta_{t,\nu} \\
& = \Delta_t,
\end{align*}
where the maximum is defined to be zero if $\splits(t,d)$ is empty. (Note that in the case $\nu$ is a leaf, the above bound is valid, since all $\delta_{t,\nu}$ values are nonnegative by definition.) 

Let us now unfix the tree $t$. Applying the above bound to bound the prediction difference of all trees $t$, and using the fact that $\lambda_t \geq 0$ for all $t$, it follows that the difference between the objective $Z^*_{MIO,d}$ of $(\xb,\yb)$ and the objective $Z_d$ of $(\xb, \tilde{\yb})$ can be written as 
\begin{align*}
Z^*_{MIO,d} - Z_d & = \sum_{t=1}^T \sum_{\ell \in \leaves(t)} \lambda_t \cdot p_{t,\ell} \cdot y_{t,\ell} - \sum_{t=1}^T \sum_{\ell \in \leaves(t)} \lambda_t \cdot p_{t,\ell} \cdot \tilde{y}_{t,\ell} \\
& =  \sum_{t=1}^T  \lambda_t \cdot \left( \sum_{\ell \in \leaves(t)} p_{t,\ell} \cdot y_{t,\ell} -  \sum_{\ell \in \leaves(t)} p_{t,\ell} \cdot \tilde{y}_{t,\ell} \right) \\
& \leq \sum_{t=1}^T  \lambda_t \cdot \Delta_t.
\end{align*}
From here, it immediately follows that 
$$ Z^*_{MIO,d} - \sum_{t=1}^T  \lambda_t \cdot \Delta_t \leq Z_d.$$
Combining this with the fact that $(\xb, \tilde{\yb})$ is a feasible solution for \problemeqref{prob:TEOMIO} and Proposition~\ref{prop:depth_ordering} leads to the inequality,
$$ Z^*_{MIO,d} - \sum_{t=1}^T  \lambda_t \cdot \Delta_t \leq Z_d \leq Z^*_{MIO} \leq Z^*_{MIO,d},$$
as required. \Halmos

\subsection{Proof of Proposition~\ref{prop:SubPrimalBinaryOptimal}}
\label{appendix:proof_SubPrimalBinaryOptimal}

\noindent \textbf{Feasibility}. Let us first show that the proposed solution $\yb_t$ is feasible. By construction, we have that $y_{t,\ell} \geq 0$ for all $\ell \in \leaves(t)$ and that $\sum_{\ell \in \leaves(t)} y_{t,\ell} = 1$, so constraints~\eqref{prob:SubPrimal_ysumtoone} and \eqref{prob:SubPrimal_ynonnegative} are satisfied. This leaves the left and right split constraints~\eqref{prob:SubPrimal_left} and \eqref{prob:SubPrimal_right}. 

For constraint~\eqref{prob:SubPrimal_left}, let $s \in \splits(t)$. If $\ell^* \in \leftleaves(s)$, then it must be that $s \in \leftsplits(\ell^*)$ and so by the definition of \getLeaf, it must be that $\sum_{j \in \cond(s)} x_{\var(s),j} =1$. Therefore, we have:
\begin{align*}
\sum_{\ell \in \leftleaves(s)} y_{t,\ell} & = \sum_{\substack{ \ell \in \leftleaves(s): \\ \ell \neq \ell^*} } y_{t,\ell} + y_{t,\ell^*} \\
& = 0 + 1 \\
& \leq \sum_{j \in \cond(s)} x_{\var(s),j} \\
& = 1. 
\end{align*}
Otherwise, if $\ell^* \notin \leftleaves(s)$, then 
\begin{align*}
\sum_{\ell \in \leftleaves(s)} y_{t,\ell} = 0,
\end{align*}
which is automatically less than or equal to $\sum_{j \in \cond(s)} x_{\var(s),j}$ (the latter can only be 0 or 1). 

For constraint~\eqref{prob:SubPrimal_right}, let $s \in \splits(t)$. If $\ell^* \in \rightleaves(s)$, then $s \in \rightsplits(\ell^*)$ and by the definition of \getLeaf, it must be that $1 - \sum_{j \in \cond(s)} x_{\var(s),j} = 1$. Therefore, applying similar reasoning as above, we get
\begin{align*}
\sum_{\ell \in \rightleaves(s)} y_{t,\ell} & = \sum_{\substack{ \ell \in \rightleaves(s): \\ \ell \neq \ell^*} } y_{t,\ell} + y_{t,\ell^*} \\
& = 0 + 1 \\
& \leq 1 - \sum_{j \in \cond(s)} x_{\var(s),j} \\
& = 1.
\end{align*}

Otherwise, if $\ell^* \notin \rightleaves(s)$, then again, we have $\sum_{\ell \in \leftleaves(s)} y_{t,\ell} = 0$, which is automatically less than or equal to $1 - \sum_{j \in  \cond(s)} x_{\var(s),j}$ (again, it can only be 0 or 1). This establishes that $\yb_t$ is a feasible solution to the subproblem~\eqref{prob:SubPrimal}. \\

\noindent \textbf{Unique feasible solution}. To establish that the proposed solution is the only feasible solution, we proceed as follows. We will show that if a solution $\yb_t$ is a feasible solution of \problemeqref{prob:SubPrimal}, then it must be equal to the solution of the statement of the proposition. 

Let $\ell \in \leaves(t)$ such that $\ell \neq \ell^*$. Then by Lemma~\ref{lemma:leafdecomp}, we have that 
$$\ell \in \bigcup_{s \in \leftsplits(\ell^*)} \rightleaves(s) \cup \bigcup_{s \in \rightsplits(\ell^*)} \leftleaves(s). $$
Moreover, the collection of sets in the union above is disjoint. Therefore, either $\ell \in \rightleaves(\bar{s})$ for some $\bar{s} \in \leftsplits( \ell^* )$ or $\ell \in \leftleaves(\bar{s})$ for some $\bar{s} \in \rightsplits( \ell^* )$. 

In the former case -- that is, $\ell \in \rightleaves(\bar{s})$ for some $\bar{s} \in \leftsplits(\ell^*)$ -- we have by constraint~\eqref{prob:SubPrimal_left} and constraint~\eqref{prob:SubPrimal_ynonnegative} that
\begin{align*}
1 - \sum_{j \in \cond(\bar{s})} x_{\var(\bar{s}),j} & \geq \sum_{\ell' \in \rightleaves(\bar{s})} y_{t,\ell'} \\
& \geq y_{t,\ell}.
\end{align*}
Therefore, $y_{t,\ell}$ is upper bounded by $1 - \sum_{j \in \cond(\bar{s})} x_{\var(\bar{s}),j}$ and lower bounded by 0 (by constraint~\eqref{prob:SubPrimal_ynonnegative}). Since $\bar{s} \in \leftsplits(\ell^*)$ and from the definition of $\getLeaf$, it must be that $\sum_{j \in \cond(\bar{s})} x_{\var(\bar{s}),j} = 1$, or equivalently, $1 - \sum_{j \in \cond(\bar{s})} x_{\var(\bar{s}),j} = 0$. Therefore, $y_{t,\ell}$ must be equal to zero. 

Similarly, if $\ell \in \leftleaves(\bar{s})$ for some $\bar{s} \in \rightsplits(\ell^*)$, then by constraints~\eqref{prob:SubPrimal_right} and \eqref{prob:SubPrimal_ynonnegative} we have that
\begin{align*}
\sum_{j \in \cond(\bar{s})} x_{\var(\bar{s}),j} & \geq \sum_{\ell' \in \leftleaves(\bar{s})} y_{t,\ell'} \\
& \geq y_{t,\ell}.
\end{align*}
Therefore, $y_{t,\ell}$ is upper bounded by $\sum_{j \in \cond(\bar{s})} x_{\var(\bar{s}),j}$ and lower bounded by 0 (by \eqref{prob:SubPrimal_ynonnegative}). Since $\bar{s} \in \leftsplits(\ell^*)$ and from the definition of $\getLeaf$, it must be that $\sum_{j \in \cond(\bar{s})} x_{\var(\bar{s}),j} = 0$. Therefore, $y_{t,\ell}$ must be equal to zero.

From the above reasoning, we have shown that $y_{t,\ell} = 0$ for every $\ell \neq \ell^*$. By constraint~\eqref{prob:SubPrimal_ysumtoone}, it must be that $y_{t,\ell^*} = 1 - \sum_{\ell \neq \ell^*} y_{t,\ell} = 1 - 0 = 1$. The resulting solution is therefore \emph{exactly the same} as the one proposed in the proposition; it follows that the proposed solution $\yb_t$ is the only feasible solution to \problemeqref{prob:SubPrimal}. \\

\noindent \textbf{Optimality}. Since $\yb_t$ is the only feasible solution of \problemeqref{prob:SubPrimal}, it must also be its optimal solution. This completes the proof. \Halmos

\subsection{Proof of Proposition~\ref{prop:SubDualBinaryOptimal}}
\label{appendix:proof_SubDualBinaryOptimal}

First, let us check the dual objective of the proposed solution $(\alphab_t ,\betab_t, \gamma_t)$. We have
\begin{align*}
& \sum_{s \in \splits(t)} \alpha_{t,s} \left[ \sum_{j \in \cond(s)} x_{\var(s), j} \right] + 
\sum_{s \in \splits(t)} \beta_{t,s} \left[ 1 - \sum_{j \in \cond(s)} x_{\var(s), j} \right] + 
\gamma_t \\
& = \sum_{s \in \rightsplits(\ell^*)} \alpha_{t,s} \left[ \sum_{j \in \cond(s)} x_{\var(s), j} \right] + 
\sum_{s \in \leftsplits(\ell^*)} \beta_{t,s} \left[ 1 - \sum_{j \in \cond(s)} x_{\var(s), j} \right] + 
p_{t,\ell^*} \\
& = 0 + 0 + p_{t,\ell^*} \\
& = p_{t,\ell^*},
\end{align*}
where the first step follows by the definition of $(\alphab_t, \betab_t, \gamma_t)$; the second step follows by the fact that $\sum_{j \in \cond(s)} x_{\var(s),j} = 0$ for $s \in \rightsplits(\ell^*)$ and $1 - \sum_{j \in \cond(s)} x_{\var(s),j} = 0$ for $s \in \leftsplits(\ell^*)$; and the last two steps by algebra. 
The final value is exactly equal to the optimal primal objective value. If $(\alphab_t, \betab_t, \gamma_t)$ is feasible for the dual problem, then the proposition will be proven. 

To verify feasibility, observe that by their definition, we have $\alpha_{t,s} \geq 0$ and $\betab_{t,s} \geq 0$ for all $s \in \splits(t)$. Thus, we only need to check constraint~\eqref{prob:SubDual_constraint} for each $\ell \in \leaves(t)$. We consider two cases: \\

\noindent \textbf{Case 1}: $\ell = \ell^*$. In this case, proceeding from the left hand side of constraint~\eqref{prob:SubDual_constraint} for $\ell^*$, we have
\begin{align*}
& \sum_{s \in \leftsplits(\ell^*)} \alpha_{t,s} + \sum_{s \in \rightsplits(\ell^*)} \beta_{t,s} + \gamma_t \\
& = 0 + 0 + p_{t,\ell^*} \\
& \geq p_{t,\ell^*},
\end{align*}
where the first equality follows because $\alpha_{t,s} = 0$ for all $s \notin \rightsplits(\ell^*)$ and $\beta_{t,s} = 0$ for all $s \notin \leftsplits(\ell^*)$ (note that $\leftsplits(\ell^*) \cap \rightsplits(\ell^*) = \emptyset$; a leaf cannot be both to the left of and to the right of the same split), and also because $\gamma_t = p_{t,\ell^*}$ by definition. \\

\noindent \textbf{Case 2}: $\ell \neq \ell^*$. In this case, by Lemma~\eqref{lemma:leafdecomp}, we know that $\ell$ satisfies 
\begin{equation*}
\ell \in \bigcup_{s \in \leftsplits(\ell^*)} \rightleaves(s) \cup \bigcup_{s \in \rightsplits(\ell^*)} \leftleaves(s).
\end{equation*}
Lemma~\ref{lemma:leafdecomp} states that each set in the above union is disjoint. Therefore, we have that $\ell \in \rightleaves(\bar{s})$ for some $\bar{s} \in \leftsplits(\ell^*)$ or $\ell \in \leftleaves(\bar{s})$ for some $\bar{s} \in \rightsplits(\ell^*)$. We now show that the inequality holds in either of these two scenarios.

If $\ell \in \rightleaves(\bar{s})$ for some $\bar{s} \in \leftsplits(\ell^*)$, then we have 
\begin{align*}
\sum_{s \in \leftsplits(\ell)} \alpha_{t,s} + \sum_{s \in \rightsplits(\ell)} \beta_{t,s} + \gamma_t & \geq \beta_{t,\bar{s}} + \gamma_t \\
& = \max\{0 , \max_{\ell' \in \rightleaves(\bar{s})} ( p_{t,\ell'} - p_{t,\ell^*} ) \} + p_{t,\ell^*} \\
& \geq p_{t,\ell} - p_{t,\ell^*} + p_{t,\ell^*} \\
& = p_{t,\ell},
\end{align*}
where the first step follows because $\bar{s}$ must belong to $\rightsplits(\ell)$ and because by definition, $\alphab_t$ and $\betab_t$ are nonnegative; the second step follows by the definition of $\beta_{t,s}$ for $s \in \leftsplits(\ell^*)$; the third step by the definition of the maximum; and the last step by algebra. 

If $\ell \in \leftleaves(\bar{s})$ for some $\bar{s} \in \rightsplits(\ell^*)$, then we have
\begin{align*}
\sum_{s \in \leftsplits(\ell)} \alpha_{t,s} + \sum_{s \in \rightsplits(\ell)} \beta_{t,s} + \gamma_t & \geq \alpha_{t, \bar{s}} + \gamma_t \\
& = \max\{0, \max_{\ell' \in \leftleaves(\bar{s})} ( p_{t,\ell'} - p_{t,\ell^*} ) \} + p_{t,\ell^*} \\
& \geq p_{t,\ell} - p_{t,\ell^*} + p_{t,\ell^*} \\
& = p_{t,\ell},
\end{align*}
which follows by logic similar to the first case ($\ell \in \rightleaves(\bar{s})$ for some $\bar{s} \in \leftsplits(\ell^*)$). 

Thus, we have established that $(\alphab_t, \betab_t,\gamma_t)$ is a feasible solution for the dual \problemeqref{prob:SubDual} and achieves the same objective as the primal optimal solution. By weak LO duality, the solution $(\alphab_t, \betab_t, \gamma_t)$ must therefore be optimal. \Halmos

\subsection{Proof of Proposition~\ref{prop:TreeTraversalFeasibility}}
\label{appendix:proof_TreeTraversalFeasibility}

For the $\Rightarrow$ direction, the implication is immediate, because if $(\xb,\yb)$ satisfies constraints~\eqref{prob:TEOMIO_left} and \eqref{prob:TEOMIO_right} for all $s \in \splits(t)$, it will also satisfy them for arbitrary subsets of $\splits(t)$. 

Thus, we only need to establish the $\Leftarrow$ direction. Let $s \in \leftsplits(\ell^*)$. For any tree $t$, any $s \in \leftsplits(\ell^*)$ and any $\ell \in \rightleaves(s)$,
\begin{align*}
y_{t,\ell} & \leq \sum_{\ell \in \rightleaves(s)} y_{t,\ell'}  \\
& \leq 1 - \sum_{j \in \cond(s)} x_{\var(s), j} \\
& = 1 - 1 \\
& = 0,
\end{align*}
where the first step follows since $y_{t,\ell'}$ is nonnegative for any $\ell'$, the second step by the hypothesis of the implication, and the remaining two steps by algebra. Since $y_{t,\ell}$ is nonnegative, it must be that $y_{t,\ell} = 0$. Similarly, for any $s \in \rightsplits(\ell^*)$,  we can show that for each $\ell \in \leftleaves(s)$, $y_{t,\ell} = 0$. 

It therefore follows that for any $\ell \in \bigcup_{s \in \rightsplits(\ell^*)} \leftleaves(s) \cup \bigcup_{s \in \leftsplits(\ell^*)} \rightleaves(s)$, $y_{t,\ell}$. Invoking Lemma~\ref{lemma:leafdecomp}, we have that $y_{t,\ell} = 0$ for all $\ell \in \leaves(t) \setminus \{ \ell^* \}$. Since $\yb$ is assumed to satisfy constraint~\eqref{prob:TEOMIO_ysumtoone}, it must be that $y_{t,\ell^*} = 1$. 

From here, we can see that $\yb_t$, the collection of $y_{t,\ell}$ values corresponding to tree $t$, is defined exactly as in the statement of Proposition~\ref{prop:SubPrimalBinaryOptimal}. Thus, invoking Proposition~\ref{prop:SubPrimalBinaryOptimal}, we can assert that $\yb$ satisfies 
\begin{align*}
& \sum_{\ell \in \leftleaves(s)} y_{t,\ell} \leq \sum_{j \in \cond(s)} x_{\var(s), j}, \\
& \sum_{\ell \in \rightleaves(s)} y_{t,\ell} \leq 1- \sum_{j \in \cond(s)} x_{\var(s), j}, 
\end{align*}
for all $s \in \splits(t)$. This concludes the proof. \Halmos

\clearpage
\pagebreak

\section{Local search procedure}
\label{appendix:localsearch}

We provide the pseudocode of our local search procedure for approximately solving \problemeqref{prob:TreeEnsembleOpt_abstract} below as Algorithm~\ref{algorithm:local_search}. Before we define the procedure, we define the set $\bar{\Xcal}_i$ as 
$$ \bar{\Xcal}_i = \Xcal_i$$
for categorical variables, and
$$ \bar{\Xcal}_i = \{  a_{i,j} \, | \, j \in \{1,\dots, K_i\}\} \cup \{ a_{i,K_i} + 1\}$$
for numeric variables, where $a_{i,j}$ is the $j$th smallest split point of variable $i$ in the tree ensemble model. The $\bar{\Xcal}_i$ are simply the domains of each independent variable defined in a way that will be helpful in defining our local search. For numeric variables, $\bar{\Xcal}_i$ consists of the $K_i$ split points of variable $i$ and one extra point, $a_{i,K_i}+1$. The extra point $a_{i,K_i}+1$ is arbitrary. We can use any value here, as long as it is strictly larger than the largest split point of variable $i$, as this will allow us to choose to be on the right-hand side of a split with query $X_i \leq a_{i,K_i}$.

\begin{algorithm}
\begin{algorithmic}
\REQUIRE Tree ensemble model $f_1(\cdot), \dots , f_T(\cdot)$, $\lambda_1, \dots, \lambda_T$; finite domains $\bar{\Xcal}_1, \dots, \bar{\Xcal}_n$. 
\STATE Select $\Xb = (X_1, \dots, X_n)$ uniformly at random from $\prod_{i'=1}^n \bar{\Xcal}_{i'}$. \\
\STATE Initialize $Z \leftarrow \sum_{t=1}^T \lambda_t f_{t} (\Xb)$. \\
\STATE Initialize $\untestedVars = \{ 1,\dots, n\}$. \\
\WHILE{ $| \untestedVars | > 0 $ }
	\STATE Select $i \in \untestedVars$. \\
	\STATE Set $\Mcal \leftarrow \{ \Xb' \in \prod_{i'=1}^n \bar{\Xcal}_{i'} \, | \, X'_{j} = X_{j} \ \text{for}\ j \neq i \}$. \\
	\STATE Set $\Xb^*  \leftarrow \arg \max_{\Xb' \in \Mcal }  \sum_{t=1}^T \lambda_t f_{t}  ( \Xb')$. \\
	\STATE Set $Z_c \leftarrow  \sum_{t=1}^T \lambda_t f_{t}  ( \Xb^*) $. \\
	\IF{$Z_c > Z$}
		\STATE Set $Z \leftarrow Z_c$. \\
		\STATE Set $\Xb \leftarrow \Xb^*$. \\
		\STATE Set $\untestedVars \leftarrow \{1, \dots, i-1, i+1, \dots, n\}$. \\
	\ELSE
		\STATE Set $\untestedVars \leftarrow \untestedVars \setminus \{ i \}$. \\
	\ENDIF
\ENDWHILE
\RETURN Locally optimal solution $\Xb$ with objective value $Z$. 
\end{algorithmic}
\caption{Local search procedure. \label{algorithm:local_search}}
\end{algorithm}

\clearpage
\pagebreak

\section{Case study 2: customized pricing} 
\label{sec:pricing}

In this section, we apply our approach to customized pricing. Section~\ref{subsec:pricing_background} provides the background on the data, while Section~\ref{subsec:pricing_models} describes our random forest model as well as two alternative models based on hierarchical Bayesian regression. Section~\ref{subsec:pricing_predictive_accuracy} compares the models in terms of out-of-sample predictions of profit. Finally, Section~\ref{subsec:pricing_optimization} formulates the profit optimization problem and compares the three models.

\subsection{Background}
\label{subsec:pricing_background}

We consider the data set from \cite{montgomery1997creating}, which was accessed via the \texttt{bayesm} package in R \citep{rossi2012bayesm}. This data set contains price and sales data for eleven different refrigerated orange juice brands for the Dominick's Finer Foods chain of grocery stores in the Chicago area.

In this data set, each observation corresponds to a given store in the chain at a given week. The data span 83 stores and a period of 121 weeks. Each observation consists of: the week $t$; the store $s$; the sales $q_{t,s,1}, \dots, q_{t,s,11}$ of the eleven orange juice brands; the prices $p_{t,s,1}, \dots, p_{t,s,11}$ of the eleven orange juice brands; dummy variables $d_{t,s,1}, \dots, d_{t,s,11}$, where $d_{t,s,i} = 1$ if orange juice brand $i$ had any in-store displays (such as in-store coupons) at store $s$ in week $t$; and dummy variables $f_{t,s,1}, \dots, f_{t,s,11}$, where $f_{t,s,i} = 1$ if brand $i$ was featured/advertised in store $s$ in week $t$. We use $\pb$, $\db$ and $\fb$ to denote vectors of prices, deal dummies and feature dummies, and subscripts to denote the observation; for example, $\pb_{t,s} = ( p_{t,s,1}, \dots, p_{t,s,11})$ is the vector of brand prices at store $s$ in week $t$. 

The data set also include 11 covariates for the 83 stores corresponding to demographic information of each store's neighborhood, such as the percentage of the population with a college degree and the percentage of households with more than five members; we denote these covariates as $z_{s,1}, \dots, z_{s,11}$. We denote the vector of covariates for store $s$ as $\zb_s$; for notational convenience later, we assume that $z_{s,0} = 1$ in $\zb_s$. For more details, we refer the reader to \cite{montgomery1997creating}.

\subsection{Models}
\label{subsec:pricing_models}

The standard approach for modeling this type of data in marketing is to posit a regression model within a hierarchical Bayes (HB) framework \citep{rossi2005bayesian}. We will consider two different HB specifications. In the first model, which we will refer to as HB-LogLog, we assume that the logarithm of sales of a given product is linear in the logarithm of prices of all eleven products, and linear in the deal and feature dummies. For a fixed focal brand $i$, the regression function is given below:
\begin{equation}
\log(q_{t,s,i}) = \beta_{s,0} + \betab_{s,i,\pb}^T \log(\pb_{t,s}) + \betab_{s,i,\db}^T \db_{t,s} + \betab_{s,i, \fb}^T  \fb_{t,s}  + \epsilon_{t,s,i}, \label{eq:HBLogLog_regression}
\end{equation}
where $\epsilon_{t,s,i}$ follows a univariate normal distribution with mean 0 and variance $\tau_s$:
\begin{equation}
\epsilon_{t,s,i} \sim \Normal(0, \tau_{s,i}),  \label{eq:epsilon_first_stage_prior}
\end{equation}
and the vector of regression coefficients $\betab_{s,i} = (\beta_{s,i,0}, \betab_{s,i,\pb}, \betab_{s,i,\db}, \betab_{s,i,\fb} )$ follows a multivariate normal distribution with mean $\Delta_{i}^T \zb_s$ and covariance matrix $V_{\betab,i}$:
\begin{equation}
\betab_{s,i} \sim \Normal( \Delta_{i}^T \zb_s, V_{\betab,i} ), \label{eq:beta_first_stage_prior}
\end{equation}
where $\Delta_i$ is a 12-by-34 matrix. (Each row corresponds to one of the store level covariates, and each column corresponds to one of the $\beta$ coefficients.) This model assumes that the regression coefficients of each store $\betab_{s,i}$ follow a normal distribution whose mean depends on the store-specific covariates $\zb_s$. Specifically, the mean of $\betab_{s,i}$ is a linear function of $\zb_s$; recall that $z_{s,0} = 1$, so that the first row of $\Delta$ specifies the intercept of the mean of $\betab_{s,i}$. We note that log-log models like equation~\eqref{eq:HBLogLog_regression} are commonly used in demand modeling. The above distributions in equations~\eqref{eq:epsilon_first_stage_prior} and \eqref{eq:beta_first_stage_prior} are sometimes referred to as first-stage priors.

For brand $i$'s sales, we model $\Delta_i$, $\tau_{1,i}, \dots, \tau_{83,i}$ and $V_{\betab,i}$ as random variables with the following prior distributions:
\begin{equation}
\vecop(\Delta_i) \, | \, V_{\betab,i} \sim \Normal(\vecop(\bar{\Delta}), V_{\betab,i} \otimes A^{-1} ), \label{eq:Delta_prior}
\end{equation}
\begin{equation}
V_{\betab,i} \sim \IW(\nu, V), \label{eq:V_prior}
\end{equation}
\begin{equation}
\tau_{s,i} \sim C / \chi^2_{\nu_{\epsilon}}, \quad \forall \ s \in \{1,\dots, 83\}, \label{eq:tau_prior}
\end{equation}
where $\vecop(\bar{\Delta})$ is the elements of the matrix $\bar{\Delta}$ (also 12-by-34) stacked into a column vector; $\otimes$ denotes the Kronecker product; $IW(\nu,V)$ is the inverse Wishart distribution with degrees of freedom $\nu$ and scale matrix $V$; and $\chi^2_{\nu_e}$ is a chi-squared distributed random variable with $\nu_e$ degrees of freedom. The matrices $\bar{\Delta}$, $V$ and $A$ and the scalars $\nu$, $\nu_{\epsilon}$ and $C$ are the prior hyperparameters. The distributions~\eqref{eq:Delta_prior}, \eqref{eq:V_prior} and \eqref{eq:tau_prior} are sometimes referred to as second-stage priors.

The second type of HB model that we will consider, which we will refer to as HB-SemiLog, we assume the same specification above except that instead of equation~\eqref{eq:HBLogLog_regression}, where we assume the logarithm of sales is linear in the logarithm of prices, we assume it is linear in the actual prices themselves:
\begin{equation}
\log(q_{t,s,i}) = \beta_{s,0} + \betab_{s,i,\pb}^T \pb_{t,s} + \betab_{s,i,\db}^T \db_{t,s} + \betab_{s,i, \fb}^T  \fb_{t,s}  + \epsilon_{t,s,i}, \label{eq:HBLogLinear_regression}
\end{equation}
This particular HB model is essentially the same model proposed in \cite{montgomery1997creating}. One minor difference is that the above model allows for cross-brand promotion effects (e.g., brand 1 being on a deal or being featured can affect sales of brand 2). In our comparison of predictive accuracies in the next section, we will test both the above regression model, which allows for cross-brand promotion effects, and a simpler regression model with only own-brand promotion effects (\ie, equation~\eqref{eq:HBLogLog_regression} and \eqref{eq:HBLogLinear_regression} excludes $d_{t,s,i'}$ and $f_{t,s,i'}$ for brands $i'$ different from $i$, as in \citealt{montgomery1997creating}).

In addition to these two HB models, we also consider a random forest model, which we denote by RF. The random forest model that we estimate will be different from the previous two HB models. Rather than predicting the sales $q_{t,s,i}$ or log sales $\log(q_{t,s,i})$, we will instead predict profit, that is, 
\begin{equation}(p_{t,s,i} - c_i) \cdot q_{t,s,i},
\end{equation}
where $c_i$ is the unit cost of brand $i$. For each brand $i$, we estimate the random forest model using $p_1,\dots,p_{11}, d_1, \dots, d_{11}, f_1, \dots, f_{11}, z_1, \dots z_{11}$ as the predictor variables, where $z_1, \dots, z_{11}$ are the store-level covariates of the store at which the prediction will be made. Letting $F_i$ be the random forest prediction of profit for brand $i$, the overall profit that is predicted from the price and promotion decision $\pb, \db, \fb$ at a store with covariate vector $\zb$ is simply the sum of these individual brand-level predictions:
$$ F_{total}(\pb, \db, \fb, \zb) = \sum_{i=1}^{11} F_i(\pb, \db, \fb, \zb).$$
The choice to predict profit as opposed to sales or log of sales is motivated by tractability. By predicting the profit from each brand, the total profit function $F_{total}$ can be directly used within the tree ensemble formulation. In contrast, if $F_i(\cdot)$ was a prediction of sales, then our objective function would be
$$ \sum_{i=1}^{11} (p_{i} - c_i) \cdot F_{i}(\pb, \db, \fb, \zb),$$
which would require modifications to our MIO formulation -- one would need to model the tree ensemble behavior of each $F_i$, and then model the product of the price $p_i$ and the predicted sales. Although our formulation can be suitably modified to do this, the resulting model is harder to solve than the basic problem~\eqref{prob:TEOMIO}. %

\subsection{Predictive accuracy results}

\label{subsec:pricing_predictive_accuracy}

Our first experiment does not consider optimization, but merely compares the out-of-sample predictive accuracy of the three models -- HB-LogLog, HB-SemiLog and RF. To do so, we proceed as follows. We take the whole data set of store-week observations, and we split it randomly into a training set and a test set. We use the training set to estimate the HB-LogLog, HB-SemiLog and RF models for each brand. For each model, we predict the profit from brand $i$ with the store-week observations in the test set, and compute the $R^2$ of the test set profit predictions. %
We repeat this procedure with ten random splits of the observations into train and test sets. In each split, 80\% of the observations are randomly chosen for the training set, and the remaining 20\% are used for the test set. 

For HB-LogLog, we estimate it using Markov chain Monte Carlo (MCMC) using the R package \texttt{bayesm}. We use default values for the prior hyperparameters provided by \texttt{bayesm}. We run MCMC for 2000 iterations to obtain samples from the posterior distribution of $\betab_s$ for each store $s$, and samples from the posterior distribution of $\tau_s$ for each store $s$. Due to mild autocorrelation in the draws of $\betab_s$ and $\tau_s$, we thin the draws by retaining every fifth draw. Of these thinned draws, we retain the last $J = 100$ samples. We index the draws/posterior samples by $j = 1,\dots,J$. For an arbitrary price and promotion decision $\pb, \db, \fb$, we compute the predicted sales of brand $i$ at store $s$ by computing an approximation to the posterior expectation of sales $q_{s,i}$ under HB-LogLog as
$$ \hat{q}_{s,i} = \frac{1}{J} \sum_{j=1}^{J} \exp \left( \beta^{(j)}_{s,i,0} + (\betab^{(j)}_{s,i,\pb})^T \log(\pb) + (\betab^{(j)}_{s,i, \db})^T \db + (\betab^{(j)}_{s,i,\fb})^T \fb + \tau_{s,i}^{(j)}  / 2  \right),$$
where $\log(\pb)$ is the component-wise logarithm of $\pb$, and quantities with the superscript $(j)$ correspond to the $j$th posterior sample from the appropriate posterior distribution (that of either $\betab_{s,i}$ or $\tau_{s,i}$). With this prediction of sales, we predict the profit of brand $i$ as $(p_i - c_i) \cdot \hat{q}_{s,i}$. 

For HB-SemiLog, we proceed in the same way as for HB-LogLog, except that the predicted sales are computed as 
$$ \hat{q}_{s,i} = \frac{1}{J} \sum_{j=1}^{J} \exp \left( \beta^{(j)}_{s,i,0} + (\betab^{(j)}_{s,i,\pb})^T \pb + (\betab^{(j)}_{s,i, \db})^T \db + (\betab^{(j)}_{s,i,\fb})^T \fb + \tau_{s,i}^{(j)} / 2  \right).$$

For RF, we use the R package \texttt{ranger} to estimate one model for each brand's profit. We use default parameters, with the exception of the number of trees which we vary in $\{20, 50, 100, 500\}$. For the RF model, we emphasize again that, unlike HB-LogLog and HB-SemiLog which first predict log sales and then translates this to revenue, the RF model directly predicts profits of each brand. 

In addition to the above three models, we also consider slightly modified versions of the above three models where we do not allow for cross promotion effects (\ie, for each brand $i$'s sales or profit prediction, we leave out $d_{i'}$ and $f_{i'}$ for $i' \neq i$ as independent variables). In our presentation below, we distinguish these models from the ones above by using the suffix ``-Own'' -- thus, HB-LogLog-Own, HB-SemiLog-Own and RF-Own are the log-log, semi-log and random forest models with only own-brand promotion effects.

The prediction tasks we will consider involve predicting profit, which requires us to specify the unit cost $c_i$ of each brand. Since the data set does not include this information, we will test two different sets of values for $\cb = (c_1, \dots, c_{11})$. In the first set, we set $c_i = 0$ for each brand $i$; we are thus effectively predicting the revenue from each brand. In the second set, we set $c_i = 0.9 \times p_{i,\min}$, where $p_{i,\min}$ is the lowest price at which brand $i$ was offered in the whole data set. 

Table~\ref{table:orange_juice_task1_Rsq} displays the test set/out-of-sample $R^2$ values for the profit predictions of each brand in the first prediction task ($c_i = 0$ for all brands $i$), averaged over the ten random splits of the data. From this table, we can see that both the HB-LogLog and HB-SemiLog models are very inaccurate for some brands, achieving $R^2$ values that are close to zero. In one case, namely brand 10, the out-of-sample $R^2$ is even negative, indicating that the model is worse than a naive model that just predicts the average training set profit. For the log-log model, when we remove cross-brand promotion effects and move from HB-LogLog to HB-LogLog-Own, the out-of-sample $R^2$ exhibits an absolute improvement ranging from 0.04 to 0.68, with an average over all brands of 0.15. For the semi-log model, when we remove cross-brand promotion effects and move from HB-SemiLog to HB-SemiLog-Own, the $R^2$ exhibits an absolute improvement ranging from 0.03 to 0.85, with an average of 0.17. With regard to the semi-log model, this finding is consistent with that of \cite{montgomery1997creating}, where an own-brand promotion effect specification yielded a lower Schwartz information criterion value and higher out-of-sample accuracy than the cross-brand specification. 

\begin{table}
\small
\centering
\begin{tabular}{lccccccccccc} \toprule   %
& \multicolumn{11}{c}{Brand $R^2$} \\
  Model & 1 & 2 & 3 & 4 & 5 & 6 & 7 & 8 & 9 & 10 & 11 \\  \midrule 
  HB-SemiLog & 0.34 & 0.74 & 0.41 & 0.07 & 0.38 & 0.53 & 0.47 & 0.34 & 0.25 & --0.43 & 0.55 \\[0.1em] 
  HB-SemiLog-Own & 0.48 & 0.80 & 0.52 & 0.35 & 0.43 & 0.62 & 0.53 & 0.42 & 0.28 & 0.42 & 0.65 \\[0.5em]
  HB-LogLog & 0.39 & 0.75 & 0.47 & 0.08 & 0.38 & 0.55 & 0.61 & 0.38 & 0.45 & --0.17 & 0.58 \\[0.1em] 
  HB-LogLog-Own & 0.50 & 0.80 & 0.59 & 0.36 & 0.45 & 0.63 & 0.66 & 0.45 & 0.49 & 0.50 & 0.66 \\[0.5em]
  RF, $T = 20$ & 0.75 & 0.84 & 0.69 & 0.83 & 0.83 & 0.69 & 0.76 & 0.64 & 0.68 & 0.79 & 0.73 \\ 
  RF, $T = 50$ & 0.77 & 0.85 & 0.70 & 0.84 & 0.84 & 0.70 & 0.77 & 0.66 & 0.69 & 0.79 & 0.74 \\ 
  RF, $T = 100$ & 0.77 & 0.85 & 0.70 & 0.84 & 0.84 & 0.71 & 0.77 & 0.67 & 0.70 & 0.80 & 0.74 \\ 
  RF, $T = 500$ & \bfseries 0.78 & \bfseries 0.85 & \bfseries 0.70 & \bfseries 0.84 & \bfseries 0.84 & \bfseries 0.71 & 0.77 & \bfseries 0.67 & \bfseries 0.70 & \bfseries 0.80 & \bfseries 0.74 \\[0.25em]
  RF-Own, $T = 20$ & 0.73 & 0.84 & 0.68 & 0.81 & 0.79 & 0.70 & 0.75 & 0.64 & 0.68 & 0.74 & 0.72 \\ 
  RF-Own, $T = 50$ & 0.74 & 0.84 & 0.69 & 0.82 & 0.80 & 0.70 & 0.77 & 0.65 & 0.68 & 0.76 & 0.73 \\ 
  RF-Own, $T =100$ & 0.74 & 0.84 & 0.69 & 0.82 & 0.80 & 0.71 & 0.77 & 0.66 & 0.68 & 0.76 & 0.73 \\ 
  RF-Own, $T = 500$ & 0.75 & 0.84 & 0.70 & 0.82 & 0.80 & 0.71 & \bfseries 0.77 & 0.66 & 0.69 & 0.77 & 0.74 \\ \bottomrule
   \end{tabular}
  \caption{Comparison of out-of-sample profit prediction $R^2$ for the different models and brands, averaged over ten random splits of the data, for the first prediction task ($c_i = 0$ for all brands $i$). The best $R^2$ value for each brand is indicated in bold. \label{table:orange_juice_task1_Rsq} }
\end{table}

Comparing the HB models to the RF model with 500 trees, we can see that RF provides a significant improvement in predictive accuracy. For example, for brand 8, the highest $R^2$ attained by HB-LogLog, HB-LogLog-Own, HB-SemiLog and HB-SemiLog-Own, is 0.45. In contrast, the $R^2$ attained by RF with 500 trees is 0.67, which is an absolute improvement of 0.25. Over all of the brands, the improvement of RF over the best HB model for each brand ranges from 0.05 (brand 2) to as much as 0.48 (brand 4). 

Within the family of RF models, Table~\ref{table:orange_juice_task1_Rsq} gives us a sense of how the number of trees affects the predictive accuracy. In particular, while the out-of-sample accuracy decreases as the number of trees is decreased, we can see that the loss in accuracy is very modest. For example, for brand 1, the $R^2$ is 0.7778 with 500 trees, which is reduced to 0.7663 when we use 50 trees; note that this is still higher than any of the HB models. This suggests that RF can still achieve an improvement over the HB models even with a smaller number of trees. 

We can also determine the impact of cross promotion effects within RF. Note that unlike HB-LogLog and HB-SemiLog, where the out-of-sample $R^2$ improves once cross-promotion effects ($d_{i'}$ and $f_{i'}$ for $i'$ different to the focal brand $i$) are removed, the opposite happens with the random forest models: RF-Own has slightly lower $R^2$ values than RF. %

Table~\ref{table:orange_juice_task2_Rsq} presents the out-of-sample $R^2$ values for the profit predictions of each brand in the second prediction task ($c_i = 0.9 \times p_{i,\min}$ for all brands $i$), averaged over the ten random splits of the data. The same insights about the relative performance of the three different families of models derived from Table~\ref{table:orange_juice_task1_Rsq} hold for this case. Overall, these results provide evidence that a random forest model for profit predictions can outperform state-of-the-art models for this type of data. While a detailed comparison of random forests and hierarchical Bayesian models is beyond the scope of the present paper, we believe these result are encouraging and underscore the potential of tree ensemble models, such as random forests, to be used for profit/revenue prediction and for making marketing decisions. 

\begin{table}
\small
\centering
\begin{tabular}{lrrrrrrrrrrr} \toprule
& \multicolumn{11}{c}{Brand $R^2$} \\
  Model & 1 & 2 & 3 & 4 & 5 & 6 & 7 & 8 & 9 & 10 & 11 \\  \midrule 
  HB-SemiLog & 0.25 & 0.75 & 0.25 & 0.08 & 0.33 & 0.53 & 0.40 & 0.32 & 0.15 & -0.62 & 0.55 \\[0.1em]
  HB-SemiLog-Own & 0.39 & 0.81 & 0.39 & 0.31 & 0.38 & 0.63 & 0.53 & 0.44 & 0.30 & 0.32 & 0.66 \\[0.5em] 
  HB-LogLog & 0.32 & 0.76 & 0.30 & 0.13 & 0.36 & 0.55 & 0.49 & 0.34 & 0.31 & -0.30 & 0.58 \\[0.1em] 
  HB-LogLog-Own & 0.44 & 0.81 & 0.43 & 0.32 & 0.42 & 0.64 & 0.59 & 0.44 & 0.43 & 0.42 & 0.67 \\[0.5em] 
  RF, $T = 20$ & 0.72 & 0.84 & 0.62 & 0.84 & 0.82 & 0.71 & 0.75 & 0.65 & 0.65 & 0.77 & 0.74 \\ 
  RF, $T = 50$ & 0.74 & 0.84 & 0.63 & 0.84 & 0.83 & 0.72 & 0.76 & 0.66 & 0.67 & 0.79 & 0.75 \\ 
  RF, $T = 100$ & 0.75 & 0.85 & 0.63 & 0.84 & 0.83 & 0.72 & 0.76 & 0.67 & 0.67 & 0.79 & 0.76 \\ 
  RF, $T = 500$ & \bfseries 0.75 & \bfseries 0.85 & \bfseries 0.63 & \bfseries 0.84 & \bfseries 0.83 & \bfseries 0.72 & 0.76 & \bfseries 0.67 & \bfseries 0.67 & \bfseries 0.79 & \bfseries 0.76 \\[0.25em] 
  RF-Own, $T = 20$ & 0.69 & 0.83 & 0.61 & 0.81 & 0.76 & 0.71 & 0.75 & 0.64 & 0.64 & 0.73 & 0.73 \\ 
  RF-Own, $T = 50$ & 0.70 & 0.84 & 0.62 & 0.81 & 0.77 & 0.72 & 0.76 & 0.66 & 0.66 & 0.74 & 0.74 \\ 
  RF-Own, $T = 100$ & 0.71 & 0.84 & 0.63 & 0.82 & 0.78 & 0.72 & 0.76 & 0.67 & 0.66 & 0.75 & 0.75 \\ 
  RF-Own, $T = 500$ & 0.72 & 0.84 & 0.63 & 0.82 & 0.78 & 0.72 & \bfseries 0.77 & 0.67 & 0.66 & 0.75 & 0.75 \\  \bottomrule
   \end{tabular}
   
    \caption{Comparison of out-of-sample profit prediction $R^2$ for the different models and brands, averaged over ten random splits of the data, for the second prediction task ($c_i = 0.9p_{i,\min}$ for all brands $i$). The best $R^2$ value for each brand is indicated in bold. \label{table:orange_juice_task2_Rsq} }
\end{table}

\subsection{Optimization results}
\label{subsec:pricing_optimization}

We now turn our attention to optimization. Using the complete data set, we estimate the HB-LogLog-Own model, HB-SemiLog-Own model and RF model with 50 trees per brand. We fix the price of each brand $c_i$ to $0.9 \times p_{i,\min}$. We restrict the price vector $\pb$ to a set $\mathcal{P}$. We will specifically consider the following choice of $\Pcal$:  
\begin{equation}
\Pcal = \prod_{i=1}^{11}  \Pcal_i
\end{equation}
where
\begin{equation}
\Pcal_i = \left\{ \delta \cdot \left \lceil \frac{ p_{i,q25}}{ \delta} \right \rceil,  \delta \cdot  \left \lceil \frac{p_{i, q25}}{ \delta} \right \rceil +  \delta, \delta \cdot  \left \lceil \frac{p_{i, q25}}{ \delta} \right \rceil +  2\delta,  \ \dots, \  \delta \cdot \left \lfloor \frac{p_{i,q75}}{ \delta} \right \rfloor \right\}.
\end{equation}
In the above definition of $\Pcal_i$, $\delta$ is a discretization parameter -- for example, $\delta = 0.05$ indicates that prices go up in increments of \$0.05 -- and $p_{i,q25}$ and $p_{i, q75}$ are respectively the 25th percentile and 75th percentile of the prices observed for brand $i$ in the whole data set. In words, the above expression restricts brand $i$'s price to go up in increments of $\delta$, starting at the smallest multiple of $\delta$ above $p_{i,q25}$ and ending at the largest multiple of $\delta$ below $p_{i,q75}$; for example, if $p_{i, q25} = 1.44$ and $p_{i,q75} = 1.78$, then brand $i$'s prices would be restricted to $\{1.45, 1.50, 1.55, 1.60, 1.65, 1.70, 1.75\}$. We will consider values of the discretization parameter $\delta \in \{0.05, 0.10, 0.20\}$. The reason we restrict each product's price to lie between the 25th and 75th percentile is to ensure the prices are sufficiently ``inside'' the range of previously offered prices, and to prevent them from taking extreme values. For the purpose of this experiment, we do not consider optimization of $\db$ and $\fb$, so we fix each $d_i$ and $f_i$ to zero. 

For HB-LogLog-Own, our profit optimization problem for store $s$ can thus be formulated as 
\begin{equation}
\underset{ \pb \in \mathcal{P} }{\text{maximize}} \quad \sum_{i=1}^{11}  (p_i - c_i) \cdot \left[ \frac{1}{J} \sum_{j=1}^{J} \exp \left( \beta^{(j)}_{s,i,0} + (\betab^{(j)}_{s,i,\pb})^T \log(\pb) + \tau_{s,i}^{(j)} / 2  \right) \right]. \label{eq:orange_juice_HBLogLog_opt}
\end{equation}

For HB-SemiLog-Own, our profit optimization problem for store $s$ is 
\begin{equation}
\underset{ \pb \in \mathcal{P} }{\text{maximize}} \quad \sum_{i=1}^{11}  (p_i - c_i) \cdot \left[ \frac{1}{J} \sum_{j=1}^{J} \exp \left( \beta^{(j)}_{s,i,0} + (\betab^{(j)}_{s,i,\pb})^T \pb + \tau_{s,i}^{(j)} / 2  \right) \right]. \label{eq:orange_juice_HBLogLinear_opt}
\end{equation}
We solve both \problemeqref{eq:orange_juice_HBLogLog_opt} and \eqref{eq:orange_juice_HBLogLinear_opt} using local search from ten different randomly chosen starting points. We remark here that local search is not guaranteed to obtain provably optimal solutions. However, for the coarsest price increment $\delta = 0.20$, it is possible to solve both \eqref{eq:orange_juice_HBLogLog_opt} and \eqref{eq:orange_juice_HBLogLinear_opt} by complete enumeration; in doing so, we find that the solutions returned by local search are provably optimal for both of these problems and for all 83 stores. Thus, it is reasonable to expect that at finer discretization levels ($\delta = 0.05$ or 0.10), the solutions should remain close to optimal.

For RF, our profit optimization problem for store $s$ is 
\begin{equation}
\underset{ \pb \in \mathcal{P} }{ \text{maximize} } \quad \sum_{i=1}^{11} F_i( \pb, \zerob, \zerob, \zb_s), \label{eq:orange_juice_RF_opt}
\end{equation}
where $F_i$ is the random forest prediction function for the profit from brand $i$, $\zerob$ is a vector of zeros of the appropriate dimension (for the two arguments above, both are of length 11), and $\zb_s$ is the vector of store-level covariates of store $s$. Regarding \problemeqref{eq:orange_juice_RF_opt}, we observe that:
\begin{enumerate}
\item The random forest defining each $F_i$ may contain not only splits on $\pb$, but also splits on $\db$, $\fb$ and $\zb$. However, $\db$, $\fb$ and $\zb$ are fixed to $\zerob$, $\zerob$ and $\zb_s$, and are not decision variables. 
\item While the actual split points on $p_1, \dots, p_{11}$ could take any value, each $p_i$ in the optimization problem~\eqref{eq:orange_juice_RF_opt} is restricted to values in $\Pcal_i$. 
\end{enumerate}
These two observations are valuable because we can use them to simplify the tree model. In particular, the two observations can be used to identify splits/leaves that are unreachable in each tree, and to thus remove redundant splits. For example, if $z_1 = 10.5$ in $\zb_s$ and a split in a tree has the query ``Is $z_1 \leq 10.1$?'', then all of the splits and leaves to the left of that split (the ``yes'' branch) can be removed because $z_1$ does not satisfy the query. As another example, suppose $\delta = 0.05$ in the definition of $\Pcal$, and we take the ``yes'' branch for a split with query ``Is $p_1 \leq 1.69$?''; if we then encounter the query ``Is $p_1 \leq 1.66$?'', we know that we cannot reach the node on the ``no'' branch because this would imply $1.66 < p_1 \leq 1.69$, but $p_1$ is restricted to multiples of $\delta = 0.05$. 

Once we have identified all such splits/leaves that cannot be reached, we can remove them and ``collapse'' the tree to obtain a much simpler, store-specific tree that is only in terms of $\pb$. Using these collapsed trees, we formulate problem~\eqref{eq:orange_juice_RF_opt} using formulation~\eqref{prob:TEOMIO} and solve it using the split generation approach. Due to the large number of MIO problems that need to be solved (83 stores by 3 discretization levels), we warm start each MIO by solving \problemeqref{eq:orange_juice_RF_opt} using local search from ten random starting points, and supplying the best of those solutions as the initial solution to Gurobi. In addition, we also used slightly modified parameters for Gurobi, which we report in Table~\ref{table:gurobi_params}.

\begin{table}[ht]
\centering
\begin{tabular}{lr} \toprule
Parameter & Value \\ \midrule
\texttt{Heuristics} & 0 \\
\texttt{Cuts} & 0 \\
\texttt{VarBranch} & 1 \\
\texttt{InfUnbdInfo} & 1 \\
\texttt{PrePasses} & 1 \\ \bottomrule
\end{tabular}
\caption{Gurobi parameters for RF profit optimization problem in Section~\ref{subsec:pricing_optimization}. \label{table:gurobi_params}}
\end{table}

Table~\ref{table:orange_juice_times} displays the average and maximum time to solve problems~\eqref{eq:orange_juice_HBLogLog_opt}, \eqref{eq:orange_juice_HBLogLinear_opt} and \eqref{eq:orange_juice_RF_opt}, for different values of $\delta$, where the average and maximum are taken over the 83 stores in the chain. (For RF, the time includes both the local search warm start time as well as the split generation MIO time.) From this table we can see that the two HB optimization problems are solved quite quickly, requiring no more than 5 seconds across all levels of discretization and all 83 stores. The RF problem~\eqref{eq:orange_juice_RF_opt} was solved in about 11.4 seconds on average per store for the finest discretization level $\delta = 0.05$ and in all instances within 30 seconds. Although the RF problem is not solved quite as fast as the two HB-based problems, the increase is modest, and reasonable given the fact that the MIO approach provides a provably optimal solution, whereas the local search solution method needed for the two HB problems does not. %

\begin{table}[ht]
\centering

\begin{tabular}{lcccccc} \toprule
Price & \multicolumn{2}{c}{RF} & \multicolumn{2}{c}{HB-SemiLog-Own} & \multicolumn{2}{c}{HB-LogLog-Own} \\
Increment $\delta$ & Avg. & Max. & Avg. & Max. & Avg. & Max. \\  \midrule
  0.05 & 11.4 & 26.3 & 1.9 & 2.5 & 2.0 & 2.6 \\ 
  0.10 & 7.3 & 16.5 & 1.2 & 1.6 & 1.2 & 4.9 \\ 
  0.20 & 1.8 & 2.3 & 0.7 & 0.9 & 0.7 & 1.0 \\ \bottomrule
   \end{tabular}

\caption{ Average and maximum computation times (in seconds) of store-level optimal prices for the different models. \label{table:orange_juice_times} }
\end{table}

Aside from computation times, it is also interesting to compare the quality of the price prescriptions themselves. The challenge with comparing the prices produced by each model is that the ground truth model that maps prices to expected profit is unknown to us: we cannot simply plug the prices into a ``true'' model and see which method performs best. To deal with this, we take the following approach. Let $\phi(\pb; m, s)$ be the predicted profit of the price vector $\pb$ at store $s$ according to model $m$, where $m \in \{ \text{RF}, \text{HB-LogLog-Own}, \text{HB-SemiLog-Own} \}$. Let $\pb_{m,s}$ be the optimal price vector if we assume model $m$ for store $s$. Lastly, let $\vartheta(m_1, m_2)$ denote the relative performance of the optimal price vector derived from model $m_1$ when evaluated according to $m_2$:
\begin{equation*}
\vartheta(m_1, m_2) = \frac{1}{83} \sum_{s=1}^{83} \frac{ \phi(\pb_{m_1,s}; m_2, s) }{ \phi(\pb_{m_2,s}; m_2, s) }. 
\end{equation*}
The metric $\vartheta(m_1, m_2)$ measures how well our prices do if we assume reality behaves according to $m_1$, but reality actually behaves according to $m_2$; this performance is measured relative to the best possible profit under $m_2$ (obtained by using the price vector $\pb_{m_2,s}$), and is averaged over all of the stores. For example, $\vartheta( \text{RF}, \text{HB-LogLog-Own})$ would mean we plug the random forest price vector into the HB-LogLog-Own model, and measure how much of the maximum profit under HB-LogLog-Own it recovers. Note that by definition, we always have $\vartheta(m, m) = 100\%$. In essence, the idea of this approach is to be agnostic to the model, and to measure how much is lost if we are wrong about what the underlying model is. For example, suppose that there were only two models, $m_1$ and $m_2$. If it turned out that $\vartheta(m_1, m_2) = 90\%$ but $\vartheta(m_2, m_1) = 60\%$, then assuming $m_1$ is ``safer'' than assuming $m_2$: assuming $m_1$ and being wrong (\ie, reality follows  $m_2$) means only missing out on 10\% of the highest possible profit, whereas assuming $m_2$ and being wrong means missing out on 40\% of the highest possible profit. 

Table~\ref{table:orange_juice_MMC} shows the results of this multi-model comparison. From this table, we can see that the HB-SemiLog-Own price vectors perform very well under the HB-LogLog-Own model, and vice versa. Compared to RF, we can see that the RF price vectors are about as robust as the HB price vectors. For example, with $\delta = 0.20$, assuming the ground truth is RF but reality behaves according to HB-SemiLog-Own means capturing only 88.1\% of the best possible profit under HB-SemiLog-Own (a loss of 11.9\%). Conversely, assuming the ground truth is HB-SemiLog-Own but reality behaving according to RF means capturing only 87.6\% of the best possible profit under RF (a loss of 12.4\%). Given the inherent uncertainty in how reality will behave, we can see that optimizing against the RF model is no more risky than optimizing under either of the HB models.

\begin{table}
\begin{tabular}{clccc} \toprule
Price & & \multicolumn{3}{c}{ $\vartheta(m_1, m_2)$ } \\
Increment $\delta$ &  & $m_2 = $ HB-SemiLog-Own & $m_2 = $ HB-LogLog-Own & $m_2 = $ RF \\ \midrule
0.05 & $m_1 = $ HB-SemiLog-Own & 100.0 & 99.5 & 80.6 \\
 & $m_1 = $ HB-LogLog-Own & 99.5 & 100.0 & 83.8 \\
 & $m_1 = $ RF & 83.3 & 82.9 & 100.0 \\ \midrule
0.1 & $m_1 = $ HB-SemiLog-Own & 100.0 & 99.5 & 83.1 \\
 & $m_1 = $ HB-LogLog-Own & 99.5 & 100.0 & 85.8 \\
 & $m_1 = $ RF & 84.0 & 83.5 & 100.0 \\ \midrule
0.2 & $m_1 = $ HB-SemiLog-Own & 100.0 & 99.6 & 87.6 \\
 & $m_1 = $ HB-LogLog-Own & 99.7 & 100.0 & 89.5 \\
 & $m_1 = $ RF & 88.1 & 88.0 & 100.0 \\ \bottomrule
 \end{tabular}
\caption{Multi-model comparison for optimal prices under the three different models. \label{table:orange_juice_MMC}}
 \end{table}

\begin{figure}
\begin{center}
\includegraphics[width=\textwidth]{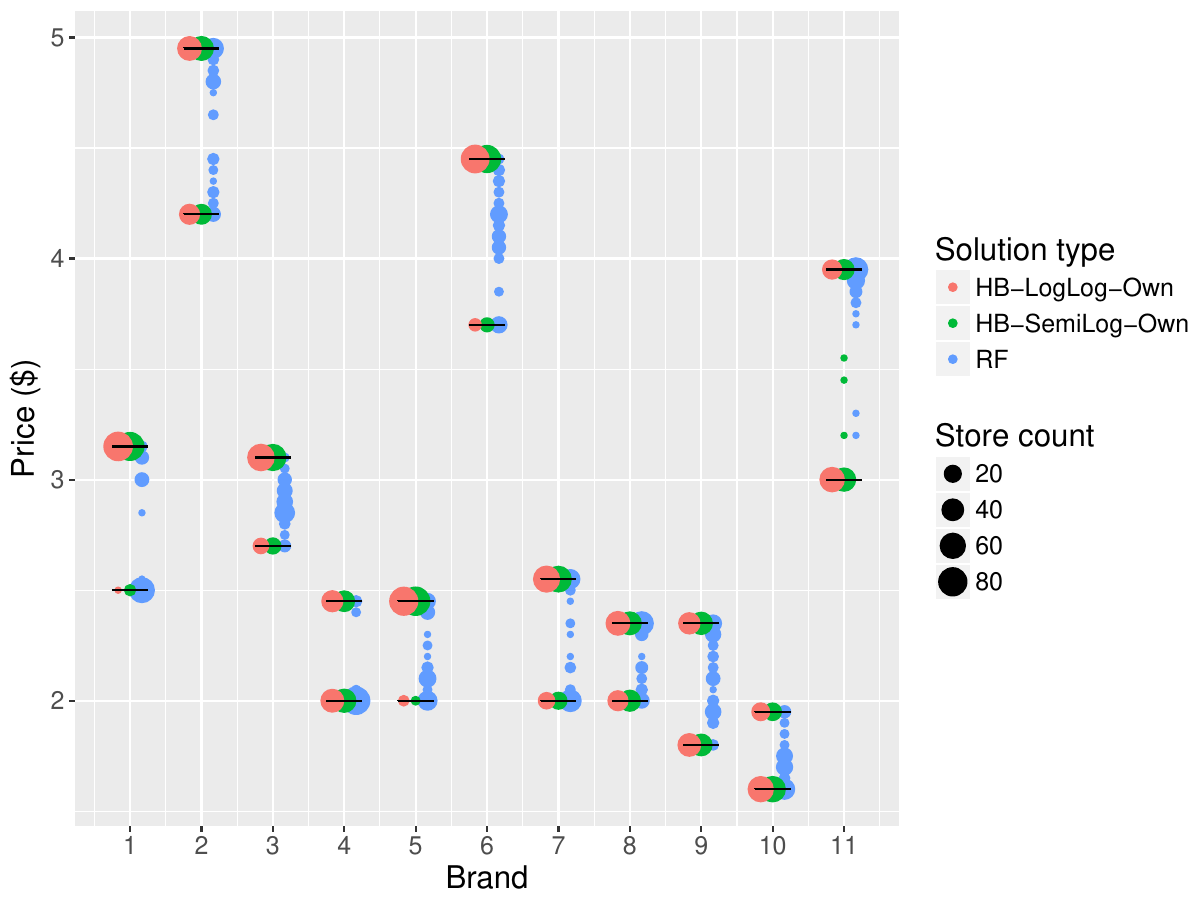}
\end{center}
\caption{Distribution of prices under RF, HB-LogLog-Own and HB-SemiLog-Own. The lowest and highest allowable price for each brand are indicated by black horizontal bars. The size of each circle indicates the number of stores for which the corresponding price was prescribed by the corresponding solution. \label{plot:orange_juice_priceDist}}
\end{figure}

Finally, it is interesting to compare the prices themselves. Figure~\ref{plot:orange_juice_priceDist} shows the optimal prices under the three different models, over the 83 different stores, for each of the eleven brands; the solutions correspond to the increment $\delta = 0.05$. From this plot, we can see that in the overwhelming majority of cases, the optimal price vector $\pb$ under HB-LogLog-Own and HB-SemiLog-Own involves setting each brand's price to the lowest or highest allowable price. In contrast, the price vectors obtained from RF are more interesting: for many stores, there are many brands whose prices are not set to the lowest or highest allowable price. Instead, the RF solutions cover a wider range of values within the price intervals of many brands, and are not as extreme as the prices found for HB-LogLog-Own and HB-SemiLog-Own. The extreme nature of prices under log-log and semi-log aggregate demand models has been highlighted previously in the marketing community \citep[see, e.g.,][]{reibstein1984optimal,anderson2001structural}. %

Overall, from a prescriptive standpoint, we see that the RF profit optimization problem can be solved quickly to provable optimality (Table~\ref{table:orange_juice_times}), leading to prices that are not extreme (Table~\ref{plot:orange_juice_priceDist}) and still achieve high profits under HB models (Table~\ref{table:orange_juice_MMC}). We believe that these results, combined with the strong out-of-sample predictive accuracy shown in Section~\ref{subsec:pricing_predictive_accuracy}, underscore the potential benefit of random forests and tree ensemble models in customized pricing, as well as the value of our optimization methodology in transforming such tree ensemble models into pricing decisions.

\end{document}